\documentclass[11pt]{amsart}
\usepackage[latin1]{inputenc}
\usepackage{amssymb}
\usepackage[colorlinks, bookmarks=true]{hyperref}
\usepackage{amsfonts,amscd,amsmath,amssymb,enumerate,cleveref}

\usepackage[T1]{fontenc}
\usepackage{mathptmx}

\newtheorem{theorem}{Theorem}[section]
\newtheorem{proposition}[theorem]{Proposition}
\newtheorem{lemma}[theorem]{Lemma}
\newtheorem{corollary}[theorem]{Corollary}
\theoremstyle{definition}
\newtheorem{remark}[theorem]{Remark}

\numberwithin{equation}{section}

\newcommand{\unit}{\mathbf{1}}
\newcommand{\element}{\mathcal{E}\ell}

\begin{document}
	
\title[Preservers of Operator Commutativity]{Preservers of Operator Commutativity}

\author[G.M. Escolano]{Gerardo M. Escolano}

\address{Departamento de An{\'a}lisis Matem{\'a}tico, Facultad de Ciencias, Universidad de Granada \hyphenation{Gra-nada}, 18071 Granada, Spain.}
\email{gemares@ugr.es}

\author[A.M. Peralta]{Antonio M. Peralta}

\address{Departamento de An{\'a}lisis Matem{\'a}tico, Facultad de Ciencias, Universidad de Granada, 18071 Granada, Spain.
Instituto de Matem{\'a}ticas de la Universidad de Granada (IMAG).}
\email{aperalta@ugr.es}

\author[A.R. Villena]{Armando R. Villena}

\address{ Departamento de An{\'a}lisis Matem{\'a}tico, Facultad de Ciencias, Universidad de Granada, 18071 Granada, Spain.
Instituto de Matem{\'a}ticas de la Universidad de Granada (IMAG).}
\email{avillena@ugr.es}

\subjclass[2010]{Primary 16W10; 47B49; 17C10; 17C50;  16R50.}

\keywords{Preservers; Jordan algebras; JBW$^*$-algebras; associating traces; associating linear maps;  functional identities; preservers of operator commutativity.}

\begin{abstract} Let $\mathfrak{M}$ and $\mathfrak{J}$ be JBW$^*$-algebras admitting no central summands of type $I_1$ and $I_2,$ and let  $\Phi: \mathfrak{M} \rightarrow \mathfrak{J}$ be a linear bijection preserving operator commutativity in both directions, that is, $$[x,\mathfrak{M},y] = 0  \Leftrightarrow [\Phi(x),\mathfrak{J},\Phi(y)] = 0,$$ for all $x,y\in \mathfrak{M}$, where the associator of three elements $a,b,c$ in $\mathfrak{M}$ is defined by $[a,b,c]:=(a\circ b)\circ c - (c\circ b)\circ a$. We prove that under these conditions there exist a unique invertible central element $z_0$ in $\mathfrak{J}$, a unique Jordan isomorphism  $J: \mathfrak{M} \rightarrow \mathfrak{J}$, and a unique linear mapping  $\beta$ from $\mathfrak{M}$ to the centre of $\mathfrak{J}$ satisfying  $$ \Phi(x) = z_0 \circ J(x) + \beta(x), $$ for all $x\in \mathfrak{M}.$ Furthermore, if $\Phi$ is a symmetric mapping (i.e., $\Phi (x^*) = \Phi (x)^*$ for all $x\in \mathfrak{M}$), the element $z_0$ is self-adjoint, $J$ is a Jordan $^*$-isomorphism, and $\beta$ is a symmetric mapping too.
	
In case that $\mathfrak{J}$ is a JBW$^*$-algebra admitting no central summands of type $I_1$, we also address the problem of describing the form of all symmetric bilinear mappings $B : \mathfrak{J}\times \mathfrak{J}\to \mathfrak{J}$ whose trace is associating (i.e., $[B(a,a),b,a] = 0,$ for all $a, b \in \mathfrak{J})$ providing a complete solution to it. We also determine the form of all associating linear maps on $\mathfrak{J}$.	
\end{abstract}

\maketitle

\section{Introduction}\label{sec: intro}

In the mathematical formulation of quantum mechanics, (bounded) observables are given by self-adjoint operators on a complex Hilbert space $H$ (cf., for example, \cite[page 75]{MackeyBook63}, \cite{GriffithsBook1995}). Some pairs of quantum observables may not be simultaneously measurable, a property referred to as complementarity, what is mathematically expressed by non-commutativity of their corresponding self-adjoint operators. Perhaps for this reason, linear maps preserving commutativity are among the most extensively studied operators in the setting of preserver problems on associative algebras (see, for example, \cite{BannMathi97,Bres1992,Bres1993,BreMei1993,ChoiJafRadjavi87,Miers88,OmlaRadjSem2001,Petek97}).  The contribution by M. Bre\v{s}ar and C.R. Miers describing the linear bijections between von Neumann algebras preserving commutativity in both senses constitute one of the most influencing achievements in this line. The result shows that if $M$ and $N$ are von Neumann algebras with no central summands of type $I_1$ or $I_2$, and $\Theta : M\to N$ is a  bijective additive map which preserves commutativity in both directions then there exist a central invertible element $c\in N$, a Jordan isomorphism $J:  M \to  N$, and an additive map $f$ from $M$ to the centre of $N$, such that $$\Theta (x) = c J(x) + f (x), \hbox{ for all } x\in M.$$ It is additionally proved that if $\Theta$ is a symmetric mapping (i.e., $\Theta (x^*) = \Theta (x)^*$, $x\in M$), then $J$ is a Jordan $^*$-isomorphism.\smallskip
  
In words of Alfsen and Shultz \cite[Introduction of chapter 6]{AlfsenShultz2003} ``\emph{When a C$^*$-algebra or a von Neumann algebra is used as an algebraic model of quantum mechanics, then it is only the self-adjoint part of the algebra that represents observables. However, the self-adjoint part of such an algebra is not closed under the given associative product, but only under the Jordan product $a\circ b =\frac12 (a b + b a)$. Therefore it has been proposed to model quantum mechanics on Jordan algebras rather than associative algebras \cite{JordanvNeumannWigner34, vNeumann36}.}'' Jordan norm-closed self-adjoint subalgebras of C$^*$-algebras are called special \emph{JB$^*$-algebras} or \emph{JC$^*$-algebras}. There are examples of Jordan algebras, called \emph{exceptional Jordan algebras}, which cannot be embedded as a Jordan subalgebras of an associative algebra, an example is given by the Jordan algebra  $H_3(\mathbb{O}),$ of all Hermitian $3\times 3$ matrices with entries in the complex octonions (cf. \cite[Corollary 2.8.5]{HOS}).\smallskip

A real or complex Jordan-Banach  algebra is a Banach space $\mathfrak{J}$ equipped with a bilinear mapping $(a,b) \rightarrow a \circ b$ (called the Jordan product) satisfying the following axioms: 
\begin{itemize}
	\item[(J1)] $\| a \circ b\| \leq \|a\|\|b\|$, for all $a,b \in \mathfrak{J}$; 
	\item[(J2)] $a\circ b = b\circ a$, for all $a, b \in \mathfrak{J} \ \ $  (commutativity);
	\item[(J3)] $(a^2 \circ b)\circ a = (a \circ b)\circ a^2$, for all $a,b \in \mathfrak{J} \ \ $ (Jordan identity).
\end{itemize}

There are two closely related types of Jordan-Banach  algebras, JB-algebras and JB$^*$-algebras. A \emph{JB-algebra} is a real Jordan-Banach  algebra $\mathfrak{J}$ in which the norm satisfies the following two additional conditions: 
\begin{itemize}
	\item[(JB1)] $\|a^2\| = \|a\|^2$  for all $a \in \mathfrak{J}$; 
	\item[(JB2)] $\|a^2\| \leq \|a^2 + b^2\|$ for all $a,b \in \mathfrak{J}$.
\end{itemize}
If $\mathfrak{J}$ is unital, it is clear from (JB1) that $\|\textbf{1}\| = 1$.\smallskip

The Jordan-model analogue to C$^*$-algebras, introduced by I. Kaplansky in 1976, is known as JB$^*$-algebra. A complex Jordan-Banach  algebra $\mathfrak{J}$ equipped with an involution $^*$ is said to be a \emph{JB$^*$-algebra} if the following condition holds:
\begin{itemize}
\item[(JB$^*$1)] $\|a\|^3 = \| U_a(a^*)\|$ for all $a \in \mathfrak{J}$.
\end{itemize}
If $\mathfrak{J}$ is unital, it clearly follows that $\textbf{1}^* = \textbf{1}$. It is known that the involutions of each JB$^*$-algebra $\mathfrak{J}$ satisfies $\|a^*\| = \|a\|$ for all $a \in \mathfrak{J}$ (cf. \cite[Lemma 4]{Youngson78}).\smallskip

These two Jordan mathematical models are intrinsically connected in the following way: the set, $\mathfrak{J}_{sa},$ of all \emph{self-adjoint} elements in a JB$^*$-algebra $\mathfrak{J}$, i.e. $\mathfrak{J}_{sa} : = \{a \in \mathfrak{J} \ : \, a^* = a\}$, is a JB-algebra \cite[see Proposition 3.8.2]{HOS}. Conversely, by a deep result due to J.D.M. Wright (cf. \cite{Wright1977}), each JB-algebra corresponds to the self-adjoint part of a (unique) JB$^*$-algebra.\smallskip

A \emph{JBW$^*$-algebra} (resp., a \emph{JBW-algebra}) is a JB$^*$-algebra (resp., a JB-algebra) which is also a dual Banach space. Thus, JBW$^*$-algebras can be considered as the Jordan analogue of von Neumann algebras. Additional details and some structure results are revisited along subsections \ref{subsec:background}, \ref{subsec:operator commutativity}, \ref{subsec: structure results}, and  \ref{subsec: spin factors}.\smallskip

The Jordan product of each JB$^*$-algebra $\mathfrak{J}$ is commutative by definition. This property makes indistinguishable the left and right multiplication operators. For each $a \in \mathfrak{J}$ we denote by $M_a$ the Jordan multiplication operator by the element $a$, that is, $M_a : \mathfrak{J} \rightarrow \mathfrak{J}, \, M_a (b) = a\circ b$. We say that elements $a$ and $b$ in a Jordan algebra $\mathfrak{J}$ \emph{operator commute} if the operators $M_a, M_b$ commute, that is, $(a\circ c)\circ b = a \circ (c\circ b)$ for every $c \in \mathfrak{J}$, equivalently, the associator $[a,c,b] = (a\circ c)\circ b - a \circ (c\circ b)$ is zero for all $c\in \mathfrak{J}$. D.M. Topping showed in \cite[Proposition 1]{Topping65} that if $\mathfrak{J}$ is a JC$^*$-algebra, regarded as a JB$^*$-subalgebra of a C$^*$-algebra $A$, two self-adjoint elements $a$ and $b$ in $\mathfrak{J}$ operator commute if, and only if, they commute with respect to the associate product of the C$^*$-algebra $A$. Subsection \ref{subsec:operator commutativity} is devoted to revisit the basic background on operator commutativity. We also add some new result proving that if $a$ and $b$ are elements in a JB$^*$-algebra $\mathfrak{J}$ such that the JB$^*$-subalgebra of $\mathfrak{J}$ that they generate is a JB$^*$-subalgebra of some C$^*$-algebra $A$, then $a$ and $b$ commute in the usual sense as elements of $A$ as soon as they operator commute in $\mathfrak{J}$. Consequently, two elements $a,b$ in a JC$^*$-algebra $\mathfrak{J}$ acting on a C$^*$-algebra $A$ operator commute in $\mathfrak{J}$ if, and only if, they commute in $A$ (see \Cref{P Topping for non-self adjoint}).\smallskip

M. Bre\v{s}ar, D. Eremita, and the third author of this note studied linear bijections preserving operator commutativity between certain classes of Jordan algebras in \cite{BreEreVil}. The result, which is a bit more general, reads as follows: Let $\mathfrak{M}$ and $\mathfrak{J}$ be non-quadratic, unital, centrally closed, prime non-degenerate Jordan algebras over a field $\mathbb{K}$ of characteristic different from $2$, $3$, and $5$. Then for every bijective linear mapping $\Phi: \mathfrak{M}\to \mathfrak{J}$  such that $[\Phi(a^2), \mathfrak{J}, \Phi(a)]=0,$ for all $a \in\mathfrak{J}$, there exist a non-zero $\alpha\in \mathbb{K},$ a Jordan isomorphism $J: \mathfrak{M}\to \mathfrak{J},$ and a linear mapping $\beta: \mathfrak{M}\to \mathbb{K}$ satisfying  $$\Phi (x) = \alpha J(x) + \beta (x) \unit,$$ for all $x\in \mathfrak{M}$ (cf. \cite[Theorem 5.3]{BreEreVil}).\smallskip

The main purpose of this study is to extend the just quoted result by Bre\v{s}ar, Eremita and Villena to linear bijections preserving operator commutativity in both directions between JBW$^*$-algebras admitting no central summands of type $I_1$ and $I_2$, \`a la Bre\v{s}ar--Meiers. Our study culminates in \Cref{t linear bijections preserving oper commut} where we establish that if $\mathfrak{M}$ and $\mathfrak{J}$ are JBW$^*$-algebras with no central summands of type $I_1$ and $I_2,$ and  $\Phi: \mathfrak{M} \rightarrow \mathfrak{J}$ is a linear bijection preserving operator commutativity in both directions, that is, $$[x,\mathfrak{M},y] = 0  \Leftrightarrow [\Phi(x),\mathfrak{J},\Phi(y)] = 0,$$ for all $x,y\in \mathfrak{M}$, then there exist a unique invertible central element $z_0$ in $\mathfrak{J}$, a unique Jordan isomorphism  $J: \mathfrak{M} \rightarrow \mathfrak{J}$, and a unique linear mapping  $\beta$ from $\mathfrak{M}$ to the centre of $\mathfrak{J}$ satisfying  
$$ \Phi(x) = z_0 \circ J(x) + \beta(x), $$ for all $x\in \mathfrak{M}.$ Furthermore, if $\Phi$ is a symmetric mapping (i.e., $\Phi (x^*) = \Phi (x)^*$ for all $x\in \mathfrak{M}$), the element $z_0$ is self-adjoint, $J$ is a Jordan $^*$-isomorphism, and $\beta$ is a symmetric mapping too (cf. \Cref{c linear bijections preserving oper commut symmetric}). \smallskip

In a spin factor $V$, two elements $a,b$ operator commute if and only if $b$ is a linear combination of $a$ and the unit element (see \Cref{r associator-commutator of an element in a spin}). Therefore, any linear mapping on a spin factor preserves operator commutativity. This is essentially the reason to avoid JBW$^*$-algebras containing summands of type $I_2$ in \Cref{t linear bijections preserving oper commut} (see \Cref{r optimality of the hypotheses in the previous theorem}).\smallskip

We also obtain a version for linear preservers of operator commutativity between JBW-algebras in \Cref{t linear bijections preserving oper commut JBW-algebras}, where we prove for every couple of JBW$^*$-algebras admitting no central summands of type $I_1$ and $I_2$, $\mathfrak{M}$ and $\mathfrak{J}$, and every linear bijection $\Phi: \mathfrak{M}_{sa} \rightarrow \mathfrak{J}_{sa}$ preserving operator commutativity in both directions, there exist a unique invertible element $z_0$ in $Z(\mathfrak{J}_{sa})$, a unique Jordan isomorphism $J: \mathfrak{M}_{sa} \rightarrow \mathfrak{J}_{sa}$, and a unique linear mapping  $\beta: \mathfrak{M}_{sa}\to Z(\mathfrak{J}_{sa})$ satisfying  
$$ \Phi(x) = z_0 \circ J(x) + \beta(x), $$ for all $x\in \mathfrak{M}_{sa}.$  The conclusion also holds when $\mathfrak{M}$ and $\mathfrak{J}$ are von Neumann algebras without central summands of type $I_1$ and $I_2$ and ``operator commutativity'' means ``commutativity'' (cf. \Cref{c self-adjoint parts of von Neumann algebras}), a result which improves the known conclusions in this setting.\smallskip

The just presented results have required a time-consuming effort and a great deal of novelty in order to combine tools and techniques from several branches such as C$^*$-algebra theory, Jordan algebra theory and preservers. The arguments are technically complex, but we have tried to make the work self-contained. It is perhaps worth to comment the route map we followed. In words of M. Bresar \cite{BresarReviewZB}, ``It has turned out that for rather large classes of rings all these (and some similar) problems can be solved using a unified approach based on a characterization of commuting traces of biadditive maps.'' In the Jordan setting, ``commuting traces'' are replaced by ``associating traces''.\smallskip

A mapping $q$ on a  linear space $X$ is said to be a \emph{trace} if there exists a symmetric bilinear mapping $B: X \times X \to X$ such that $q(x) = B(x, x)$ for all $x \in X$. Let $F: \mathfrak{J}\to \mathfrak{J}$ be a mapping on a Jordan algebra. We shall say that $F$ is \emph{associating} if for every $x\in \mathfrak{J}$, $F(x)$ and $x$ operator commute as elements in $\mathfrak{J}$, that is, $[F(x),\mathfrak{J},x]=0$. In \Cref{t linear associating maps on JBW*-algebras without type I1} we prove that every linear associating mapping $T$ on a  JBW$^*$-algebra $\mathfrak{J}$ admitting no direct summands of type $I_1$, can be expressed in the form: 
\begin{equation*}
	T(x) = \lambda\circ x + \mu(x), \hbox{ for all } x\in \mathfrak{J},
\end{equation*}
where $\lambda \in Z (\mathfrak{J})$, and $\mu: \mathfrak{J} \rightarrow Z (\mathfrak{J})$ is linear. Furthermore, for every associating trace $q(x) = B(x,x)$ on  $\mathfrak{J}$, there exist a unique $\lambda \in Z (\mathfrak{J})$, a unique linear mapping $\mu: \mathfrak{J}\rightarrow Z (\mathfrak{J})$, and a unique bilinear map $\nu : \mathfrak{J}\times \mathfrak{J} \rightarrow Z (\mathfrak{J})$ satisfying 
\begin{equation*}
	B(x,x) = \lambda \circ x^2 + \mu(x) \circ x + \nu(x,x)\circ \mathbf{1}, 
\end{equation*} for all $x \in \mathrm{J}$. If $\mathfrak{J}$ is a JBW$^*$-algebra of type $I_2$ the element $\lambda$ is always zero (cf. \Cref{teo bilinear mapping JBW* algebra}).

\subsection{Background }\label{subsec:background}  We say that a Jordan algebra $\mathfrak{J}$ is \emph{unital} if there exists an element $\mathbf{1} \in \mathfrak{J}$ (called the unit of $\mathfrak{J}$) such that $\mathbf{1}\circ a = a$ for all $a \in \mathfrak{J}$. Given elements $a,c \in \mathfrak{J}$, the symbol $U_{a,b}$ will stand for the linear mapping on $\mathfrak{J}$ defined by $U_{a,c}(b) := (a\circ b)\circ c + (b\circ c)\circ a - (a\circ c)\circ b$ ($b\in \mathfrak{J}$). We simply write $U_a$ for $U_{a,a}$.\smallskip

As we commented in the introduction, a natural example of Jordan algebra is provided by any associative algebra $A$ equipped with the natural Jordan product given by $a\circ b := \frac{1}{2}(ab + ba)$. Any linear subspace of an associative algebra which is closed under the Jordan product is a Jordan algebra. Such Jordan algebras are called \emph{special}. Given an element $a$ in an associative algebra $A$, we shall write $L_a$ (respectively, $R_a$) for the left (respectively, right) multiplication operator by the element $a$, that is, $L_a (x ) = a x $ (respectively, $R_a (x) = x a$). Jordan algebras which cannot be embedded as a Jordan subalgebras of an associative algebra are called \emph{exceptional}. A widely known example of exceptional Jordan algebra is the algebra $H_3(\mathbb{O})$ of all Hermitian $3\times 3$ matrices with entries in the complex octonions (see \cite[Corollary 2.8.5]{HOS} for more information). \smallskip

An element $a$ in a unital Jordan-Banach  algebra $\mathfrak{J}$ is called \emph{invertible} if we can find $b \in \mathfrak{J}$ satisfying $a\circ b = \unit$ and $a^2\circ b = a$. The element $b$ is unique, it is called the inverse of $a$ and denoted by $a^{-1}$ (cf. \cite[3.2.9]{HOS}).\smallskip

Along this note, the symbol ``$*$'' will mainly stand for the involution of a JB$^*$-algebra. For each algebra $A$ equipped with an algebra involution $^*$, we shall write $A_{sa}$ for the set of all symmetric elements for this involution. The reader should be warned that we shall also write $X^*$ for the dual space of a Banach space $X$, and $X_*$ for a predual of $X$ in case that the latter is a dual Banach space.\smallskip

The reader should know that every JBW-algebra is unital (see \cite[Lemma 4.1.7]{HOS}). It is known that a JB$^*$-algebra $\mathfrak{J}$ is a JBW$^*$-algebra if, and only if, $\mathfrak{J}_{sa}$ is a JBW-algebra \cite{Edwards80}. It is also established in the just quoted reference (see also \cite[Lemma 2.2]{BoHamKal2017}) that  the following assertions hold: 
    \begin{enumerate}[\rm(i)] \label{lemma:relJBW}
        \item[\rm(i)] $\mathfrak{J}_{sa}$ is weak$^*$-closed in $\mathfrak{J}$.
        \item[\rm(ii)] The operator $\phi: (\mathfrak{J}_{*})_{sa} \rightarrow (\mathfrak{J}_{sa})_* $ defined by $\phi(\omega) = \omega |_{\mathfrak{J}_{sa}}$ is an onto linear isometry of real Banach spaces, where $\mathfrak{J}_{*sa} = \{\varphi\in \mathfrak{J}_{sa} :\varphi (a^*) =\overline{\varphi(a)},\ \forall a\in \mathfrak{J} \}$ is the self-adjoint part of the predual, $\mathfrak{J}_*$, of $\mathfrak{J}$.
        \item[\rm(iii)] The operator $\psi : \mathfrak{J}_{sa} \times \mathfrak{J}_{sa} \rightarrow \mathfrak{J}$ defined by $\psi(x,y) = x + iy$ is a onto real-linear weak$^*$-to-weak$^*$ homeomorphism.
    \end{enumerate} 
    
 For the basic theory of JB- and JB$^*$-algebras, the reader is referred to the monographs \cite{AlfsenShultz2003, HOS} and \cite{CabGarPalVol1}. \smallskip
 
Let $B(H)$ denote the C$^*$-algebra of all bounded linear operators on a complex Hilbert space $H$. A \emph{JC-algebra} is a JB-algebra that is isometrically isomorphic to a norm closed Jordan subalgebra of $B(H)_{sa}$ \cite[see Proposition 1.35]{AlfsenShultz2003}. There exist JB-algebras which are not JC-algebras, for instance the algebra $H_3(\mathbb{O})$ of all Hermitian $3\times 3$ matrices with entries in the complex octonions $\mathbb{O}$ (see \cite[Corollary 2.8.5]{HOS}). A \emph{JC$^*$-algebra} is a JB$^*$-algebra which materialises as a norm-closed self-adjoint Jordan subalgebra of a C$^*$-algebra, and hence, by the Gelfand-Neimark theorem, a norm-closed self-adjoint Jordan subalgebra of some $B(H)$. \smallskip

A \emph{JW-algebra} is a weak$^*$-closed real Jordan subalgebra of some $B(H)_{sa}$. JW-algebras were first studied by D.M. Topping \cite{Topping65} and E. St{\o}rmer \cite{Stormer66}. A \emph{JW$^*$-algebra} is a JC$^*$-algebra which is also a dual Banach space, or equivalently, a weak$^*$-closed JB$^*$-subalgebra of some von Neumann algebra.

\subsection{Operator commutativity and the centre}\label{subsec:operator commutativity} The product of a Jordan algebra $\mathfrak{J}$ is, by definition, commutative, so every pair of elements in $\mathfrak{J}$ commute if we only employ the ``usual'' sense. However, if we assume that an associative algebra $A$ is equipped with its natural Jordan product $a\circ b = \frac12 ( ab + ba )$, it can be easily checked that if $a,b\in A$ commute with respect to the associative product (i.e., $ a b = b a$), then the corresponding Jordan multiplication operators $M_a$ and $M_b$ commute. This is the motivation for the usual notion of operator commutativity in a Jordan algebra. According to the standard sources, elements $a$ and $b$ in a Jordan algebra $\mathfrak{J}$ are said to \emph{operator commute} if the operators $M_a, M_b$ commute (i.e., $(a\circ c)\circ b = a \circ (c\circ b)$ for every $c \in \mathfrak{J}$). For example, the Jordan identity is equivalent to say that, for each element $a$ in a Jordan algebra $\mathfrak{J}$,   $a^2$ and $a$ operator commute. The \emph{associator} of three elements $a,b,c$ in a Jordan algebra $\mathfrak{J}$, defined by $[a,c,b] := (a\circ c) \circ b - (b\circ c)\circ a$, testes the operator commutativity of $a$ and $b$ since $a$ and $b$ operator commute in $\mathfrak{J}$ if, and only if, $[a,\mathfrak{J},b]=0.$ \smallskip

On an associative algebra $A$ we shall write $[\cdot,\cdot]$ for the usual Lie bracket given by $[a,b] := a b - ba $, $\forall a,b\in A$. \smallskip

The reader should be warned that, in a general Jordan algebra, operator commutativity of $a$ and $b$ is not always a necessary nor sufficient condition to the property that $a$ and $b$ generate a commutative and associative subalgebra  of $\mathfrak{J}$ (cf. \cite[2.5.1 and Example 2.5.2]{HOS}). \smallskip

The \emph{centre} of a Jordan algebra $\mathfrak{J}$ consists of all elements $z\in \mathfrak{J}$ that operator commute with every element of $\mathfrak{J}$. The symbol $Z(\mathfrak{J})$ will stand for the centre of $\mathfrak{J}$, and its elements are called central. The centre of a JB$^*$-algebra $\mathfrak{J}$ is a commutative C$^*$-algebra, and contains the identity of $\mathfrak{J}$ if it exists (see \cite[Proposition 1.52]{AlfsenShultz2003}. The centre of a JBW$^*$-algebra is a commutative von Neumann (see \cite[Proposition 2.36]{AlfsenShultz2003}, \cite{Edwards80}). If $Z(\mathfrak{J})$ consists of scalar multiples of the identity alone, $\mathfrak{J}$ is called a JBW$^*$-factor.\smallskip

We recall that a \emph{(Jordan) factor representation} of a JB$^*$-algebra $\mathfrak{J}$ is a Jordan $^*$-homomorphism $\pi$ from $\mathfrak{J}$ onto a weak$^*$-dense Jordan subalgebra of a JBW$^*$-factor. Observe that every JB$^*$-algebra admits a faithful or separating family of Jordan factor representations (cf. \cite[Corollary 5.7]{AlfShulStor79GelfandNeumark}). 

\begin{remark}\label{Rem: central representation} Let $a$ be an element in a JB$^*$-algebra $\mathfrak{J}$. The following statements are equivalent:  
	\begin{enumerate}[$(a)$]
		\item $a \in Z(\mathfrak{J})$,
		\item $\pi(a) \in Z(\mathfrak{J}_{\pi}),$ for every Jordan factor representation 
		$\pi : \mathfrak{J} \rightarrow \mathfrak{J}_{\pi}$,
		\item $\pi_i(a) \in Z(\mathfrak{J}_{\pi_i}),$ for every representation $\pi_i$ in a faithful family of Jordan factor representations  $\{\pi_i : \mathfrak{J} \rightarrow \mathfrak{J}_{\pi} \}_{i}$. 
	\end{enumerate}
Namely, the implication $(a) \Rightarrow (b)$ follows by just applying the weak$^*$-density of $\pi(\mathfrak{J})$ in $\mathfrak{J}_{\pi}$ and the separate weak$^*$-continuity of the Jordan product of the latter JBW$^*$-algebra. The implication $(b) \Rightarrow (c)$ is clear. To see $(c)\Rightarrow (a)$, we assume the existence of a faithful family of Jordan factor representations  $\{\pi_i : \mathfrak{J} \rightarrow \mathfrak{J}_{\pi_i}\}_i$ satisfying that $\pi_i (a)$ lies in $Z\left(\mathfrak{J}_{\pi_i}\right)$ for all $i$. Define $\displaystyle \pi_0: \mathfrak{J} \rightarrow \stackrel{\infty}{\bigoplus}_{i} \mathfrak{J}_{\pi_i}$ by $\pi_0(a) = (\pi_i(a))_i$, which is clearly a Jordan $^*$-monomorphism. The hypothesis implies that $\pi_0(a) \in  Z(\stackrel{\infty}{\bigoplus}_i \mathfrak{J}_{\pi_i})$. Since $\pi_0[a, \mathfrak{J}, \mathfrak{J}]= [\pi_0(a), \pi_0(\mathfrak{J}),\pi_0(\mathfrak{J})] = 0$, we deduce that $[a, \mathfrak{J}, \mathfrak{J}]= 0$, and thus $a \in Z(\mathfrak{J})$. 
\end{remark}

In connection with the above remark, it is perhaps worth to note that every Jordan $^*$-monomorphism between JB$^*$-algebras is automatically an isometry, while every Jordan $^*$-homomorphisms is automatically continuous and non-expansive (see \cite[Corollaries 1.4 and 1.5]{Wright1977}).\smallskip

If $\mathfrak{J}$ is a JC$^*$-algebra, regarded as a JB$^*$-subalgebra of a C$^*$-algebra $A$, a result originally due to Topping, shows that two self-adjoint elements $a$ and $b$ in $\mathfrak{J}$ operator commute if, and only if, they commute with respect to the associate product of the C$^*$-algebra $A$ (cf.  \cite[Proposition 1]{Topping65}, see also the subsequent rediscoverings in \cite[Proposition 1.49]{AlfsenShultz2003}, \cite[Lemma 5.1]{HO1983} and  \cite[Proposition 3.2]{vanWetering2020}). The study of operator commutativity for self-adjoint elements in a JB$^*$-algebra $\mathfrak{J}$ was revisited by van de Wetering in \cite[Theorem 3.13]{vanWetering2020}, who showed, via Shirshov-Cohn theorem, that the following statements are equivalent for any $a,b\in \mathfrak{J}_{sa}$:
\begin{enumerate}[$(a)$]
	\item $a$ and $b$ operator commute.
	\item $a$ and $b$ generate an associative JB$^*$-algebra.
	\item $a$ and $b$ generate an associative JB$^*$-algebra of mutually operator commuting elements.
	\item $a$ and $a^2$ operator commute with $b$ and $b^2$.
\end{enumerate} In Topping's result, as well as in the just quoted characterization obtained by van de Wetering, the elements are required to be self-adjoint. Furthermore, in Topping's result the elements $a$ and $b$ are assumed in a JC$^*$-algebra and we have a C$^*$-algebra dominating the Jordan algebra and providing an associative product. Next we refine the argument in Topping's original paper to improve the previous results by showing that if two (non-necessarily self-adjoint) elements in a JB$^*$-algebra operator commute, then they commute inside any C$^*$-algebra containing the JB$^*$-subalgebra that they generate. The conclusion is new even for JC$^*$-algebras. 

\begin{proposition}\label{P Topping for non-self adjoint} Let $a,b$ be elements in a JB$^*$-algebra $\mathfrak{J}$, and let $\mathfrak{J}_{a,b}$ denote the JB$^*$-subalgebra of $\mathfrak{J}$ generated by $a$ and $b$. Suppose we can find a C$^*$-algebra $A$ containing $\mathfrak{J}_{a,b}$ as a JB$^*$-subalgebra. Assume that $a$ and $b$ operator commute in $\mathfrak{J}$. Then $a$ and $b$ commute as elements in $A$. In particular, two elements $a,b$ in a JC$^*$-algebra $\mathfrak{J}$ acting on a C$^*$-algebra $A$ operator commute in $\mathfrak{J}$ if, and only if, they commute in $A$.  
\end{proposition}  

\begin{proof} Since $a$ and $b$ operator commute in $\mathfrak{J}$ we have  $[a,\mathfrak{J},b] =0$. Let us denote the product of $A$ by mere juxtaposition.  Suppose $B$ denotes the C$^*$-subalgebra of $A$ generated by $\mathfrak{J}_{a,b}$. By rewriting the hypotheses in terms of the associative product of $B$ we arrive to $$ 0=-4 [a,\mathfrak{J}_{a,b},b] =  [[a,b], \mathfrak{J}_{a,b}].$$ It is well-known that the mapping $x\mapsto [[a,b],x]$ is a (continuous) derivation on $B$. Therefore, given $y,z\in \mathfrak{J}_{a,b}$ we have $$[[a,b],y z ] = [[a,b], y ] z + y [[a,b],z ] = 0,$$ and thus, a basic argument mixing linearity and continuity implies that $[[a,b], B ] = 0.$ We have therefore shown that $[a,b]$ is a central element in the C$^*$-algebra $B$. In particular $[[a,b],a] = 0,$ and hence Kleinecke's theorem \cite{Kleinecke1957} assures that $[a,b]$ is a quasi-nilpotent element, that is, it has zero spectral radius in $B$, which is equivalent to say that $[a,b]=0$ since the latter is a central element. Therefore, $a$ and $b$ commute in $B$ (and in $A$). The rest is clear.
\end{proof}

\subsection{Some structure results of \texorpdfstring{JBW$^*$-}{JBW*-}algebras}\label{subsec: structure results} 

An element $p$ in a JB-algebra is called a \emph{projection} if $p^2= p \circ p = p$. Similarly, an element $p$ in a JBW$^*$-algebra $\mathfrak{J}$ is called a projection if $p^* = p =p\circ p = p$, i.e., $p$ is a projection in $\mathfrak{J}$ if, and only if, it is a projection in  $\mathfrak{J}_{sa}$. Two projections $p, q$ in $\mathfrak{J}$ are called \emph{orthogonal} ($p\perp q$ in short) if $p\circ q = 0$ (see  \cite[Lemma 4.2.2]{HOS} for equivalent conditions). \smallskip

An element $s$ in a unital JB$^*$-algebra $\mathfrak{J}$ is said to be a \emph{symmetry} if $s=s^*$ and $s^2= \mathbf{1}$. It can be proved that for each symmetry $s \in \mathfrak{J}$ the map $U_s$ is a Jordan $^*$- automorphism (cf. \cite[Proposition  2.34]{AlfsenShultz2003}). The group Int$(\mathfrak{J})$ generated by all Jordan  $^*$-automorphisms of the form $U_s$ with $s$ running in the set of all symmetries in $\mathfrak{J}$ is called the group of all \emph{inner automorphisms} of $\mathfrak{J}$, and its elements are called inner automorphisms of $\mathfrak{J}$. \smallskip

Two projections $p,q$ in a JBW-algebra $\mathfrak{J}$ are called \emph{equivalent} ($p \sim q$ in short) if there is an inner automorphism $\alpha$  of $\mathfrak{J}$ such that $q = \alpha(p)$. If $\alpha$ can be written as $\alpha = U_{s_1}U_{s_2}\cdots U_{s_n}$ we write $p \sim_n q$. If $n = 1$ we say $p$ and $q$ are \emph{exchanged by a symmetry}. \smallskip

A projection $p$ in a JBW$^*$-algebra $\mathfrak{J}$ is called \emph{abelian} if the algebra $\mathfrak{J}_p = U_p(\mathfrak{J}):= \{U_p (x) \, : \, x \in \mathfrak{J}\} $ is associative. If $\mathfrak{J}$ is a JW$^*$-algebra (resp. a JW-algebra) it follows from \cite[1.49]{AlfsenShultz2003} that $p$ is abelian if, and only if, $J_p$ consist of mutually operator commuting elements. \smallskip

Given an associative algebra $A$, we denote by $A^{op}$ the \emph{opposite algebra} of $A$, that is, the algebra formed by reversing the order of the product in $A$. If $A$ is C$^*$-algebra, $A^{op}$ is a C$^*$-algebra with respect to the same norm and involution. \smallskip

The next lemma states a property which is probably well known by experts on JB$^*$-algebras. We included it in this note for completeness reasons. We recall first that given a JBW$^*$-algebra $\mathfrak{J}$, there exists a unique von Neumann algebra $W^*(\mathfrak{J})$ together with a normal Jordan  homomorphism $\Psi: \mathfrak{J}\rightarrow W^*(\mathfrak{J})$ satisfying: 
\begin{enumerate}[$(i)$] \item $ \Psi(\mathfrak{J})$ generates $W^*(\mathfrak{J})$ as a von Neumann algebra. 
\item If $\mathcal{N}$ is a von Neumann algebra and $\Phi: \mathfrak{M} \rightarrow \mathcal{N} $ is a normal homomorphism then there is a normal $^*$-homomorphism $\widehat{\Phi}: W^*(\mathfrak{J}) \rightarrow \mathcal{N}$ such that $\widehat{\Phi}\Psi = \Phi$. 
\item There is a $^*$-anti-automorphism $\tau$ of order-2 of $W^{*}(\mathfrak{J})$ such that $\tau (\Psi(x)) = \Psi(x)$ for all $x\in \mathfrak{J}$, i.e. $\Psi(\mathfrak{J}) \subseteq H(W^{*}(\mathfrak{J}),\tau) =\{x\in W^{*}(\mathfrak{J}) : \tau (x) =x\}.$
\end{enumerate} The von Neumann algebra $W^*(\mathfrak{J})$ is called the \emph{universal von Neumann algebra for $\mathfrak{J}$} \cite[Theorem 7.1.9]{HOS}. It is further known (see \cite[Remark 7.2.8]{HOS}) that the Jordan $^*$-homomorphism $\Psi$ is an isometric Jordan $^*$-monomorphism when $\mathfrak{J}$ is special (i.e. it is a JW$^*$-algebra). If $\mathfrak{J}$ is a JW$^*$-algebra without spin part, we actually have $\mathfrak{J} = H(W^{*}(\mathfrak{J}),\tau)$ (cf. \cite[Proposition 7.3.3]{HOS}).

\begin{lemma}\label{l HWtau factor implies W factor} Let $\tau$ be a $^*$-anti-homomorphism of order-2 on a von Neumann algebra $W$. Then the JW$^*$- algebra $H(W,\tau) =\{x\in W : \tau (x) =x \}$ is a factor if, and only if, one of the next statements holds \begin{enumerate}[$(a)$]\item $W$ is a factor; 
\item $H(W,\tau)$ is $^*$-isomorphic to a factor von Neumann algebra $N$ which is a weak$^*$-closed ideal of $W = N\oplus^{\infty} \tau(N)$, and $H(W,\tau)=\{(a,\tau(a)) : a\in N\}$.    
\end{enumerate}
\end{lemma}

\begin{proof} Suppose that $H(W,\tau)$ is a factor. Let $Z$ stand for the centre of $W$. Clearly $\tau (Z) = Z$, and $\tau|_{Z}: Z\to Z$ is an order-2 $^*$-automorphism. Lemma 7.3.4 in \cite{HOS} assures the existence of two projections $q,p\in Z$ such that $p+q+\tau(q)= \unit$, every subprojection of $p$ in $Z$ is $\tau$-invariant. Observe that $p,q,\tau(q)$ are mutually orthogonal. Since the projections $p,q+\tau(q)$ lie in $H(W,\tau)$ and the latter is a factor, we must have $p=0$ or $q =0$.\smallskip

Suppose first that $p=0$. It is clear that $N = q W$ is a weak$^*$-closed ideal of $W,$ $W = N\oplus^{\infty} \tau(N)$, and $H(W,\tau)=\{(a,\tau(a)) : a\in N\}$. Consider the von Neumann algebra $\tilde{W} = N\oplus^{\infty} \tau(N)^{op}$. Observe that $H(W,\tau)$ is a weak$^*$-closed C$^*$-subalgebra of $\tilde{W}$. Finally, the mapping $a\mapsto (a,\tau (a))$ is an isometric C$^*$-isomorphism from $N$ onto $H(W,\tau)$ when the latter is regarded as a von Neumann subalgebra of $\tilde{W}$, which implies that $N$ is a factor.\smallskip

If $q=0,$ $p = \unit$ and every projection in $Z$ is $\tau$-symmetric, and hence a projection in $H(W,\tau)$. Since the latter is a factor, the von Neumann algebra $W$ is a factor too.\smallskip

To see the reciprocal implication, we observe that if $(b)$ holds, then $H(W,\tau)$ is clearly a factor. If $W$ is a factor, by observing that every central projection in $H(W, \tau)$ is a central projection in $W$, the desired conclusion is clear (see  \cite[Lemma 7.3.2]{HOS}). 
\end{proof}

\subsection{Spin factors}\label{subsec: spin factors}

An important example of JBW$^*$-algebra, which also has an important role in this paper, is given by the so-called (complex) spin factors. The notion of spin factor dates back to the analytic classification of bounded symmetric domains in complex Banach spaces by authors like E. Cartan,  W. Kaup, L. Harris (cf. \cite{Harris1974,Harris1981} and the introduction of \cite{HervIsi1992} or \cite[\S 3]{FriedmanBook}). The construction is as follows. Let $H$ be a complex Hilbert space. A \emph{spin factor} (also known as a \emph{spinor} or \emph{Cartan factor of type IV}) is a norm closed complex subspace $V$ of $B(H)$ such that dim$(V)>2$, $a^*\in V,$ and $a^2 \in \mathbb{C} \unit,$ for all $a\in V$, where $\unit$ denotes the unit of $B(H)$. For example, the space $S_2(\mathbb{C}),$ of all symmetric $2\times 2$ complex matrices is an spin factor. \smallskip

It is further known that there exists an inner product $\langle\cdot| \cdot\rangle$ on $V$ satisfying 
$$ a b^* + b^* a= 2\langle a |b\rangle \unit, \hbox{ for all } a,b\in V.$$ The norm $\| x\|_2^2:=\langle a |a\rangle$ given by the inner product is equivalent to the operator norm $\|\cdot\|$, and both are related by $$\| a \|^2 =  \|a\|_2^2 + \left( \|a\|_2^4
- | \langle a |a^*\rangle |^2 \right)^{\frac12}.$$ It is also known that $V$ is closed for the triple product given by $$\{a,b,c\} = \frac12 (a b^* c + c b^* a) = \langle a |b\rangle c + \langle c |b\rangle a - \langle a |c^* \rangle b^* \ \ (a,b,c\in V).$$ For each norm-one element $u\in V$ with $u^* =u$, the Banach space $V$ is a unital JW$^*$-algebra with unit $u$, Jordan product $a\circ_u b = \{a,u,b\} = \frac12 (a u b + b u a),$  and involution $a^{*_u} : =\{u, a, u\} = u a^* u$ (cf. \cite[page 358]{Harris1981}). Every element $a$ in the JW$^*$-algebra $(V,\circ_u, *_u)$ satisfies $$a\circ_u a = \{a,u,a\} = 2 \langle a |u \rangle a - \langle a |a^*\rangle u,$$ that is, the square of every element $a$ in $(V,\circ_u, *_u)$ is a linear combination of $a$ and the unit element. 

\begin{remark}\label{r quadratic} We recall that a non-necessarily associative algebra $\mathfrak{J}$ is called \emph{quadratic} if it is unital and the square of every element $a$ in $\mathfrak{J}$ lies in the linear hull of $\{1,a\}$ (cf. \cite[2.2.5]{HOS} or \cite[\S 2.5.9]{CabGarPalVol1}). It is known that every quadratic JB$^*$-algebra is a spin factor (see \cite[Theorem 3.5.5 and Corollary 3.5.7]{CabGarPalVol1}).
\end{remark}

It is worth to note that if $V$ is a spin factor, the self-adjoin part of the associated JW$^*$-algebra $V_{sa} = \{a\in V : a^{*_u} = a \}$ is a real spin factor in the sense of \cite[\S 6]{HO1983}.

\begin{remark}\label{r associator-commutator of an element in a spin} It should be also noted that for each element $x$ in a spin factor $V$, the associator of $x$, that is the collection of all $y\in V$ operator commuting with $x$ reduces to $\mathbb{C}\unit \oplus \mathbb{C} x.$ Namely, let ${V}$ be a spin factor with inner product $\langle \cdot|\cdot \rangle,$ involution $a\mapsto \overline{a}$, unit $\unit$ (i.e. a norm-one element with $\overline{\unit} =\unit$), and Jordan product  $$a\circ b = \{a,\unit,b\} = \langle a |\unit \rangle b + \langle b |\unit \rangle a - \langle a |\overline{b} \rangle \unit, \hbox{ and } a^* = \{\unit, a, \unit\} = 2 \langle \unit |a\rangle \unit - \overline{a},$$ for all $a,b\in V.$ Observe that $V = \mathbb{C}\unit\oplus^{\ell_2} \{\unit\}^{\perp_2}$, where $ \{\unit\}^{\perp_2}$ stands for the orthogonal complement of $\{\unit\}$ in the underlying Hilbert space. Furthermore, for each $a,b\in \{\unit\}^{\perp_2},$ we have $a\circ b = - \langle a |\overline{b} \rangle \unit$. Pick $\lambda\unit + a, \mu \unit + b$ two operator commuting elements in $V$ (with $a,b\in \{\unit\}^{\perp_2}$). By assumptions, for each $c\in \{\unit\}^{\perp_2}$ we have $$\begin{aligned}
((\lambda\unit + a) \circ c)\circ (\mu \unit + b) &= (\lambda\unit + a)\circ (c\circ (\mu \unit + b)), \\
\lambda \mu c + \lambda c \circ b + \mu a \circ c + (a\circ c) \circ b  &= \lambda \mu c + \lambda c \circ b + \mu a \circ c + (b\circ c) \circ a, \\
- \langle a |\overline{c} \rangle b &= - \langle b |\overline{c} \rangle a,
\end{aligned}$$ and thus the arbitrariness of $c$ gives $b\in \mathbb{C} a,$ and consequently $\mu \unit + b$ lies in $\mathbb{C} \unit +\mathbb{C} (\lambda\unit + a),$ as desired.
\end{remark}

\section{Elementary operators on \texorpdfstring{JB$^*$}{JB*}-algebras}\label{sec: elementary operators} 

Let $\mathfrak{J}$ be a JB$^*$-algebra.
We write $\mathcal{B}(\mathfrak{J})$ for the Banach algebra of all bounded linear operators on the Banach space $\mathfrak{J}$, and
we write $\element(\mathfrak{J})$ for the subalgebra of $\mathcal{B}(\mathfrak{J})$
generated by the multiplication operators $M_a$ with $a$ running in 
$\mathfrak{J}$. The elements of $\element(\mathfrak{J})$ are called \emph{elementary operators} on $\mathfrak{J}$. 
Furthermore, we will denote by 
$\element_{\mathfrak{J}_{sa}}(\mathfrak{J})$ the real subalgebra of
$\element(\mathfrak{J})$ generated by the multiplication operators
$M_a$ with $a$ running in $\mathfrak{J}_{sa}$.
It is worth pointing out that each operator 
$\mathcal{E}\in \element_{\mathfrak{J}_{sa}}(\mathfrak{J})$
satisfies the condition
\begin{equation}\label{ele14}
\mathcal{E}(x^*)=\mathcal{E}(x)^*,\quad\forall x\in\mathfrak{J}.
\end{equation}
If $\mathfrak{J}$ is a JBW$^*$-algebra, then the subset of 
$\mathcal{B}(\mathfrak{J})$ consisting of all weak$^*$-continuous 
operators on $\mathfrak{J}$ is a subalgebra of 
$\mathcal{B}(\mathfrak{J})$ which, on account of the separate weak$^*$ continuity of the product of $\mathfrak{J}$, contains the
multiplication operators on $\mathfrak{J}$. This implies that 
each elementary operator on $\mathfrak{J}$ is weak$^*$-continuous.

\begin{remark}\label{rem:elementary Z-linear} Let $a,z$ be elements in a Jordan algebra $\mathfrak{J}$, with $z$ 
central. Clearly, $M_a M_z = M_a M_z,$ and thus, by linearity, $
\mathcal{E} M_z = M_z \mathcal{E},$ for every elementary operator 
$\mathcal{E}\in \element(J)$. 
\end{remark}

\begin{lemma}{\label{lemma:elem_homo}}
Let $\mathfrak{J}$  and $\mathfrak{H}$ be two Jordan algebras, and let $\phi\colon \mathfrak{J} \rightarrow \mathfrak{H}$ be a Jordan homomorphism. The the following statements hold.
    \begin{enumerate}[\rm(i)]
        \item
        For each $\mathcal{E}\in\element(\mathfrak{J})$ there exists
        $\mathcal{F}\in\element(\mathfrak{H})$ such that 
        $\phi\mathcal{E}=\mathcal{F}\phi$.
        \item
        Suppose that $\mathfrak{H}$ is a JBW$^*$-algebra and 
        that $\phi(\mathfrak{J})$ is weak$^*$-dense in $\mathfrak{H}$.
        Then for each $\mathcal{E}\in\element(\mathfrak{J})$ there exists a unique 
        $\Phi(\mathcal{E})\in\element(\mathfrak{H})$ 
        such that $\phi\mathcal{E}=\Phi(\mathcal{E})\phi$. Further,
        the map $\Phi\colon\element(\mathfrak{J})\to\element(\mathfrak{H}),$ $\mathcal{E}\mapsto \Phi(\mathcal{E})$ 
        is a homomorphism.
        \end{enumerate}
\end{lemma}

\begin{proof}
(i) It is immediate to check that the set
\[
\mathcal{L}=
\bigl\{
\mathcal{E}\in\element(\mathfrak{J}): 
\phi\mathcal{E} =\mathcal{F}\phi \textup{ for some } 
\mathcal{F}\in\element(\mathfrak{H})
\bigr\}
\]
is a subalgebra of $\element(\mathfrak{J})$. 
Moreover,
for each $a\in\mathfrak{J}$, we have
$\phi M_a=M_{\phi(a)}\phi$, and therefore $M_a\in\mathcal{L}$.
Consequently, $\mathcal{L}=\element(\mathfrak{J})$.\smallskip

\noindent(ii) Suppose that $\mathcal{F}_1,\mathcal{F}_2\in\element(\mathfrak{H})$ are
such that $\mathcal{F}_1\phi=\mathcal{F}_2\phi$.
From the weak$^*$-continuity of both maps
$\mathcal{F}_1$ and $\mathcal{F}_2$ and
the weak$^*$-density of $\phi(\mathfrak{J})$ in $\mathfrak{H}$ we 
deduce that $\mathcal{F}_1=\mathcal{F}_2$.
We conclude from (i) that, for each 
$\mathcal{E}\in\element(\mathfrak{J})$, there exists a unique 
$\Phi(\mathcal{E})\in\element(\mathfrak{H})$ 
such that $\phi\mathcal{E}=\Phi(\mathcal{E})\phi$. 
Hence, we can define a map 
$\Phi\colon\element(\mathfrak{J})\to\element(\mathfrak{H})$ 
through the condition
$\phi\mathcal{E}=\Phi(\mathcal{E})\phi$ for each 
$\mathcal{E}\in\element(\mathfrak{J})$,
and routine verifications show that $\Phi$ is a homomorphism.
\end{proof}

Our next result is one of the key tools in our arguments. The statement combines non-associative Jordan algebras, and in particular JBW$^*$-algebras with a rich geometric--algebraic structure, and associative algebras regarded as Jordan algebras via their natural Jordan product. 

\begin{proposition}\label{cor:3_lib}
Let $\mathfrak{J}$ be a JBW$^*$-algebra with no direct summands of type $I_1$ or $I_2$. 
Then there exist $u,v\in\mathfrak{J}_{sa}$ and 
$\mathcal{E}_0,\mathcal{E}_1,\mathcal{E}_2\in\element_{\mathfrak{J}_{sa}}(\mathfrak{J})$ 
with the following properties:
\begin{enumerate}[\rm(i)]
\item
$\mathcal{E}_i(u^j)=\delta_{ij}\mathbf{1},$
for all $i,j\in\{0,1,2\}$,
\item
for each unital complex associative algebra $\mathfrak{A}$ and each unital Jordan
homomorphism $\phi\colon\mathfrak{J}\to\mathfrak{A}$, the elements
\[
[\phi(u)^2,\phi(v)],
[\phi(u),\phi(v)^2], \hbox{ and }
[\phi(u),\phi(v)]
\]
are linearly independent,  where the brackets stand for the usual Lie product on $\mathfrak{A}$.
\end{enumerate} We can further assume that $\| \mathcal{E}_i\| \leq 10$ for every $i\in \{0,1,2\}.$
\end{proposition}

\begin{proof} 
It is well known that $\mathfrak{J}_{sa}$ is a JBW-algebra, and further $\mathfrak{J}$ and 
$\mathfrak{J}_{sa}$ have the same lattice of projections (cf. \cite{Edwards80,Wright1977}). Therefore, $\mathfrak{J}_{sa}$ is a 
JBW-algebra with no direct summands of type $I_1$ or $I_2$. 
The proof of the Proposition will be divided into a sequence of cases according to the structure of $\mathfrak{J}_{sa}$. \smallskip

\noindent\emph{Case 1:} We assume first that $\mathfrak{J}_{sa}$ is such that there exist pairwise orthogonal projections $p_1,p_2,p_3,p_4$ which are pairwise exchangeable by symmetries such that
\[ p_1+p_2+p_3+p_4=\mathbf{1}.\] 
It should be pointed out that this condition holds in each of the following cases:
	\begin{enumerate}[$(1)$]
		\item
		$\mathfrak{J}_{sa}$ is a JBW-algebra with no direct summand of type $I$ (see \cite[Proposition 5.2.15]{HOS}),
		\item
		$\mathfrak{J}_{sa}$ is a JBW-algebra of type $I_{\infty}$ (see \cite[Proposition 3.24]{AlfsenShultz2003}),
		\item
		$\mathfrak{J}=\mathfrak{J}_1\oplus^{\infty}\mathfrak{J}_2$, where $\mathfrak{J}_1$ and $\mathfrak{J}_2$ are JBW$^*$-algebras as those considered in (1) and (2), respectively.
	\end{enumerate}
 
\noindent (i)	Take symmetries $s,s',s''\in\mathfrak{J}_{sa}$ such that
$$ U_s(p_1)=p_2,\ U_{s'}(p_1)=p_3, \hbox{ and } U_{s''}(p_1)=p_4.$$
	Define $u,v\in\mathfrak{J}_{sa}$ by
 \[
 u=2s\circ p_1, \hbox{ and } v=2s'\circ p_1,
 \]
 so that, by the commented Shirshov--Cohn theorem, 
 \[
 u^2=p_1+p_2, \hbox{ and } v^2=p_1+p_3.
 \]
Then, $$U_{p_3}(\mathbf{1})=p_3,\ U_{p_3}(u)=0,\hbox{ and } U_{p_3}(u^2)=0,$$
	where in the middle term we compute the value in the JC-subalgebra generated by $p_1-p_2,$ $s$ and $\unit,$ and the fact $U_{\unit-p_1-p_2} U_{p_3} U_{\unit-p_1-p_2} = U_{U_{\unit-p_1-p_2} (p_3)} = U_{p_3},$ whence
	\begin{align*}
		U_{s'}U_{p_3}(\mathbf{1})&=p_1,&U_{s'}U_{p_3}(u)&=0,&U_{s'}U_{p_3}(u^2)&=0,\\
		U_sU_{s'}U_{p_3}(\mathbf{1})&=p_2,&U_sU_{s'}U_{p_3}(u)&=0,&U_sU_{s'}U_{p_3}(u^2)&=0,\\
		U_{s''}U_{s'}U_{p_3}(\mathbf{1})&=p_4,&U_{s''}U_{s'}U_{p_3}(u)&=0,&U_{s''}U_{s'}U_{p_3}(u^2)&=0.
	\end{align*}
	By defining
	\[
	\mathcal{E}_0=U_{p_3}+U_{s'}U_{p_3}+U_sU_{s'}U_{p_3}+U_{s''}U_{s'}U_{p_3},
	\]
	we have
 $\mathcal{E}_0\in\element_{\mathfrak{J}_{sa}}(\mathfrak{J})$ and
	\[
	\mathcal{E}_0(\mathbf{1})=\mathbf{1},\quad
	\mathcal{E}_0(u)=0,\quad
	\mathcal{E}_0(u^2)=0.
	\]
	Observe that
	\begin{align*}
		(U_{p_1}-U_{s'}U_{p_3})(\mathbf{1})&=0, & (U_{p_1}-U_{s'}U_{p_3})(u)&=0, & (U_{p_1}-U_{s'}U_{p_3})(u^2)&=p_1,
	\end{align*}
	and therefore
	\begin{align*}
		U_s(U_{p_1}-U_{s'}U_{p_3})(\mathbf{1})&=0, & U_s(U_{p_1}-U_{s'}U_{p_3})(u)&=0, & U_s(U_{p_1}-U_{s'}U_{p_3})(u^2)&=p_2,\\
		U_{s'}(U_{p_1}-U_{s'}U_{p_3})(\mathbf{1})&=0, & U_{s'}(U_{p_1}-U_{s'}U_{p_3})(u)&=0, & U_{s'}(U_{p_1}-U_{s'}U_{p_3})(u^2)&=p_3,\\
		U_{s''}(U_{p_1}-U_{s'}U_{p_3})(\mathbf{1})&=0, & U_{s''}(U_{p_1}-U_{s'}U_{p_3})(u)&=0, & U_{s''}(U_{p_1}-U_{s'}U_{p_3})(u^2)&=p_4.
	\end{align*} By defining
	\[
	\mathcal{E}_2=
	(U_{p_1}-U_{s'}U_{p_3})+
	U_s(U_{p_1}-U_{s'}U_{p_3})+
	U_{s'}(U_{p_1}-U_{s'}U_{p_3})+
	U_{s''}(U_{p_1}-U_{s'}U_{p_3}),
	\] we obtain
 $\mathcal{E}_2\in\element_{\mathfrak{J}_{sa}}(\mathfrak{J})$ and
	\[
	\mathcal{E}_2(\mathbf{1})=0,\quad
	\mathcal{E}_2(u)=0,\quad
	\mathcal{E}_2(u^2)=\mathbf{1}.
	\]
	Finally, note that
$$M_u(\mathbf{1})=u,\ M_u(u)=u^2, \hbox{ and } M_u(u^2)=u.$$
	Consequently, by defining $\mathcal{E}_1=\mathcal{E}_2M_u,$ we get
	\[
	\mathcal{E}_1(\mathbf{1})=0,\quad
	\mathcal{E}_1(u)=\mathbf{1},\hbox{ and }
	\mathcal{E}_1(u^2)=0.
	\]	

\noindent (ii) Let $\mathfrak{A}$ be a unital complex associative algebra, and
let $\phi\colon\mathfrak{J}\to\mathfrak{A}$ be a unital Jordan homomorphism.
Note that, for all $i,j\in\{1,2,3\}$ with $i\ne j$, we have
\[
\phi(p_i)^2=
\phi(p_i\circ p_i)=\phi(p_i),
\]
$$\begin{aligned}
   \phi(p_i)\phi(p_j)&=
\phi\bigl(\underbrace{U_{p_j}p_i}_{=0}\bigr)+\phi(p_i)\phi(p_j)=
U_{\phi(p_j)}\phi(p_i)+\phi(p_i)\phi(p_j)^2 \\
&=
2\bigl(\phi(p_i)\circ\phi(p_j)\bigr)\phi(p_j)=
2\phi(\underbrace{p_i\circ p_j}_{=0})\phi(p_j)=
2\phi(0)\phi(p_j)=0, 
\end{aligned}$$
\[
\phi(s)^2=\phi(s\circ s)=\phi(\mathbf{1})=\mathbf{1},
\]
\[
\phi(s')^2=
\phi(s'\circ s')=
\phi(\mathbf{1})=\mathbf{1},
\]
\[
\phi(s)\phi(p_1)=
\phi(s)\phi(p_1)\phi(s)^2=
U_{\phi(s)}\phi(p_1)\phi(s)=
\phi\bigl(U_{s}p_1\bigr)\phi(s)=
\phi(p_2)\phi(s),
\]
\[
\phi(p_1)\phi(s)=
\phi(s)^2\phi(p_1)\phi(s)=
\phi(s)U_{\phi(s)}\phi(p_1)=
\phi(s)\phi\bigl(U_{s}p_1\bigr)=
\phi(s)\phi(p_2),
\]
\[
\phi(s')\phi(p_1)=
\phi(s')\phi(p_1)\phi(s')^2=
U_{\phi(s')}\phi(p_1)\phi(s')=
\phi\bigl(U_{s'}p_1\bigr)\phi(s')=
\phi(p_3)\phi(s'),
\]
\[
\phi(p_1)\phi(s')=
\phi(s')^2\phi(p_1)\phi(s')=
\phi(s')U_{\phi(s')}\phi(p_1)=
\phi(s')\phi\bigl(U_{s'}p_1\bigr)=
\phi(s')\phi(p_3).
\]
Therefore
\begin{equation}\label{eq 2805}\left\{
\begin{aligned}
\relax[\phi(u)^2,\phi(v)] &=
[\phi(p_1)+\phi(p_2),\phi(s')\phi(p_1)+\phi(p_1)\phi(s')] \\ &=
\phi(p_1)\phi(s')\phi(p_1)+
\phi(p_1)^2\phi(s')+
\phi(p_2)\phi(s')\phi(p_1) \\ &+
\phi(p_2)\phi(p_1)\phi(s')
-\phi(s')\phi(p_1)^2-
\phi(p_1)\phi(s')\phi(p_1) \\ &-
\phi(s')\phi(p_1)\phi(p_2)-
\phi(p_1)\phi(s')\phi(p_2) \\ &=
\phi(p_1)\phi(s')-\phi(s')\phi(p_1) +
\underbrace{\phi(p_2)\phi(p_3)}_{=0}\phi(s') \\ 
+ &
\underbrace{\phi(p_2)\phi(p_1)}_{=0}\phi(s') -\phi(s')\underbrace{\phi(p_1)\phi(p_2)}_{=0}
-\phi(s')\underbrace{\phi(p_3)\phi(p_2)}_{=0} \\ &=
\phi(p_1)\phi(s')-\phi(s')\phi(p_1),
\end{aligned}  \right.  
\end{equation}
and similarly
\begin{equation}\label{eq 2805b}
[\phi(u),\phi(v)^2]=
\phi(s)\phi(p_1)-\phi(p_1)\phi(s).
\end{equation}
On the other hand,
\begin{equation}\label{eq 2805c}\left\{
\begin{aligned}
\relax[\phi(u),\phi(v)] &=
[\phi(s)\phi(p_1)+\phi(p_1)\phi(s),\phi(s')\phi(p_1)+\phi(p_1)\phi(s')] \\ &=
\phi(s)\phi(p_1)\phi(s')\phi(p_1)+
\phi(s)\phi(p_1)^2\phi(s') \\ &+
\phi(p_1)\phi(s)\phi(s')\phi(p_1)+
\phi(p_1)\phi(s)\phi(p_1)\phi(s') \\ &-\phi(s')\phi(p_1)\phi(s)\phi(p_1)
-\phi(p_1)\phi(s')\phi(s)\phi(p_1) \\ &-\phi(s')\phi(p_1)^2\phi(s)
-\phi(p_1)\phi(s')\phi(p_1)\phi(s) \\ &=
\phi(s)\underbrace{\phi(p_1)\phi(p_3)}_{=0}\phi(s')+
\phi(s)\phi(p_1)\phi(s') \\ &+
\phi(s)\underbrace{\phi(p_2)\phi(p_3)}_{=0}\phi(s')+
\phi(s)\underbrace{\phi(p_2)\phi(p_1)}_{=0}\phi(s') \\ &-\phi(s')\underbrace{\phi(p_1)\phi(p_2)}_{=0}\phi(s)
-\phi(s')\underbrace{\phi(p_3)\phi(p_2)}_{=0}\phi(s) \\ &-\phi(s')\phi(p_1)\phi(s)
-\phi(s')\underbrace{\phi(p_3)\phi(p_1)}_{=0}\phi(s) \\ &=
\phi(s)\phi(p_1)\phi(s')-\phi(s')\phi(p_1)\phi(s).
\end{aligned}  \right.  
\end{equation}

Observe that
$$\begin{aligned}
 \mathbf{1}&=\phi(\mathbf{1})=
\phi(p_1)+\phi(p_2)+\phi(p_3)+\phi(p_4) \\ &= \phi(p_1)+
\phi\bigl(U_sp_1\bigr)+
\phi\bigl(U_{s'}p_1\bigr)+
\phi\bigl(U_{s''}p_1\bigr) \\ &=
\phi(p_1)+
\phi(s)\phi(p_1)\phi(s)+
\phi(s')\phi(p_1)\phi(s')+
\phi(s'')\phi(p_1)\phi(s'').
\end{aligned}$$ The previous identity implies that $\phi(p_1)\ne 0$, and, since
$$
\begin{aligned}
\phi(p_1) = \phi(s)\phi(p_2)\phi(s)  = \phi(s')\phi(p_3)\phi(s') = \phi(s'')\phi(p_4)\phi(s''),
\end{aligned}
$$
it follows that
\[
\phi(p_2),\phi(p_3),\phi(p_4)\ne 0.
\]

Suppose that $\lambda,\mu,\nu\in\mathbb{C}$ are such that
\[
\lambda[\phi(u)^2,\phi(v)]+
\mu[\phi(u),\phi(v)^2]+
\nu[\phi(u),\phi(v)]=0.
\]
Then, having in mind \eqref{eq 2805}, \eqref{eq 2805b} and \eqref{eq 2805c}, we get 
$$\begin{aligned}
 0 &=
\bigl(\lambda[\phi(u)^2,\phi(v)]+
\mu[\phi(u),\phi(v)^2]+
\nu[\phi(u),\phi(v)]\bigr)\phi(p_1)\\ &=
\lambda\bigl(\phi(p_1)\phi(s')\phi(p_1)-\phi(s')\phi(p_1)^2\bigr) \\ &+
\mu\bigl(\phi(s)\phi(p_1)^2-\phi(p_1)\phi(s)\phi(p_1)\bigr)\\ &+
\nu\bigl(\phi(s)\phi(p_1)\phi(s')\phi(p_1)-\phi(s')\phi(p_1)\phi(s)\phi(p_1)\bigr)\\ &=
\lambda\bigl(\underbrace{\phi(p_1)\phi(p_3)}_{=0}\phi(s')-\phi(s')\phi(p_1)\bigr)\\ &+
\mu\bigl(\phi(s)\phi(p_1)-\underbrace{\phi(p_1)\phi(p_2)}_{=0}\phi(s)\bigr)\\ &+
\nu\bigl(\phi(s)\underbrace{\phi(p_1)\phi(p_3)}_{=0}\phi(s')-\phi(s')\underbrace{\phi(p_1)\phi(p_2)}_{=0}\phi(s)\bigr)\\ &=
-\lambda\phi(s')\phi(p_1)+\mu\phi(s)\phi(p_1)
=-\lambda\phi(p_3)\phi(s')+\mu\phi(p_2)\phi(s),   
\end{aligned}$$
whence
\[
0=\phi(p_3)\bigl(-\lambda\phi(p_3)\phi(s')+
\mu\phi(p_2)\phi(s)\bigr)
=
-\lambda\phi(p_3)\phi(s')
\]
and
\[
0=\bigl(-\lambda\phi(p_3)\phi(s')\bigr)\phi(s')=
-\lambda\phi(p_3),
\]
which gives $\lambda =0$.
Therefore
\[
\mu\phi(p_2)\phi(s)=0,
\]
and
\[
\mu\phi(p_2)=\mu\phi(p_2)\phi(s)^2= \mu 
\bigl(\phi(p_2)\phi(s)\bigr)\phi(s)=0,
\]
which yields $\mu=0$.
We thus obtain
\[
\nu[\phi(u),\phi(v)]=0.
\]
Hence
$$\begin{aligned}
0 &=\nu[\phi(u),\phi(v)]\phi(p_3)=
\nu\bigl(\phi(s)\phi(p_1)\phi(s')\phi(p_3)-
\phi(s')\phi(p_1)\phi(s)\phi(p_3)\bigr) \\ &=
\nu\phi(s)\phi(s')\underbrace{\phi(p_3)^2}_{=\phi(p_3)}-
\phi(s')\phi(s)\underbrace{\phi(p_1)\phi(p_3)}_{=0}
=
\nu\phi(s)\phi(s')\phi(p_3),
\end{aligned}$$
and therefore
\[
0=\phi(s)\phi(s')\bigl(\nu\phi(s)\phi(s')\phi(p_3)\bigr)=
\nu\phi(p_3),
\]
which gives $\nu=0$.

\noindent\emph{Case II:} We assume next that $\displaystyle \mathfrak{J}_{sa}=\bigoplus_{n\in N}^{\infty} J_n,$	where
	$N\subset\{n\in\mathbb{N}\colon n\ge 3\}$ and
	$J_n$ is a JBW-algebra of type $I_n$ for each $n\in N$ (cf. \cite[5.3.3 and Theorem 5.3.5]{HOS}).
	For each $n\in N$ take $j_n\in\mathbb{N}$ and $k_n\in\{0,1,2\}$ such that
	\[
	n=3j_n+k_n,
	\]
	and take $e_{n,1},\ldots,e_{n,n}$ pairwise orthogonal projections in $J_n$ such that
	\[
	e_{n,1}+\cdots+e_{n,n}=\mathbf{1}_{J_n}
	\]
	and such that they are pairwise exchangeable by a symmetry (cf. \cite[Definition 3.21]{AlfsenShultz2003} or \cite[5.3.3]{HOS}).
	We define pairwise orthogonal projections $p_1,p_2,p_3,q_1,q_2$ by
	\begin{equation*}
		\begin{split}
			p_1&=\sum_{n\in N}\sum_{j=0}^{j_n-1}e_{n,1+3j},\ \
			p_2=\sum_{n\in N}\sum_{j=0}^{j_n-1}e_{n,2+3j},\ \
			p_3=\sum_{n\in N}\sum_{j=0}^{j_n-1}e_{n,3+3j},\\
			q_1& =\sum_{n\in N\colon k_n\in\{1,2\}}e_{n,1+3j_n}, \ \ \hbox{ and }  \ \ \ \ 
			q_2=\sum_{n\in N\colon k_n=2}e_{n,2+3j_n}.
		\end{split}
	\end{equation*}
	Observe that $q_1$ and $q_2$ are possibly zero and 
	\[
	p_1+p_2+p_3+q_1+q_2=\mathbf{1},
	\]
	and, further, $p_1,p_2,p_3$ are pairwise exchangeable by a symmetry (cf. \cite[Lemma 5.2.9]{HOS}).\smallskip

\noindent (i)
	Take symmetries $s,s'\in\mathfrak{J}_{sa}$ such that $$ U_s(p_1)=p_2,\hbox{ and }  U_{s'}(p_1) = p_3.$$
	Define $u,v\in\mathfrak{J}_{sa}$ by
	\[
	u=2s\circ p_1, \hbox{ and } v=2s'\circ p_1,
	\]
	so that, as in the previous case, 
	\[
	u^2=p_1+p_2, \hbox{ and } v^2=p_1+p_3.
	\]
	Then
$$ U_{p_3}(\mathbf{1})=p_3, \ \ \  U_{p_3}(u)=0, \  \hbox{ and } \  U_{p_3}(u^2)=0,$$	which yields
	\begin{align*}
		U_{s'}U_{p_3}(\mathbf{1})&=p_1, & U_{s'}U_{p_3}(u)&=0 ,& U_{s'}U_{p_3}(u^2)&=0,\\
		U_sU_{s'}U_{p_3}(\mathbf{1})&=p_2 ,&U_sU_{s'}U_{p_3}(u)&=0 ,& U_sU_{s'}U_{p_3}(u^2)&=0.
	\end{align*}
	Further, $$ U_{q_1+q_2}(\mathbf{1}) =q_1+q_2 ,\ \ U_{q_1+q_2}(u)=0,\ \ \hbox{ and } \ \ U_{q_1+q_2}(u^2)=0.$$
	We define $\mathcal{E}_0\in\element_{\mathfrak{J}_{sa}}(\mathfrak{J})$ by
	\[
	\mathcal{E}_0=U_{p_3}+U_{s'}U_{p_3}+U_sU_{s'}U_{p_3}+U_{q_1+q_2},
	\]
	and we check at once that
	\[
	\mathcal{E}_0(\mathbf{1})=\mathbf{1},\quad\mathcal{E}_0(u)=0,\hbox{ and } \mathcal{E}_0(u^2)=0.
	\]
	We now define projections $r_1,r_2\in\mathfrak{J}_{sa}$ by
	\begin{equation*}
		\begin{split}
			r_1=\sum_{n\in N\colon k_n\in\{1,2\}}e_{n,2}\leq p_2, \hbox{ and } \ 
			r_2=\sum_{n\in N\colon k_n=2}e_{n,2}\leq p_2.
		\end{split}
	\end{equation*}
	Note that the projections $r_1$ and $q_1$ are exchanged by a symmetry and the projections 
	$r_2$ and $q_2$ are also exchanged by a symmetry (cf. \cite[Lemma 5.2.9]{HOS}). Take symmetries $t,t'\in \mathfrak{J}$ such that
	$ U_{t}(r_1)=q_1,$  and $U_{t'}(r_2)=q_2.$
	Observe that
$$ (U_{p_1}-U_{s'}U_{p_3})(\mathbf{1})=0, \ (U_{p_1}-U_{s'}U_{p_3})(u)=0, \hbox{ and } (U_{p_1}-U_{s'}U_{p_3})(u^2)=p_1,$$
	and therefore
	\begin{align*}
		U_s(U_{p_1}-U_{s'}U_{p_3})(\mathbf{1})&=0, & U_s(U_{p_1}-U_{s'}U_{p_3})(u)&=0, & U_s(U_{p_1}-U_{s'}U_{p_3})(u^2)&=p_2,\\
		U_{s'}(U_{p_1}-U_{s'}U_{p_3})(\mathbf{1})&=0, & U_{s'}(U_{p_1}-U_{s'}U_{p_3})(u)&=0,  &  U_{s'}(U_{p_1}-U_{s'}U_{p_3})(u^2)&=p_3.
	\end{align*}
	On the other hand,
$$(U_{r_1}-U_{t}U_{q_1})(\mathbf{1})=0,  \ (U_{r_1}-U_{t}U_{q_1})(u)=0, \hbox{ and }  (U_{r_1}-U_{t}U_{q_1})(u^2)=r_1,$$ which gives
	$$ U_t(U_{r_1}-U_{t}U_{q_1})(\mathbf{1})=0, \  U_t(U_{r_1}-U_{t}U_{q_1})(u)=0, \hbox{ and } U_t(U_{r_1}-U_{t}U_{q_1})(u^2)=q_1,$$ and furthermore
	 $$ (U_{r_2}-U_{t'}U_{q_2})(\mathbf{1})=0, \  (U_{r_2}-U_{t'}U_{q_2})(u)=0, \hbox{ and }  (U_{r_2}-U_{t'}U_{q_2})(u^2) =r_2,$$
	so that $$
		U_{t'}(U_{r_2}-U_{t'}U_{q_2})(\mathbf{1})=0, \ U_{t'}(U_{r_2}-U_{t'}U_{q_2})(u)=0, \hbox{ and } U_{t'}(U_{r_2}-U_{t'}U_{q_2})(u^2)=q_2.$$
	By defining
	$$\begin{aligned}
	 \mathcal{E}_2 =
	&(U_{p_1}-U_{s'}U_{p_3})+
	U_s(U_{p_1}-U_{s'}U_{p_3})+
	U_{s'}(U_{p_1}-U_{s'}U_{p_3}) \\
 &+U_t(U_{r_1}-U_{t}U_{q_1})+
	U_{t'}(U_{r_2}-U_{t'}U_{q_2}),   
	\end{aligned}$$
	we obtain
 $\mathcal{E}_2\in\element_{\mathfrak{J}_{sa}}(\mathfrak{J})$ and
	\[
	\mathcal{E}_2(\mathbf{1})=0,\quad
	\mathcal{E}_2(u)=0,\quad
	\mathcal{E}_2(u^2)=\mathbf{1}.
	\]
	Finally, by defining
	\[
	\mathcal{E}_1=\mathcal{E}_2M_u
	\]
	we get
 $\mathcal{E}_1\in\element_{\mathfrak{J}_{sa}}(\mathfrak{J})$ and
	\[
	\mathcal{E}_1(\mathbf{1})=0,\quad
	\mathcal{E}_1(u)=\mathbf{1},\quad
	\mathcal{E}_1(u^2)=0.
	\]

\noindent (ii)
Let $\mathfrak{A}$ be a unital complex associative algebra, and
let $\phi\colon\mathfrak{J}\to\mathfrak{A}$ a unital Jordan homomorphism.\smallskip

A similarly argument to that employed in the previous \emph{Case I} gives 
$$
\begin{aligned}
\relax[\phi(u)^2,\phi(v)]&=\phi(p_1)\phi(s')-\phi(s')\phi(p_1),\\
[\phi(u),\phi(v^2)]&=\phi(s)\phi(p_1)-\phi(p_1)\phi(s), \hbox{ and }\\
[\phi(u),\phi(v)] \ &= \phi(s)\phi(p_1)\phi(s')-\phi(s')\phi(p_1)\phi(s).
\end{aligned}$$
Next, in order to prove that  $\phi(p_1),\phi(p_2),\phi(p_3)$ are all non-zero, we shall argue in a different way. Let us assume, contrary to the desired statement, that $\phi(p_1)=0$. Then we have 
$\phi(u)=2\phi(s)\circ\phi(p_1)=0$ and thus \Cref{lemma:elem_homo}(i) implies the existence of $\mathcal{F}_1\in \element({\mathfrak{A}})$
\[
0=\mathcal{F}_1 (\phi(u))=\phi\bigl(\mathcal{E}_1(u)\bigr)=\phi(\mathbf{1})=\mathbf{1},
\] which is impossible. Now, as in \emph{Case I}, we have
\[ \phi(p_1)=\phi(s)\phi(p_2)\phi(s) =\phi(s')\phi(p_3)\phi(s'),
\]
 which implies that $\phi(p_2),$ and $ \phi(p_3)$ must be both non-zero.
 The linear independence of $
[\phi(u)^2,\phi(v)],$ $ [\phi(u),\phi(v)^2],$ and $
[\phi(u),\phi(v)]$ follows as in \emph{Case I} too.\smallskip

\noindent\emph{General case.} 
By \cite[Theorems 5.1.5 and 5.3.5]{HOS},
$\mathfrak{J}_{sa}$ admits a unique decomposition in the form 
\[
\mathfrak{J}_{sa}=
\underbrace{J_{I_1}\oplus J_{I_2}\oplus\ldots}_{J_2}
\oplus
\underbrace{J_{I_\infty}
\oplus
J_{II}
\oplus
J_{III}}_{J_1},
\]
where each $J_k$ is either $0$ or a JBW-algebra in the class studied in \emph{Case $k$} above. This implies that $\mathfrak{J}$ splits into a direct sum
\[
\mathfrak{J} = 
\mathfrak{J}_1\oplus\mathfrak{J}_2,
\]
where $\mathfrak{J}_1$ and $\mathfrak{J}_2$ are JBW$^*$-algebras such that
${\mathfrak{J}_1}_{sa}=J_1$ and ${\mathfrak{J}_2}_{sa}=J_2$, where according to the observation made at the beginning of the discussion
of \emph{Case I},  $\mathfrak{J}_1$ fits the condition required in that case. Since $\mathfrak{J}$ contains no direct summands of type $I_1$ nor $I_2$, it follows that $J_{I_1}=J_{I_2}=0$ and therefore
$\mathfrak{J}_2$ fits the assumptions in \emph{Case II}. We denote by $\unit_1$ and $\unit_2$ the units of $\mathfrak{J}_1$ 
and $\mathfrak{J}_2$, respectively.\smallskip

Let $u_1,v_1\in \mathfrak{J}_1$ and 
$\mathcal{E}_{1,0},\mathcal{E}_{1,1}$, $\mathcal{E}_{1,2}$ denote the
elements and the elementary operators given by \emph{Case I} for 
$\mathfrak{J}_1$, and let $u_2,v_2\in\mathfrak{J}_2$ and 
$\mathcal{E}_{2,0},\mathcal{E}_{2,1}$, 
$\mathcal{E}_{2,2}$ the corresponding elements and elementary operators
for  
$\mathfrak{J}_2$ whose existence is guaranteed by \emph{Case II}. 
Setting 
\[
u=u_1+u_2, \  v=v_1+v_2\in\mathfrak{J}_{sa},
\]
and 
	\[ \mathcal{E}_i=\mathcal{E}_{1,i} U_{\unit_1}+\mathcal{E}_{2,i} 
 U_{\unit_2}\in\element_{\mathfrak{J}_{sa}}(\mathfrak{J}) \quad 
 (i\in \{0,1,2\}),
	\]
 we get the dessired conclusion in (i).\smallskip

We shall finally prove the statement in (ii). For this purpose,  let $\mathfrak{A}$ be a unital complex associative algebra, and
let $\phi\colon \mathfrak{J}\to \mathfrak{A}$ be a unital Jordan homomorphism.
Then 
\[
\phi(\unit_k)^2=\phi(\unit_k^2)=\phi(\unit_k),\quad k\in\{1,2\},
\]
and
\[
\phi(\unit_1)\phi(\unit_2)=\phi(\unit_2)\phi(\unit_1)=0.
\]
We define subalgebras $\mathfrak{A}_1$ and $\mathfrak{A}_2$ of $\mathfrak{A}$ by
\[
\mathfrak{A}_k=\bigl\{a\in \mathfrak{A}\colon a\phi(\unit_k)=\phi(\unit_k)a=a\bigr\},
\quad k\in\{1,2\},
\]
and we observe that 
\begin{equation}\label{Aort}
\mathfrak{A}_1 \mathfrak{A}_2=\mathfrak{A}_2 \mathfrak{A}_1=\{0\}.
\end{equation}
Further, for each $k\in\{1,2\}$ and each $x\in \mathfrak{J}_k$,
\[
\phi(x)=\phi\bigl(U_{\unit_k}x)=U_{\phi(\unit_k)}\phi(x)=
\phi(\unit_k)\phi(x)\phi(\unit_k),
\]
which immediately yields $\phi(x)\in \mathfrak{A}_k$,
and $\phi$ gives a  homomorphism from $\mathfrak{J}_k$ to $\mathfrak{A}_k$.
Since 
\[
\phi(\unit_1)+\phi(\unit_2)=\phi(\unit)=\unit,
\]
it follows that the elements $\phi(\unit_1)$ and $\phi(\unit_2)$ cannot be both zero.
Suppose that $\phi(\unit_1),\phi(\unit_2)\ne 0$.
Then each $\mathfrak{A}_k$ is unital with unit $\phi(\unit_k)$, $\phi$ yields a unital homomorphism from $\mathfrak{J}_k$ to $\mathfrak{A}_k$, and, consequently, 
\[
[\phi(u_k)^2,\phi(v_k)],
[\phi(u_k),\phi(v_k)^2], \hbox{ and }
[\phi(u_k),\phi(v_k)],
\]
are linearly independent.
Since
$$\begin{aligned}
\relax[\phi(u)^2,\phi(v)]& =
[\phi(u_1)^2,\phi(v_1)]+
[\phi(u_2)^2,\phi(v_2)],\\
[\phi(u),\phi(v)^2]& =
[\phi(u_1),\phi(v_1)^2]+
[\phi(u_2),\phi(v_2)^2],\\
[\phi(u),\phi(v)]& =
[\phi(u_1),\phi(v_1)]+
[\phi(u_2),\phi(v_2)],
\end{aligned}$$ by using \eqref{Aort} we obtain the desired linear independence in statement (ii).
We now suppose that $\phi(\unit_1)=0$. Then $\phi(\mathfrak{J}_1)=0$,
$\mathfrak{A}_2$ is a unital algebra, and $\phi$
gives a unital homomorphism from  $\mathfrak{J}_2$ to $\mathfrak{A}_2$. Hence
$[\phi(u)^2,\phi(v)] =
[\phi(u_2)^2,\phi(v_2)],$ $[\phi(u),\phi(v)^2] =
[\phi(u_2),\phi(v_2)^2],$ and $[\phi(u),\phi(v)]=
[\phi(u_2),\phi(v_2)],$ are linearly independent.
The same reasoning applies to the case $\phi(\unit_2)=0$. \smallskip

The statement concerning the norm of $\mathcal{E}_i$ is clear from the definitions. 
\end{proof}

The next result also considers JBW$^*$-algebras with non-necessarily zero type $I_2$ part. The conclusion is, of course, weaker than in the previous \Cref{cor:3_lib}. 

\begin{proposition}\label{cor:3_lib including spin part}
Let $\mathfrak{J}$ be a JBW$^*$-algebra with no direct summands of type $I_1$. 
Then there exist an element $u\in\mathfrak{J}_{sa}$ and 
$\mathcal{E}_0,\mathcal{E}_1\in\element_{\mathfrak{J}_{sa}}(\mathfrak{J})$ satisfying  $\mathcal{E}_i(u^j)=\delta_{ij}\mathbf{1},$ and $\| \mathcal{E}_i\|\leq 10$, for all $i,j$ in $\{0,1\}$. 
\end{proposition}

\begin{proof} By representation theory, $\mathfrak{J}$ decomposes as the orthogonal sum of two JBW$^*$-subalgebras $\mathfrak{J}= \mathfrak{J}_{spin}\oplus \mathfrak{J}_{nspin},$ where $\mathfrak{J}_{spin}$ is a JBW$^*$-algebra of type $I_2$ and $\mathfrak{J}_{nspin}$ contains no direct summands of type $I_1$ or $I_2$  (cf. \cite[Theorems 5.1.5 and 5.3.5]{HOS}). In the case of $\mathfrak{J}_{nspin}$, the desired conclusion follows from \Cref{cor:3_lib}. Arguing as in the final part of the proof of the just quoted proposition, and employing the orthogonality of $\mathfrak{J}_{spin}$ and $ \mathfrak{J}_{nspin}$, it suffices to show that the result also holds for $\mathfrak{J}_{spin}.$\smallskip

By structure theory, we can always find two orthogonal projections $p_1, p_2$ in $\mathfrak{J}_{spin}$ which are exchangeable by a symmetry $s\in \mathfrak{J}_{spin}$ (see \cite[5.3.3]{HOS}). Take $u =s$ and observe that $2 p_1\circ s = 2 p_2\circ s  = s.$ By defining $\mathcal{E}_1 = M_s M_{2 p_2} M_{2 p_1}$ and $\mathcal{E}_0 = \mathcal{E}_1 M_{s},$ it is easy to check that $\mathcal{E}_i(u^j)=\delta_{ij}\mathbf{1},$ for all $i,j$ in $\{0,1\}$.
\end{proof}

In order to obtain a refinement of \Cref{lemma:elem_homo}, we recall some basic results on the strong$^*$-topology of a JBW$^*$-algebra $\mathfrak{J}$. For each positive normal functional $\varphi$ in the predual, $\mathfrak{J}_{*}$, of $\mathfrak{J}$, the sesquilinear form $(x,y)\mapsto \varphi (x\circ y^*)$ is positive, and defines a preHilbertian seminorm $\|a\|_{\varphi}^2 := \varphi (x\circ x^*)$ on $\mathfrak{J}$.  The \emph{strong$^*$ topology} of $\mathfrak{J},$ denoted by $S^*(\mathfrak{J},\mathfrak{J}_*)$, is the topology generated by all the seminorms $\|\cdot \|_{\varphi}$ with $\varphi$ running in the set of all normal states of $\mathfrak{J}$ (cf. \cite{BarFri90}). In case that $\mathfrak{J}$ is a von Neumann algebra regarded as a JBW$^*$-algebra, the strong$^*$ topology defined above is precisely the usual C$^*$-algebra strong$^*$ topology (see \cite[\S 3]{BarFri90}). The strong$^*$ topology is stronger than the weak$^*$ topology of $\mathfrak{J}$ (see \cite[Theorem 3.2]{BarFri90}). When $\mathfrak{J}$ is regarded as an element in the strictly wider class of JBW$^*$-triples, the strong$^*$-topology, as JBW$^*$-triple, as defined in \cite{BarFri90}, agrees with the strong$^*$-topology we have just defined (see \cite[Proposition 3]{RodPa91} or \cite[\S 4]{PeRo2001}). It is further known that the strong$^*$-topology of $\mathfrak{J}$ is compatible with the duality $(\mathfrak{J},\mathfrak{J}_*)$ (cf. \cite[Corollary 9]{PeRo2001}). Consequently, a convex subset of $\mathfrak{J}$ has the same closure with respect to the strong$^*$ and the weak$^*$ topologies. Moreover, a linear map between JBW$^*$-algebras is strong$^*$-continuous if, and only if, it is weak$^*$-continuous (compare \cite[Corollary 3]{RodPa91} or \cite[comments in page 621]{PeRo2001}). A consequence of this fact implies that every elementary operator on a JBW$^*$-algebra is strong$^*$-continuous.\smallskip 

One additional property of the strong$^*$-topology remains to be recalled. The so-called \emph{Kaplansky density theorem} affirms that if $\mathfrak{B}$ is a weak$^*$-dense JB$^*$-subalgebra of a JBW$^*$-algebra $\mathfrak{J}$, then the closed unit ball of $\mathfrak{B}$ is strong$^*$-dense in the closed unit ball of $\mathfrak{J}$ (cf. \cite[Corollary 3.3]{BarFri90} or \cite[Proposition 2.4]{AlfsenShultz2003}).\smallskip

After gathering the previous basic properties of the strong$^*$-topology, we can now establish the following generalization of \Cref{lemma:elem_homo}.   

\begin{lemma}{\label{lemma:elem_homo metric with norm}} Let $\mathfrak{J}$ be a JB$^*$-algebra,  $\mathfrak{M}$ a JBW$^*$-algebra, and $\pi\colon \mathfrak{J} \rightarrow \mathfrak{M}$ be a Jordan $^*$-homomorphism with weak$^*$-dense image. Then for each $\mathcal{E}\in\element(\mathfrak{J})$, there exists a unique  $\widetilde{\mathcal{E}}\in\element(\mathfrak{M})$ such that $\pi\mathcal{E}= \widetilde{\mathcal{E}} \pi$ and $\|\widetilde{\mathcal{E}}\| \leq \|\mathcal{E}\|$. 
\end{lemma}

\begin{proof} The existence and uniqueness are guaranteed by previous \Cref{lemma:elem_homo}. To prove the statement concerning the norm of $\widetilde{\mathcal{E}}$, let us fix $z$ in the closed unit ball of $\mathfrak{M}$ and an arbitrary $\varepsilon >0$. We begin by observing that, by Kaplansky density theorem \cite[Corollary 3.3]{BarFri90}, we can find a net $(\pi(x_j))_j$ in the closed unit ball of $\pi (\mathfrak{J})$ converging to $z$ in the strong$^*$-topology of $\mathfrak{M}$, and thus $(\widetilde{\mathcal{E}} \pi (x_j))_j\to \widetilde{\mathcal{E}} (z)$ in the strong$^*$-topology, and in the weak$^*$-topology of $\mathfrak{M}$.  Since $\ker (\pi)$ is a norm-closed Jordan ideal of $\mathfrak{J},$ and the mapping $\widehat{\pi}: \mathfrak{J}/\ker(\pi) \to \mathfrak{M},$ $\widehat{\pi}(x+\ker(\pi)):= \pi (x)$ is an isometric Jordan $^*$-monomorphism for all $j$, we can find $k_j\in \ker(\pi)$ such that the inequality $\|x_j + k_j \| \leq (1+\varepsilon ) \|\pi(x_j)\| \leq (1+\varepsilon )$ holds. Finally, since $$\widetilde{\mathcal{E}} (z) = w^{*}\hbox{-}\lim_j \widetilde{\mathcal{E}} \pi (x_j+ k_j) = w^{*}\hbox{-}\lim_j \pi {\mathcal{E}} (x_j+ k_j), $$ with $\| \pi {\mathcal{E}} (x_j+ k_j) \|\leq  (1+\varepsilon ) \|\mathcal{E}\|, $ it follows that $\| \widetilde{\mathcal{E}} (z)\| \leq (1+\varepsilon ) \|\mathcal{E}\|$. 
\end{proof}

\subsection{Linear dependence in prime \texorpdfstring{C$^*$}{C*}-algebras}\ \\

We devote this subsection to survey some results on linear independence of elements in a C$^*$-algebra, and stating several variants to be applied in our Jordan setting.\smallskip    

In this section we shall employ the extended centroid of a prime C$^*$-algebra regarded as a prime ring, in several of our arguments. For this purpose, we recall some basic results. Let ${R}$ be a (possibly non-unital associative) semiprime ring. Following the standard notation (see \cite[\S 2.1]{AraMathieuBook}), we say that an \emph{essentially defined double centraliser} on $R$ is a triple $(f, g,I)$, where $f$ and $g$ are mappings from an essential ring  two-sided ideal $I$ of $R$ into $R,$ $f$ is a left $R$-module homomorphism, $g$ is a right $R$-module homomorphism,
and $f(x) y = x g(y)$ for all $x, y \in I$. Two essentially defined double centralisers $(f_1, g_1,I_1)$ and $(f_2, g_2,I_2)$ are equivalent if $f_1$ and $f_2$ agree on their common domain $I_1\cap  I_2$ --note that $g_1$ and $g_2$ agree on $I_1 \cap I_2$ too. The \emph{symmetric ring of quotients of $R$}, $Q_s(R)$, is defined as the set of equivalence classes of essentially defined double centralisers on $R$. Addition and multiplication operations on $Q_s(R)$ are defined by pointwise addition and composition on the appropriate domains, that is, $$[(f_1, g_1,I_1)]+ [(f_2, g_2,I_2)] := [((f_1+f_2)|_{I_1\cap I_2}, (g_1+g_2)|_{I_1\cap I_2},{I_1\cap I_2})], $$ and  
$$[(f_1, g_1,I_1)] [(f_2, g_2,I_2)] := [(f_2 f_1)|_{I_1 I_2}, (g_1 g_2)|_{I_1 I_2},{I_1 I_2})].$$ Endowed with these operations, $Q_s(R)$ becomes a semiprime ring with identity $[(id, id, R)]$. The ring $R$ embeds into $Q_s(R)$ via the mapping $a\hookrightarrow [(R_a,L_a,R)]$ ($a\in R$). The \emph{extended centroid of $R$}
is defined as the centre of $Q_s(R)$, and it is usually denoted by $C(R)$ \cite[Definition 2.1.1 and Definition 2.1.4]{AraMathieuBook}. The extended centroid of $R$ can be identified with the ring of equivalence classes of essentially defined bimodule homomorphisms on $R$ (cf. comments after Definition 2.1.4 in  \cite{AraMathieuBook}).\smallskip

For later purposes, we present an argument to check that the extended centroid of a prime C$^*$-algebra regarded as a ring is the complex field. 

\begin{remark}\label{r extended centorid of a prime C*-algebra} Let $A$ be a prime C$^*$-algebra, regarded as a prime ring with respect to its natural product. Let $I$ be a non-zero (not necessarily closed) two-sided ring ideal of $A$,
and let $f\colon I\to A$ be a ring $A$-bimodule homomorphism, that is, an additive map such that
\[
f(ax)=f(a)x, \ f(xa)=xf(a), \ \quad\forall a\in I, \ \forall x\in A.
\]
Clearly, the linear span, $\widetilde{I},$ of $I$ is a two-sided algebra ideal of $A$. We define $\widetilde{f}\colon\widetilde{I}\to A$ by
$\widetilde{f}\Big(\sum\alpha_j a_j\Big):=\sum\alpha_jf(a_j)$.
If $\sum \alpha_ja_j=\sum \beta_k b_k$, then $ x \Big(\sum\alpha_j a_j\Big) = x \Big(\sum\beta_k b_k\Big)$, for every $x\in A$. Now, by applying that $f$ is a ring $A$-bimodule homomorphism we get
\[
\begin{split}
x\Big(\sum\alpha_jf(a_j)-\sum\beta_kf(b_k)\Big)
&=
\sum x\alpha_jf(a_j)-\sum x\beta_kf(b_k)\\
&=
f\Big(\sum x\alpha_ja_j\Big)-
f\Big(\sum x\beta_kb_k\Big) \\
& =
f\Big(x \Big(\sum\alpha_j a_j\Big) \Big)-
f\Big(x \Big(\sum\beta_k b_k\Big) \Big) =0,
\end{split}
\] which, by the primeness of $A$ and the arbitrariness of $x\in A$, gives
$\sum\alpha_jf(a_j)=\sum\beta_k f(b_k)$,
and thus $\widetilde{f}$ is a well-defined $A$-bimodule homomorphism. 
Since $A$ is prime, by \cite[Proposition 2.5]{Mathieu1989}, there exists $\alpha\in\mathbb{C}$ such that $\widetilde{f}(x)=\alpha x$ for each $x\in\widetilde{I}$, and hence
\[
f(a)=\alpha a,\hbox{ for all }a\in I.
\] This shows, in particular, that $C(A) = \mathbb{C}$ when $A$ is regarded as a prime ring. 
\end{remark}

We continue by recalling another useful tool. Let $n\in\mathbb{N}$ with $n\ge 2$. The noncommutative polynomial $c_n$ in the non-necessarily commutative indeterminates $\xi_1,\ldots,\xi_n,\eta_1,\ldots,\eta_{n-1}$ defined by
$$ c_n(\xi_1,\ldots,\xi_n,\eta_1,\ldots,\eta_{n-1})=
\sum_{\sigma\in S_n}\!\!\!\text{sign}(\sigma)
\xi_{\sigma(1)}\eta_1\xi_{\sigma(2)}\eta_2\cdots\xi_{\sigma(n-1)}\eta_{n-1}\xi_{\sigma(n)},$$
is called the \emph{$n$th Capelli polynomial} (cf. \cite[Definition 6.11]{Br}). Here, $S_n$ stands for the symmetric group of order $n$.\smallskip

Our next result shows how Capelli's polynomials can be employed to determine the linear independence of a set of vectors in a prime C$^*$-algebra.  

\begin{lemma}\label{independencia1} Let $A$ be a prime C$^*$-algebra, and
let $a\in A^n$, $n\ge 2$.
Then $a_1,\ldots,a_n$ are linearly dependent if and only if
\[
c_n (a_1, \ldots, a_n, x_1, \ldots, x_{n-1}) =0, \ \hbox{ for all } x_1,\ldots,x_{n-1}\hbox{ in }A. \]
\end{lemma}

\begin{proof} The desired conclusion follows from  \cite[Theorem 7.45]{Br} by just recalling the natural embedding of $A$ inside $Q_s (A)$ together with the fact that $C(A) = \mathbb{C}$ (see Remark~\ref{r extended centorid of a prime C*-algebra}).
\end{proof}

We now state some results on linear identities in the line of \cite{Martindale31969, Mathieu1989,LL}, which will be required in our arguments. 

\begin{lemma}\label{independencia2} Let $A$ be a C$^*$--algebra, and
let $a=(a_1,\ldots, a_n)$, $b=(b_1,\ldots, b_n)$ in $A^n$.
Assume that $a_1,\ldots,a_p$ are linearly independent. Then the following statements hold.
\begin{enumerate}[{\rm(i)}]
\item If $A$ is prime and $\displaystyle \sum_{j=1}^n a_j x b_j=0,$ for all $x\in A,$ then each $b_j$, with $j\in\{1,\ldots,p\}$, is a linear combination of $b_{p+1},\ldots,b_n$ (in case $p=n$, this should be understood as $b_j=0,$ for each $j\in\{1,\ldots,n\}$).
\item If $A$ is prime, $\tau$ is a ring involution on $A$, $a_1,\ldots,a_p\in H(A,\tau) = \{z\in A : \tau (z) =z\}$ and 
\[
\sum_{j=1}^na_jxb_j=0, \hbox{ for all } x\in H(A,\tau),
\] then each $b_j$, with $j\in\{1,\ldots,p\}$, is a linear combination of $b_{p+1},\ldots,b_n$ (in case $p=n$, this should be understood as $b_j=0,$ for each $j\in\{1,\ldots,n\}$).
\item If $\tau$ is a $^*$-anti-automorphism of period-$2$ on $A$, and $A$ decomposes in the form $A = N\oplus^{\infty} \tau(N)$, where $N$ is a prime C$^*$-algebra and a closed ideal of $A$, $a_1,\ldots,a_p\in H(A,\tau)$ and 
\[
\sum_{j=1}^na_jxb_j=0, \hbox{ for all } x\in H(A,\tau),
\] then each $b_j$, with $j\in\{1,\ldots,p\}$, is a linear combination of $b_{p+1},\ldots,b_n$ (in case $p=n$, this should be understood as $b_j=0,$ for each $j\in\{1,\ldots,n\}$).
\end{enumerate}
\end{lemma}

\begin{proof}
(i) By renumbering, we may suppose that $\{a_1,\ldots,a_m\}$, with $p\le m\le n$,
is a maximal linearly independent subset of $\{a_1,\ldots,a_n\}$.
Then we apply \cite[Theorem 4.1]{Mathieu1989} to the opposite algebra $A^{op}$ of $A$, to obtain that each $b_j$, with $j\in\{1,\ldots,m\}$, is a linear combination of $b_{m+1},\ldots,b_n$.\smallskip

(ii) The previous Remark~\ref{r extended centorid of a prime C*-algebra} shows that the \emph{extended centroid} of $A$ regarded as a ring is isomorphic to $\mathbb{C}$, then the desired result is a consequence of \cite[Lemma 3]{LL}.\smallskip

(iii) As in the proof of \Cref{l HWtau factor implies W factor}, $H(A, \tau)=\{ a + \tau(a) : a\in N\}$ admits a structure of C$^*$-subalgebra of $N\oplus^{\infty} \tau(N)^{op},$ making it C$^*$-isomorphic to $N$. Let $\pi$ denote the canonical projection of $A$ onto $N$, which is a C$^*$-homomorphism.  Since $\displaystyle \sum_{j=1}^na_jxb_j=0,$ for all $x\in H(A,\tau),$ we conclude that  $\displaystyle \sum_{j=1}^n \pi(a_j) \pi(x) \pi (b_j)=0,$  for all $x\in H(A,\tau),$ and $N$ is a prime C$^*$-algebra, we deduce from (i) that each $\pi (b_j)$, with $j\in\{1,\ldots,p\}$, is a linear combination of $\pi (b_{p+1}),\ldots,\pi(b_n)$, which in turn gives that each $b_j = \pi(b_j) + \tau\pi(b_j)$, with $j\in\{1,\ldots,p\}$, is a linear combination of $b_{p+1}=\pi(b_{p+1}) + \tau\pi(b_{p+1}),\ldots,b_n=\pi(b_n) + \tau\pi(b_n)$.
\end{proof}

\section{Associating traces on \texorpdfstring{JW$^*$}{JW*}-algebras without type \texorpdfstring{$I_1$}{I1} and \texorpdfstring{$I_2$}{I2} summands}\label{sec: bilinear with associating trace on JW*-algebras without spin}

The main goal of this section is to determine the general form of all associating traces on a JW$^*$-algebra $\mathfrak{J}$ without associative and spin part. 
The task will be completed in \Cref{prop:JBWspecial}, however the arguments will require to develop certain tools and techniques inspired on arguments developed in different branches. We have already justified in the introduction how the existence of associative summands makes hopeless our conclusions.\smallskip

In a first stage we focus now on the uniqueness of the ``standard form'' of a symmetric bilinear map $B$ on a JBW$^*$-algebra $\mathfrak{J}$ whose trace is associating.

\begin{proposition}\label{prop: uniq_sol} Let $\mathfrak{J}$ and $\mathfrak{M}$ be JBW$^*$-algebras, $B : \mathfrak{J} \times \mathfrak{J} \rightarrow \mathfrak{M}$ a symmetric bilinear mapping, and $\pi: \mathfrak{J} \rightarrow \mathfrak{M}$ a Jordan homomorphism having weak$^*$-dense image. Assume that there exist an element $w$ in $\mathfrak{J}$ and elementary operators $\mathcal{E}_0, \mathcal{E}_1,$ and  $\mathcal{E}_2$ on $\mathfrak{J}$ such that $\mathcal{E}_i(w^j) = \delta_{ij} \unit,$ for all $ i,j \in \{0,1,2\}$. Suppose additionally that $B$ admits the following standard representation 
$$ B(x,x) = \lambda\circ \pi(x)^2  + \mu(x)\circ \pi(x) + \nu(x,x),$$ for every $x \in \mathfrak{J}$, where $\lambda \in Z (\mathfrak{M})$, $\mu : \mathfrak{J} \rightarrow Z (\mathfrak{M})$ is a linear mapping, and $\nu: \mathfrak{J} \times \mathfrak{J}  \rightarrow Z (\mathfrak{M})$ is a symmetric bilinear mapping. Then this representation is unique and satisfies the following statements:
\begin{enumerate}[$\bullet$]
\item $\lambda = \widehat{\mathcal{E}_2}(B(w,w))$,
\item $\mu(x) = 2\widehat{\mathcal{E}_1}(B(w,x)) -2\widehat{\mathcal{E}_2}(B(w,w))\circ \widehat{\mathcal{E}_1}(\pi(w\circ x)) -\widehat{\mathcal{E}_1}(B(w,w))\circ\widehat{\mathcal{E}_1}(\pi(x)),$ for every\hyphenation{every} $x \in \mathfrak{J}$,
\item $\nu(x,x):= \widehat{\mathcal{E}_0}(B(x,x))- \lambda \circ\widehat{\mathcal{E}_0}(\pi(x)^2) - \mu(x)\circ\widehat{\mathcal{E}_0}(\pi(x)),$ for every $x \in \mathfrak{J}$, 
\end{enumerate}
where each $\widehat{\mathcal{E}_i}$ is the (unique) elementary operators in $\element(\mathfrak{M})$ satisfying $\widehat{\mathcal{E}_i}\pi = \pi\mathcal{E}_i,$ for every $i \in \{0,1,2\}$.
\end{proposition}

\begin{proof} By \Cref{lemma:elem_homo}(ii) there exist unique elementary operators $\widehat{\mathcal{E}_i} \in \element(\mathfrak{M})$ such that $  \widehat{\mathcal{E}_i}\pi = \pi\mathcal{E}_i$, $i \in \{0,1,2\}$. Thus, they satisfy 
    \begin{equation*}
     \widehat{\mathcal{E}_i}(\pi(w)^j) = \widehat{\mathcal{E}_i}(\pi(w^j)) =\pi(\mathcal{E}_i(w^j)) = \delta_{ij} \unit \hbox{ for all $i,j \in \{0,1,2\}$}.
    \end{equation*} Since the elementary operators $\widehat{\mathcal{E}_i}$ are $\mathbb{C}$-linear and actually $Z (\mathfrak{M})$-linear (cf. Remark~\ref{rem:elementary Z-linear}), it follows from the assumptions that 
	\begin{equation*}
		\widehat{\mathcal{E}_2}(B(w,w)) = \lambda\circ \widehat{\mathcal{E}_2}(\pi(w)^2)+ \mu(w)\circ \widehat{\mathcal{E}_2}(\pi(w)) + \nu(w,w)\circ \widehat{\mathcal{E}_2}(\mathbf{1}) = \lambda \circ\mathbf{1} = \lambda.
	\end{equation*}
	Similarly, by applying $\widehat{\mathcal{E}_1}$ to $B(w,w)$ we obtain from the hypotheses that 
        \begin{equation*}\begin{aligned}
        		\widehat{\mathcal{E}_1}(B(w,w)) &= \lambda\circ \widehat{\mathcal{E}_1}(\pi(w)^2)+ \mu(w)\circ \widehat{\mathcal{E}_1}(\pi(w)) + \nu(w,w)\circ \widehat{\mathcal{E}_1}(\mathbf{1}) \\ 
        		&= \mu(w) \circ\mathbf{1} = \mu(w).
        	\end{aligned}
	\end{equation*}
To deduce the second identity in the conclusions, let us take an  arbitrary $x \in \mathfrak{J}$. Then, by hypotheses, we have $$B(w+x,w+x) = 2B(w,x) + B(x,x) + B(w,w).$$ Making some computations we obtain 
	$$  2 B(w,x)=2 \lambda\circ \pi(w\circ x)  + \mu(w)\circ \pi(x) + \mu(x)\circ \pi(w) + \nu(w,x) + \nu(x,w).$$
	Now, by applying $\widehat{\mathcal{E}_1}$ on both sides of the previous identity we get 
	$$\begin{aligned}
		2\widehat{\mathcal{E}_1}(B(w,x)) & = 2\lambda \circ \widehat{\mathcal{E}_1}(\pi(w\circ x)) + \mu(w)\circ\widehat{\mathcal{E}_1}(\pi(x)) + \mu(x)\circ\widehat{\mathcal{E}_1}(\pi(w))\\ &\hspace*{0.3cm}+ (\nu(w,x) + \nu(x,w))\circ\widehat{\mathcal{E}_1}(\mathbf{1})\\
		& = 2\lambda \circ \widehat{\mathcal{E}_1}(\pi(w\circ x)) + \mu(w)\circ\widehat{\mathcal{E}_1}(\pi(x)) + \mu(x),
	\end{aligned}$$	
	which implies that
	\begin{equation*}
		\mu(x) = 2\widehat{\mathcal{E}_1}(B(w,x)) - 2\lambda\circ \widehat{\mathcal{E}_1}(\pi(w\circ x))- \mu(w)\circ\widehat{\mathcal{E}_1}(\pi(x)), \ \forall x\in \mathfrak{J}.
	\end{equation*}
	Having in mind that $\lambda = \widehat{\mathcal{E}_2}(B(w,w))$ and $\mu(w) = \widehat{\mathcal{E}_1}(B(w,w)) $ we are led to 
	$$\mu(x) = 2 \widehat{\mathcal{E}_1}(B(w,x)) - 2 \widehat{\mathcal{E}_2}(B(w,w))\circ\widehat{\mathcal{E}_1}(\pi(w\circ x)) -\widehat{\mathcal{E}_1}(B(w,w))\circ \widehat{\mathcal{E}_1}(\pi(x))$$ for every $x \in \mathfrak{J},$ as we claimed. Finally, by applying $\widehat{\mathcal{E}_0}$ to $B(x,x)$ for $x \in \mathfrak{J}$ we obtain the last desired identity for $ \nu(x,x)$. 
\end{proof}

The arguments in the proof of the previous proposition are also valid to get the following conclusion.

\begin{proposition}\label{prop: uniq_sol spin} Let $\mathfrak{J}$ and $\mathfrak{M}$ be JBW$^*$-algebras, $B : \mathfrak{J} \times \mathfrak{J} \rightarrow \mathfrak{M}$ a symmetric bilinear mapping, and $\pi: \mathfrak{J} \rightarrow \mathfrak{M}$ a Jordan homomorphism having weak$^*$-dense image. Assume that there exist an element $w$ in $\mathfrak{J},$ and elementary operators $\mathcal{E}_0$ and  $\mathcal{E}_1$ on $\mathfrak{J}$ such that $\mathcal{E}_i(w^j) = \delta_{ij} \unit,$ for all $ i,j \in \{0,1\}$. Suppose additionally that $B$ admits the following standard representation 
$$ B(x,x) =  \mu(x)\circ \pi(x) + \nu(x,x),$$
     for every $x \in \mathfrak{J}$, where $\mu : \mathfrak{J} \rightarrow Z (\mathfrak{M})$ is a linear mapping, and $\nu: \mathfrak{J} \times \mathfrak{J}  \rightarrow Z (\mathfrak{M})$ is a symmetric bilinear mapping. Then this representation is unique and satisfies the following statements:
\begin{enumerate}[$\bullet$]
\item $\mu(x) = 2\widehat{\mathcal{E}_1}(B(w,x))  -\widehat{\mathcal{E}_1}(B(w,w))\circ\widehat{\mathcal{E}_1}(\pi(x)),$ for every\hyphenation{every} $x \in \mathfrak{J}$, 
\item $\nu(x,x):= \widehat{\mathcal{E}_0}(B(x,x)) - \mu(x)\circ\widehat{\mathcal{E}_0}(\pi(x)),$ for every $x \in \mathfrak{J}$, 
\end{enumerate}
where each $\widehat{\mathcal{E}_i}$ is the (unique) elementary operators in $\element(\mathfrak{M})$ satisfying $\widehat{\mathcal{E}_i}\pi = \pi\mathcal{E}_i,$ $i \in \{0,1,2\}$.
\end{proposition}

The next lemma establishes a kind of selection principle for mappings whose trace at every point is a linear combination of a Jordan representation at that point and the unit element, that is, the ``properly quadratic'' monomial part vanishes pointwise. The result is also the first step towards the existence of the mappings $\mu(\cdot)$ and $\nu (\cdot, \cdot)$ appearing in the standard form of a symmetric bilinear mapping whose trace is associating. 

\begin{lemma}\label{lemma:lin_exist}
	Let $\mathfrak{J}, \mathfrak{M}$ be JBW$^*$-algebras where $\mathfrak{J}$ has no direct summands of type $I_1$ and $I_2$, $\pi: \mathfrak{J} \rightarrow \mathfrak{M}$ a Jordan homomorphism having weak$^*$-dense image, and  $B : \mathfrak{J} \times \mathfrak{J} \rightarrow \mathfrak{M}$ a symmetric bilinear map. Assume that there exist an element $w$ in $\mathfrak{J}$ and elementary operators $\mathcal{E}_0, \mathcal{E}_1, \mathcal{E}_2 \in \element(\mathfrak{J})$  such that $\mathcal{E}_i(w^j) = \delta_{ij}\mathbf{1},$ for all $i,j \in {0,1,2}$. Suppose additionally that  for each $x \in \mathfrak{J}$ there exist $\alpha_x, \beta_x \in \mathbb{C}$ such that  $$B(x,x) = \alpha_x \pi(x) + \beta_x \mathbf{1}.$$  Then there exist a linear map $\mu: \mathfrak{J} \rightarrow \mathbb{C}$ and a symmetric bilinear map $\nu: \mathfrak{J} \times \mathfrak{J} \rightarrow \mathbb{C}$ such that 
    $$ B(x,x) = \mu(x) \pi(x) + \nu(x,x) \mathbf{1},$$
    for every $x \in \mathfrak{J}$.
\end{lemma}

\begin{proof} By \Cref{lemma:elem_homo}(ii) we can find unique elementary operators  $\widehat{\mathcal{E}}_{j}\in \element (\mathfrak{M})$ satisfying $\pi\mathcal{E}_i = \widehat{\mathcal{E}_i}\pi,$ for all $i \in \{0,1,2\}$. In particular, $\widehat{\mathcal{E}}_{j} (\pi(w)^i) = \delta_{ij} \unit,$ $\forall i,j\in \{0,1,2\}$. \smallskip

Let us fix and arbitrary $x \in \mathfrak{J}$. The identity \begin{equation*}
		B(x+w, x+w) + B(x-w,x-w) = 2B(x,x)+ 2 B(w,w),
	\end{equation*} holds by hypotheses. 
    Thus, we have 
    $$\begin{aligned}
    	(\alpha_{x+w} + \alpha_{x-w} -2\alpha_x) \pi(x) &+ (\alpha_{x+w} - \alpha_{x-w} - 2\alpha_w)\pi(w)  \\
       &+ (\beta_{x +w} + \beta_{x-w} -\beta_x - 2\beta_t )\mathbf{1} = 0.
    \end{aligned}$$
     
If the vectors $\{\pi(x), \mathbf{1}, \pi(w)\}$ are linearly independent we deduce that   
     $$\alpha_{x+w} + \alpha_{x-w} = 2\alpha_x,\hbox{ and } 
     	\alpha_{x+w} - \alpha_{x-w} = 2\alpha_w,$$
     which implies that
     \begin{equation}
     	\alpha_{x+w} = \alpha_x + \alpha_w,
     	\label{eq:lin_indep}
     \end{equation}
     for every $x \in \mathfrak{J}$ such that $\pi(x)$ is not a linear combination of $\mathbf{1}$ and $\pi(w)$. Note that by the properties of the elementary operators $\widehat{\mathcal{E}}_{j}\in \element (\mathfrak{M}),$ the vectors $\mathbf{1}$ and $\pi(w)$ are linearly independent. \smallskip

Let us now define the maps $\mu: \mathfrak{J} \rightarrow \mathfrak{M}$, $\nu: \mathfrak{J} \times \mathfrak{J} \rightarrow \mathfrak{M}$ by
      $$\mu(x):=2\widehat{\mathcal{E}}_1(B(x,w)) -\widehat{\mathcal{E}_1}(B(w,w))\circ \widehat{\mathcal{E}}_1(\pi(x)), \hbox{ and } \nu(x,x):= B(x,x) - \mu(x)\circ \pi(x).$$ 
      By definition, $\mu$ is a linear map while $\nu$ is given by a symmetric bilinear map. It only remains to prove that both maps take values in $\mathbb{C}\mathbf{1}$. Indeed,\smallskip
      \begin{enumerate}[$(1)$]
\item If $\pi(x) \notin \hbox{Span}\{\mathbf{1}, \pi(w)\}$, it is clear from \eqref{eq:lin_indep} that $$B(x+w,x+w) = (\alpha_x + \alpha_w)\pi(x+w) + \beta_{x+w}\unit,$$  and thus, by combining the properties of $B$ and  $\widehat{\mathcal{E}}_1$ with the previous identity we get 
       $$\begin{aligned}
           2 \widehat{\mathcal{E}}_1(B(x,w)) &= \alpha_x\widehat{\mathcal{E}}_1(\pi(w)) + \alpha_w \widehat{\mathcal{E}}_1(\pi(x)) + (\beta_{x+w}-\beta_w-\beta_x)\widehat{\mathcal{E}}_1(\unit)\\
            &= \alpha_x\unit + \alpha_w \widehat{\mathcal{E}}_1(\pi(x)) = \alpha_x\unit + \widehat{\mathcal{E}_1}(B(w,w))\circ \widehat{\mathcal{E}}_1(\pi(x))  ,
       \end{aligned}$$
        and thus  $\mu(x) = 2 \widehat{\mathcal{E}}_1(B(x,w))- \widehat{\mathcal{E}_1}(B(w,w))\circ \widehat{\mathcal{E}}_1(\pi(x))  = \alpha_x\unit \in \mathbb{C}\mathbf{1}$. Observe that $\nu(x,x)=B(x,x) - \mu(x)\circ \pi(x) =B(x,x) - \alpha_x \circ \pi(x) = \beta_x \unit.$
\item If $\pi(x) \in \hbox{Span} \{\mathbf{1}, \pi(w)\}$, by the properties of the mappings $\widehat{\mathcal{E}}_i$, \linebreak $\pi(w^2)  \notin \hbox{Span} \{\mathbf{1}, \pi(w)\}$, and we can thus conclude that for every $\varepsilon > 0,$ the element $\pi(x +\varepsilon w^2)$ is not in $\hbox{Span} \{\mathbf{1}, \pi(w)\}$. By applying \eqref{eq:lin_indep} we obtain
      	\begin{equation*}
      		B(x+\varepsilon w^2,x+\varepsilon w^2 ) - \mu(x+\varepsilon w^2)\pi(x+\varepsilon w^2) \in \mathbb{C}\mathbf{1}.
      	\end{equation*}
       Note that, by what we have just proved in $(1)$, $\mu(x +\varepsilon w^2) = \alpha_{x + \varepsilon w^2} \unit.$ We also know that $\mu$ is linear. Thus 
       $$\begin{aligned}
       	B(x,x) + 2B(x,\varepsilon w^2) + \varepsilon^2B(w^2,w^2) &- (\mu(x) + \varepsilon\mu(w^2))\circ\pi(x+\varepsilon w^2) \in \mathbb{C}\mathbf{1},\\
        \mu (x + \varepsilon w^2) = & \mu(x) + \varepsilon \mu (w^2)\in \mathbb{C}\unit.
       \end{aligned}$$
      	Taking $\varepsilon \rightarrow 0$ we get 
      	\begin{equation}
      		\nu(x,x)= B(x,x) - \mu(x)\circ \pi(x) \in \mathbb{C}\mathbf{1}, \hbox{ and } \mu (x)\in \mathbb{C}\unit,
      	\end{equation}
       for every $x \in \mathfrak{J}$.
      \end{enumerate}\end{proof}

\begin{remark}\label{r identity Bresar} We gather in this remark some useful identities inspired by previous contributions in \cite[Theorem 2]{Bres1992}, \cite[pages 2891-2892]{LL} and \cite[Lemma 3.8.]{BreEreVil}. \smallskip

Let $\mathfrak{X}$ be a linear space, $A$ an associative algebra, $\pi : \mathfrak{X}\to A$ a linear mapping, and $B : \mathfrak{X} \times \mathfrak{X} \rightarrow A$ a symmetric bilinear map. Suppose $[B(x,x),\pi(x)] = 0,$ for every $x \in \mathfrak{X}$, where the brackets denote the usual Lie product on $A$. Then  \begin{equation}\label{eq new remark 1206} 2[B(x,y),\pi(y)] = -[B(y,y),\pi(x)], \hbox{ for all } x,y\in \mathfrak{X}.
\end{equation}
Namely, the desired identity is just a clear consequence of the identity $0= [B(x+ \alpha y,x+\alpha y),\pi(x+\alpha y)],$ for all $x,y\in \mathfrak{X},$ $\alpha \in \mathbb{R}$, by just expanding the polynomial on $\alpha$ via the properties of $B$.\smallskip

Similar arguments building upon identity \eqref{eq new remark 1206} lead to \begin{equation}\label{eq new remark 2 1206}  [B(x,y),\pi(z)] + [B(x,z),\pi(y)] + [B(y,z), \pi(x)] = 0, \hbox{ for all }   x,y\in \mathfrak{X}.
\end{equation}

Suppose additionally that $\mathfrak{X}$ is a Jordan algebra and $\pi$ is a Jordan homomorphism. We combine the associative structure of $A$ with the previous identites to deduce the following:
$$\begin{aligned}
\relax[B(x^2,x^2),\pi(y)] &= -2 [B(x^2,y),\pi(x^2)] \\
&= -2 [B(x^2,y),\pi(x)] \pi(x) - 2 \pi(x) [B(x^2,y),\pi(x)]   \\
&= 2 \left( [B(x^2,x),\pi(y)] + [B(y,x),\pi(x^2)] \right) \pi (x) \\ 
&\ \ \ +  2 \pi (x) \left( [B(x^2,x),\pi(y)] + [B(y,x),\pi(x^2)] \right) \\
= &2 \left( [B(x^2,x),\pi(y)] + [B(y,x),\pi(x)] \pi (x) + \pi (x) [B(y,x),\pi(x)] \right) \pi (x) \\ 
 +&  2 \pi (x) \left( [B(x^2,x),\pi(y)] + [B(y,x),\pi(x)] \pi (x) + \pi (x) [B(y,x),\pi(x)] \right) \\
&= 2 \pi (x^2)  [B(y,x),\pi(x)] + 4 \pi (x) [B(y,x),\pi(x)] \pi (x) \\
& \ \ + 2 [B(y,x),\pi(x)]  \pi (x^2)+  2 \pi (x)[B(x^2,x),\pi(y)] \\
&\ \ +  2 [B(x^2,x),\pi(y)]  \pi (x) \\
& = - \pi (x^2)  [B(x,x),\pi(y)]  - 2 \pi (x) [B(x,x),\pi(y)] \pi (x) \\
&\ \ - [B(x,x),\pi(y)] \pi (x^2) +  2 \pi (x)[B(x^2,x),\pi(y)]  \\ 
&\ \ +  2 [B(x^2,x),\pi(y)]  \pi (x),
\end{aligned}$$ where in the first equality we applied \eqref{eq new remark 1206}, while in the third and sixth equalities we applied \eqref{eq new remark 2 1206}. If we rewrite the previous identity in terms of the associative product we get \begin{equation}\label{equation crazy LeeLee}
\begin{aligned} &\pi (\unit) \pi (y) \left(2 B(x^2,x) \pi(x) - B(x^2,x^2) - B(x,x) \pi(x^2)  \right)\\  &+ \pi (x) \pi (y) \left(2 B(x^2,x) - 2 B(x,x)\pi(x)  \right) - \pi (x^2)  \pi (y) B(x,x) \\
& +  \left(B(x^2,x^2) + \pi (x^2) B(x,x) - 2 \pi (x) B(x,x^2)\right) \pi (y) \pi(\unit) \\ 
&+\left(2 \pi (x) B(x,x) -2 B(x,x^2) \right) \pi (y) \pi (x) + B(x,x) \pi (y) \pi (x^2) = 0,
\end{aligned}   
\end{equation} for all $x,y\in \mathfrak{X}.$\smallskip

We finally observe that the arguments above are valid to get the following conclusion in the Jordan setting.  Suppose $\mathfrak{X}$ is a linear space, $\pi : \mathfrak{X}\to \mathfrak{J}$ is a linear mapping, where $\mathfrak{J}$ is a Jordan algebra, and $B : \mathfrak{X} \times \mathfrak{X} \rightarrow \mathfrak{J}$ is a symmetric bilinear map such that $[B(x,x),\mathfrak{J},\pi(x)] = 0,$ for every $x \in \mathfrak{X}$. Then  \begin{equation}\label{eq new remark 1206 Jordan} 2[B(x,y),a,\pi(y)] = -[B(y,y),a,\pi(x)], \hbox{ for all } x,y\in \mathfrak{X}, a\in\mathfrak{J}.
\end{equation} 

\end{remark}

We deal next with general bilinear maps with associating trace. 

\begin{lemma}{\label{lemma:lambda_indep}}
	Let $\mathfrak{J}$ be a JBW$^*$-algebra without direct summands of type $I_1$ and $I_2$, $A$ a prime von Neumann algebra, $\pi: \mathfrak{J} \rightarrow A$ a unital Jordan homomorphism with weak$^*$-dense image, and $B : \mathfrak{J} \times \mathfrak{J} \rightarrow A$ a symmetric bilinear map. Suppose that $[B(x,x),\pi(x)] = 0,$ for every $x \in \mathfrak{J}$, where the brackets denote the usual Lie product on $A$. Assume, additionally, that for each $x \in \mathfrak{J}$ there exist $\lambda_x, \mu_x, \nu_x \in \mathbb{C}$ such that  
	$$ B(x,x) = \lambda_x\pi(x)^2  + \mu_x \pi(x) + \nu_x\mathbf{1}.$$
    Then there exist a complex value $\lambda \in \mathbb{C}$ such that $B(x,x) - \lambda \pi(x)^2 \in \hbox{Span}\{\mathbf{1}, \pi(x)\}$ for every $x \in \mathfrak{J}$.
\end{lemma}

\begin{proof}
    Let us take $x,y \in \mathfrak{J}$, by the hypotheses we get
    \begin{align*}
        B(x+y,x+y) &= \lambda_{x+y}\pi(x+y)^2 + \mu_{x+y}\pi(x+y) + \nu_{x+y}\mathbf{1}\\
                   &= \lambda_{x+y}(\pi(x)^2 + \pi(y^2) + 2 \pi(x\circ y) ) + \mu_{x+y}(\pi(x) + \pi(y)) +\nu_{x+y}\mathbf{1},
    \end{align*}
    which leads us to the following identity
    \begin{align*}
        2B(x,y) = &(\lambda_{x+y} -\lambda_x)\pi(x)^2 + (\lambda_{x+y} -\lambda_y)\pi(y)^2  \\
                &+ (\mu_{x+y} -\mu_x)\pi(x) + (\mu_{x+y} -\mu_y)\pi(y)\\
                &+(\nu_{x+y} - \nu_x - \nu_y)\mathbf{1} + 2 \lambda_{x+y}\pi(x\circ y).
    \end{align*}
  Next, we apply the commutator $[\cdot, \pi(y)]$  on the previous identity to derive that  
    \begin{align*}
        [2B(x,y), \pi(y)]  = & (\lambda_{x+y} -\lambda_x)[\pi(x)^2,\pi(y)] + (\lambda_{x+y} -\lambda_y)\underbrace{[\pi(y)^2,\pi(y)]}_{= 0}  \\
                             &+(\mu_{x+y} -\mu_x)[\pi(x),\pi(y)] + (\mu_{x+y} -\mu_y)\underbrace{[\pi(y),\pi(y)]}_{=0} \\
                             &+(\nu_{x+y} - \nu_x - \nu_y)\underbrace{[\mathbf{1},\pi(y)]}_{=0} + 2\lambda_{x+y}[\pi(x\circ y), \pi(y)].
    \end{align*}
 By \eqref{eq new remark 1206} we have $2[B(x,y),\pi(y)] = -[B(y,y),\pi(x)]$. To deal with the remaining summands on the right-hand-side of the previous identity, we allude to the associative structure of $A$ to make the following computations:
     \begin{align*}
         [2\pi(x\circ y)], \pi(y)] & = [\pi(x)\pi(y) + 
                                   \pi(y)\pi(x),\pi(y)]\\
                                   &= [\pi(x)\pi(y) ,\pi(y)] + [\pi(y)\pi(x), \pi(y)] \\
                                   & = \pi(x)\pi(y)^2 - \pi(y)\pi(x)\pi(y) + \pi(y)\pi(x)\pi(y) - \pi(y)^2\pi(x)  \\
                                   & = \pi(x)\pi(y)^2 - \pi(y)^2\pi(x) = [\pi(x), \pi(y)^2].
     \end{align*} By combining all the previous observations we get 
     \begin{align*}
         -[B(y,y),\pi(x)] &=  (\lambda_{x+y} -\lambda_x)[\pi(x)^2,\pi(y)] \\ 
         &+ (\mu_{x+y} -\mu_x)[\pi(x),\pi(y)] + \lambda_{x+y}[\pi(x),\pi(y)^2]. 
     \end{align*}
Having in mind that, by hypotheses, $B(y,y) = \lambda_y \pi(y)^2 + \mu_y\pi(y) + \nu_y\mathbf{1}$ for some $\lambda_y, \mu_y,\nu_y \in \mathbb{C},$ we arrive to 
     $$-[B(y,y),\pi(x)] = \lambda_y[\pi(x),\pi(y)^2] + \mu_y[\pi(x),\pi(y)].$$ The last two previous identities together give 
     \begin{equation}{\label{eq: lambda}}
     \begin{aligned}
         (\lambda_{x+y} - \lambda_x)[\pi(x)^2,\pi(y)] &+ (\lambda_{x+y} -\lambda_y)[\pi(x),\pi(y)^2] \\
         & + (\mu_{x+y}-\mu_x-\mu_y)[\pi(x),\pi(y)] = 0,
     \end{aligned}
     \end{equation}
     for every $x, y \in \mathfrak{J}$. \smallskip

Let us take elements $u,v \in \mathfrak{J}$ whose existence is guaranteed by \Cref{cor:3_lib}. By hypotheses we have $$B(u,u) = \lambda_u \pi(u)^2  + \mu_u \pi(u) + \nu_u \mathbf{1}.$$ The just quoted \Cref{cor:3_lib} assures the existence of elementary operators $\mathcal{E}_j \in \element(\mathfrak{J})$ satisfying $\mathcal{E}_j (u^{i}) = \delta_{i,j} \unit,$ for all $i,j\in\{0,1,2\}.$ \Cref{lemma:elem_homo}(ii) implies the existence of elementary operators $\widehat{\mathcal{E}}_j \in \element(A)$ satisfying $\widehat{\mathcal{E}}_j (\pi (u)^{i}) = \delta_{i,j} \unit,$ for all $i,j\in\{0,1,2\}.$ Consequently, the elements $\pi(u)^2,$ $\pi(u)$ and $\mathbf{1}$ are linearly independent. Set $\lambda:= \lambda_u$. We claim that \begin{equation}\label{eq B(x,x)-quadratic part lies in the previous lemma}\hbox{$B(x,x) - \lambda \pi(x)^2 \in \hbox{Span}\{\mathbf{1}, \pi(x)\},$ for every $x \in \mathfrak{J}$.}   
\end{equation} To prove the claim we shall employ the results on linear independence in terms of Capelli's polynomials. We shall distinguish two cases:
\begin{enumerate}[$(1)$]
\item \textit{Case I}: The vectors $\{ [\pi(x)^2,\pi(u)], [\pi(x),\pi(u)^2], [\pi(x),\pi(u)]\}$ are linearly independent for some $x \in \mathfrak{J}$, then we deduce from  \eqref{eq: lambda} that $\lambda_{x+u} - \lambda_x = 0$ and $\lambda_{x+u} - \lambda_u = 0$. This implies that $ \lambda_x = \lambda_u = \lambda$ by definition, and by hypotheses $$ B(x,x) = \lambda \pi(x)^2  + \mu_x \pi(x) + \nu_x\mathbf{1},$$ which proves the claim in this case.\smallskip
         
\item \textit{Case II}:  There exists $x \in \mathfrak{J}$
such that the vectors $[\pi(x)^2,\pi(u)],$ $ [\pi(x),\pi(u)^2],$ and $ [\pi(x),\pi(u)]\}$ are linearly dependent. Consider the third Capelli polynomial $c_3(\cdot,\cdot,\cdot,\cdot,\cdot),$ and define the polynomials $f,g: \mathfrak{J} \rightarrow A$ by 
         \begin{align*}
             f_{a,b}(x) & := c_3\Big(B(x,x)-\lambda \pi(x)^2, \mathbf{1}, \pi(x),a,b\Big),\\
             g_{c,d}(x) & := c_3\Big([\pi(x)^2,\pi(u)], [\pi(x),\pi(u)^2], [\pi(x),\pi(u)], c,d\Big),
         \end{align*}
 where  $a,b,c,d$ are fixed but arbitrary elements in $A$. \smallskip
 
 The element $v$ plays now its role, since by \Cref{cor:3_lib}(ii) the operators $[\pi(v)^2,\pi(u)]$, $[\pi(v),\pi(u)^2],$ and  $[\pi(v),\pi(u)]$ are linearly independent. Thus, \Cref{independencia1} assures that $g_{c,d} (v) \neq 0$ for some $c,d \in A$. 
 As $g_{c,d} : \mathfrak{J} \rightarrow A$ is continuous there exists $r>0$ such that $g_{c,d}(z) \neq 0$ for every $z \in \mathbf{B}(v,r)$. Here $\mathbf{B}(v,r)$ stands for the open unit ball in $\mathfrak{J}$ of centre $v$ and radius $r$. Since for every $z \in  \mathbf{B}(v,r)$ the operators $[\pi(z)^2,\pi(u)],[\pi(z),\pi(u)^2],[\pi(z),\pi(u)]$ must be linearly independent by \Cref{independencia1}, we can proceed as in \textit{Case I} to deduce that $B(z,z) - \lambda \pi(z)^2 \in \hbox{Span}\{\mathbf{1}, \pi(z)\},$ and hence $f_{a,b}(z) = 0,$ for all $a,b\in A$, $z \in \mathbf{B}(v,r)$ by \Cref{independencia1}. Fix two arbitrary elements $a,b\in A$. \smallskip 
 
Observe that the element $x \in \mathfrak{J}$ fixed at the beginning of this \emph{Case II} satisfies $x \neq v.$ Let us consider a polynomial  mapping from $\mathbb{C}$ to $A$ given by $\zeta  \rightarrow p_{a,b} (\zeta)=f_{a,b}\Big(v + \zeta(x-v)\Big)$,  which is clearly continuous because $f_{a,b}$ is a polynomial. Note that, by what we have just deduced in the previous paragraph, $f_{a,b}\Big(v + \zeta(x-v)\Big) = 0$ for every $\zeta \in \mathbf{B}_{\mathbb{C}}\left(0, \frac{r}{\|x-v\|}\right)$. Therefore the polynomial $p_{a,b} (\zeta)$ must be zero, that is, $f_{a,b}\Big(v + \zeta(x-v)\big) = 0,$ for any $\zeta \in \mathbb{C}$. In particular taking $\zeta = 1$ we get $f_{a,b} (x) = 0$ for every $a,b \in A$, and thus, $B(x,x) - \lambda \pi(x)^2 \in \hbox{Span}\{\mathbf{1}, \pi(x) \}$ by \Cref{independencia1}, as we wanted to prove. \end{enumerate}
 \end{proof}

We establish next our main conclusion on the generic form of those symmetric bilinear maps on special JBW$^*$-algebras without spin part whose trace is associating.  
    
\begin{theorem}\label{prop:JBWspecial} Let $\mathfrak{J}$ be a special JBW$^*$-algebra with no direct summands of type $I_1$ or $I_2$, and let $B : \mathfrak{J} \times \mathfrak{J} \rightarrow \mathfrak{J}$ be a symmetric bilinear mapping. Suppose that 
$[B(x,x), \mathfrak{J}, x] = 0$, for all $x \in \mathfrak{J}$, where the brackets denote the associator with respect to the Jordan product of $\mathfrak{J}$. Then $B$ can be uniquely written in the form: 
    \begin{equation*}
    	B(x,x) = \lambda\circ x^2 + \mu(x)\circ x + \nu(x,x), \ \ \hbox{ for all } x\in \mathfrak{J},
    \end{equation*}
    where $\lambda \in Z (\mathfrak{J})$, $\mu: \mathfrak{J} \rightarrow Z (\mathfrak{J})$ is a linear mapping and $\nu : \mathfrak{J} \times \mathfrak{J} \rightarrow Z (\mathfrak{J})$ is a symmetric bilinear map.
 \end{theorem}

\begin{proof} The proof is divided into 3 steps.\smallskip

\noindent $(1)$ \textit{First Step}:  Let us take an arbitrary Jordan factor representation $\pi : \mathfrak{J} \to \mathfrak{J}_{\pi}$ (here, $\mathfrak{J}_{\pi}$ is a factor JBW$^*$-algebra). We shall study the composition of $B$ with $\pi$, i.e., $
			\pi \circ B : \mathfrak{J} \times \mathfrak{J} \rightarrow \mathfrak{J}_{\pi}$.\smallskip
   
We observe that $\mathfrak{J}_{\pi}$ must be a special and unital JBW$^*$-algebras since $\mathfrak{J}$ so is (cf. \cite[7.2.1., Theorems 7.2.3 and 7.2.7]{HOS} and \cite{Wright1977}).\smallskip

We claim that $\mathfrak{J}_{\pi}$ cannot be a JBW$^*$-algebra factor of type $I_1$ or $I_2$, equivalently, it is a non-quadratic factor (see \cite[Theorem 3.5.5]{CabGarPalVol1} and \Cref{r quadratic}). Namely,  if $\mathfrak{J}_{\pi}$ were a factor of type $I_1$ or of type $I_2$, then it would be quadratic, or equivalently, the square of every element in $\mathfrak{J}_{\pi}$ would be a linear combination of the element itself and the unit. However, by \Cref{cor:3_lib}, there exist $u\in \mathfrak{J}$ and elementary operators $\mathcal{E}_0,\mathcal{E}_1,\mathcal{E}_2 $ in $\mathfrak{J}$ such that $\mathcal{E}_i(u^j) = \delta_{ij}\unit,$ for all $ i,j \in \{0,1,2\}$, and by \Cref{lemma:elem_homo}(ii) there are elementary operators $\widehat{\mathcal{E}_i} \in \element(\mathfrak{J}_{\pi})$ such that $\pi \mathcal{E}_ i = \widehat{\mathcal{E}}_i\pi$ for $i \in \{0,1,2\}$. This contradicts that $\pi (u^2) = \pi (u)^2$ is a linear combination of $\pi(u)$ and $\unit$.  \smallskip
     
Since $\pi$ is a Jordan homomorphism, it follows from the hypotheses on $B(\cdot,\cdot)$ that the identity  
     	\begin{equation*}
     		[\pi(B(x,x)), \pi(\mathfrak{J}), \pi(x)] = 0,
     	\end{equation*}	holds for every $x \in \mathfrak{J}$. In particular, by the weak$^*$-density of $\pi(\mathfrak{J})$ in $\mathfrak{J}_{\pi},$ and the separate weak$^*$-continuity of the Jordan product of $\mathfrak{J}_{\pi},$ we can further conclude that
        \begin{equation}
        	[\pi(B(x,x)), J_{\pi}, \pi(x)] = 0,
        	\label{eq:asocido_pi}
        \end{equation}  for every $x \in \mathfrak{J}$. Furthermore, having in mind that $\mathfrak{J}_{\pi}$ is a factor JW$^*$-algebra not of type $I_1$ or $I_2$, its universal von Neumann algebra, $W^* (\mathfrak{J}_{\pi}),$ satisfies that $\mathfrak{J}_{\pi} = H(W^* (\mathfrak{J}_{\pi}),\tau)$, where $\tau$ is a period-2 $^*$-anti-automorphism on $W^* (\mathfrak{J}_{\pi})$ and the latter either is a factor von Neumann algebra, or it can be written in the form $W^* (\mathfrak{J}_{\pi}) = N\oplus^{\infty}\tau(N)$, for some weak$^*$-closed ideal $N$ of $W^* (\mathfrak{J}_{\pi})$ which is a factor von Neumann algebra (cf. \Cref{l HWtau factor implies W factor} and the comments prior to it or \cite[Theorem 7.1.9, Remark 7.2.8 and Proposition 7.3.3]{HOS}). The product of $W^* (\mathfrak{J}_{\pi})$ will be  denoted by mere juxtaposition. Thus, according to this notation, \eqref{eq:asocido_pi} combined with \Cref{P Topping for non-self adjoint} imply that $\pi(B(x,x))$ and $\pi(x)$ commute in the von Neumann algebra $W^* (\mathfrak{J}_{\pi})$ for all $x\in \mathfrak{J}$, that is, 
        \begin{equation*}
        	[\pi(B(x,x)),\pi(x)] = \pi(B(x,x))\pi(x) - \pi(x)\pi(B(x,x))=  0, \ \forall x\in \mathfrak{J},
        \end{equation*} identity where we employed the associative product of $W^* (\mathfrak{J}_{\pi})$. \smallskip 

\noindent $(2)$ \textit{Second step}: In order to simplify the proof in this step, fix and arbitrary $x\in \mathfrak{J},$ and define $B_0: = \pi(B(x,x)),$ $ B_1:=\pi(B(x,x^2))$ and $B_2:=\pi(B(x^2,x^2))$. Observe that we are in a position to apply \Cref{r identity Bresar}, and more concretely \eqref{equation crazy LeeLee} to deduce that $$ 
\begin{aligned} \pi (\unit) \pi (y) &\left(2 B_1 \pi(x) - B_2  - B_0 \pi(x^2)  \right) + \pi (x) \pi (y) \left(2 B_1 - 2 B_0 \pi(x)  \right) \\ &+ \pi (x^2) (- \pi (y) B_0 )  +  \left(B_2 + \pi (x^2) B_0 - 2 \pi (x) B_1 \right) \pi (y) \pi(\unit) \\ 
&+\left(2 \pi (x) B_0  -2 B_1 \right) \pi (y) \pi (x) + B_0 \pi (y) \pi (x^2) = 0,
\end{aligned}$$ for all $y \in \mathfrak{J}.$ Furthermore, by weak$^*$-density of $\pi (\mathfrak{J})$ in $\mathfrak{J}_{\pi}$, and the separate weak$^*$-continuity of the product of $W^* (\mathfrak{J}_{\pi})$, the previous identity also holds when $\pi(y)$ is replaced by any element $z\in \mathfrak{J}_{\pi} = H(W^* (\mathfrak{J}_{\pi}),\tau),$ that is, 
\begin{equation}\label{eq 26 Jun} 
\begin{aligned} \pi (\unit) z &\left(2 B_1 \pi(x) - B_2  - B_0 \pi(x^2)  \right) + \pi (x) z \left(2 B_1 - 2 B_0 \pi(x)  \right) \\ &- \pi (x^2) z B_0  +  \left(B_2 + \pi (x^2) B_0 - 2 \pi (x) B_1 \right) z \pi(\unit) \\ 
&+\left(2 \pi (x) B_0  -2 B_1 \right) z \pi (x) + B_0 z \pi (x^2) = 0,
\end{aligned}    
\end{equation} for all $z\in \mathfrak{J}_{\pi}= H(W^* (\mathfrak{J}_{\pi}),\tau)$.\smallskip

Three possible cases can be considered here:
\begin{enumerate}[$(1)$]
       \item \underline{\textit{Case 1}}: $\{\mathbf{1}= \pi(\mathbf{1}), \pi(x), \pi(x^2)\} $ are $\mathbb{C}$-linearly independent. Having in mind that $\mathfrak{J}_{\pi} = H(W^* (\mathfrak{J}_{\pi}),\tau)$ is a factor, we can combine \Cref{l HWtau factor implies W factor} and \Cref{independencia2} with the identity in \eqref{eq 26 Jun}, to conclude that 	
     	\begin{equation*}
     		\pi(B(x,x)) = B_0 = \lambda_{0,x} \mathbf{1}+\lambda_{1,x}\pi(x) + \lambda_{2,x}\pi(x^2),
     	\end{equation*}
     	where $\lambda_{i,x} \in \mathbb{C}$ for every $i \in \{0,1,2\}$.\smallskip    
     	
     	\item \underline{\textit{Case 2}}: $\{\mathbf{1}, \pi(x), \pi(x^2)\} $ are $\mathbb{C}$-linearly dependent and $\pi(x) \notin \mathbb{C}\mathbf{1}$. In this case there exist $\alpha, \beta \in \mathbb{C}$ such that 
     	$$ \pi(x^2) = \alpha \pi(x) + \beta \pi(1).$$
      By making some computations, via \eqref{eq new remark 1206} in \Cref{r identity Bresar}, we obtain
       $$\begin{aligned}
        \pi(x)[B_0, \pi(y)] + [B_0,\pi(y)]\pi(x) = & -2\pi(x)[\pi(B(x,y)), \pi(x)]\\ &-2[\pi(B(x,y)), \pi(x)]\pi(x)\\ =& -2[\pi(B(x,y)), \pi(x^2)] \\= & -2[\pi(B(x,y)), \alpha\pi(x) + \beta \pi(\mathbf{1})]\\
        =& -2\alpha[\pi(B(x,y)), \pi(x)] = \alpha[B_0,\pi(y)],
        \end{aligned}
        $$ for all $y\in \mathfrak{J},$
       which leads us to the next identity
       \begin{equation} \label{26 06 (2)}
       \begin{aligned}
           &(\pi(x) - \alpha\pi(\mathbf{1}))B_0 \pi(y)\pi(\mathbf{1}) + B_0\pi(y)\pi(x) + \\ 
       	+& \pi(\mathbf{1})\pi(y)B_0(\alpha \pi(\mathbf{1}) - \pi(x)) - \pi(x)\pi(y)B_0 = 0,
       \end{aligned}
       \end{equation} for all $y\in \mathfrak{J}$ (and thus the conclusion holds when $\pi(y)$ is replaced with any $z\in \mathfrak{J}_{\pi}= H(W^* (\mathfrak{J}_{\pi}),\tau)$).
       Since $\{\pi(\unit), \pi(x)\}$ are $\mathbb{C}$-linearly independent, by a new application of \Cref{independencia2} (see also \Cref{l HWtau factor implies W factor}), applied to \eqref{26 06 (2)}, we conclude that there exist $\lambda_{0,x}, \lambda_{1,x} \in \mathbb{C}$ such that 
       \begin{equation*}
       	B_0= \pi(B(x,x)) = \lambda_{0,x}\pi(\mathbf{1}) + \lambda_{1,x} \pi(x).
       \end{equation*}
       \item \underline{\textit{Case 3}}: $\{\unit, \pi(x), \pi(x^2)\} $ are $\mathbb{C}$-linearly dependent and $\pi(x) \in \mathbb{C}\unit$. Then it is clear from \Cref{r identity Bresar}, and in particular by \eqref{eq new remark 1206}, that  $$[B_0, \pi(y)] = -2[\pi(B(x,y)), \pi(x)] = 0,  \hbox{ for every } y \in \mathfrak{J}.$$ Moreover, due to the weak$^*$-density of $\pi(\mathfrak{J})$ in $\mathfrak{J}_{\pi}$, and the separate weak$^*$-continuity of the product of $W^*(\mathfrak{J_{\pi}})$, the previous identity also holds for any $z \in \mathfrak{J}_{\pi}$, i.e., $[B_0, z] = 0$ for every $ z \in \mathfrak{J}_{\pi}$. This shows that $B_0 \in Z(\mathfrak{J}_{\pi}) = \mathbb{C}\unit$, and thus, that there exists $\lambda_{0,x} \in \mathbb{C}$ such that
       \begin{equation*}
       	B_0 = \pi(B(x,x)) = \lambda_{0,x}\pi(\mathbf{1}).
       \end{equation*}
     \end{enumerate}
     
In summary, all the previous cases lead to the conclusion that for each $x \in \mathfrak{J}$ there exist complex numbers $\lambda_{0,x}, \lambda_{1,x}, \lambda_{2,x} $ such that
      \begin{equation}\label{eq special B}
      	\pi(B(x,x)) =  \lambda_{0,x} \pi(\mathbf{1}) +\lambda_{1,x}\pi(x) + \lambda_{2,x}\pi(x^2).
      \end{equation} 
  	
\noindent $(3)$ \textit{Step 3}: We are now in a position to apply \Cref{lemma:lambda_indep} to deduce the existence of $\alpha \in \mathbb{C} $ such that $\pi(B(x,x)) - \alpha\pi(x)^2 \in \hbox{Span}\{\mathbf{1}, \pi(x)\}$ for every $x \in \mathfrak{J}$. We observe that the symmetric bilinear mapping $\mathfrak{J}\times \mathfrak{J}\to \mathfrak{J}_{\pi},$ $(x,y)\mapsto \pi(B(x,y))-\alpha\pi(x\circ y)$ satisfies the hypothesis in \Cref{lemma:lin_exist}. Therefore, by the just quoted lemma, there exist a linear map $\beta : \mathfrak{J} \rightarrow \mathbb{C}$ and a symmetric bilinear map $\gamma : \mathfrak{J} \times \mathfrak{J} \rightarrow \mathbb{C}$ such that 
  	\begin{equation*}
  		\pi(B(x,x)) = \alpha \pi(x^2) + \beta(x) \pi(x) + \gamma(x,x)\pi(\mathbf{1}),
  	\end{equation*}
     for every $x \in \mathfrak{J}$. Now observe that $\alpha, \beta(x), \gamma(x,x) \in Z (\mathfrak{J}_{\pi}) = \mathbb{C}\mathbf{1},$ for all $x \in \mathfrak{J}$, which allows us to write the above identity in terms of the natural Jordan product defined in $W^*(\mathfrak{J}_{\pi})$ as follows
     \begin{equation}\label{eq linear comb pix}
     \pi(B(x,x)) = \alpha \circ \pi(x^2) + \beta(x) \circ\pi(x) + \gamma(x,x)\circ \pi(\mathbf{1}).
     \end{equation}
     As commented in \emph{Step 1}, there exist an element $\pi(u) \in \mathfrak{J}_{\pi}$ and elementary operators $\widehat{\mathcal{E}_i}$ for all $i \in \{0,1,2\}$ such that $\widehat{\mathcal{E}_i} \pi = \pi \mathcal{E}_i$, and  $\widehat{\mathcal{E}_i}(\pi(u)^j) = \delta_{ij}\unit,$ for all $ i,j \in \{0,1,2\}$. By \Cref{prop: uniq_sol}, applied to the identity in \eqref{eq linear comb pix}, we deduce that $\alpha, \beta(x),$ and $\gamma(x,x)$ are uniquely written as follows:
    \begin{equation}\label{eq expressions for alpha, beta and gamma} \left\{\begin{aligned}
    \alpha &= \widehat{\mathcal{E}_2}(\pi(B(u,u))) = \pi(\mathcal{E}_2({B(u,u)})),\\
    \beta(x) &= 2\widehat{\mathcal{E}_1}(\pi(B(u,x))) - 2\widehat{\mathcal{E}_2}(\pi(B(u,u)))\circ \widehat{\mathcal{E}_1}(\pi(t\circ x))\\ & - \widehat{\mathcal{E}}_1(\pi(B(u,u)))\circ \widehat{\mathcal{E}}_1(\pi(x))\\ = &\pi\Big (2\mathcal{E}_1(B(u,x)) -2\mathcal{E}_2(B(u,u))\circ \mathcal{E}_1(t\circ x) -\mathcal{E}_1(B(u,u))\circ\mathcal{E}_1(x)\Big),\\
    \gamma(x,x&) = \widehat{\mathcal{E}}_0(\pi(B(x,x)))- \alpha\circ\widehat{\mathcal{E}_0}(\pi(x^2)) - \beta(x)\circ\widehat{\mathcal{E}_0}(\pi(x))\\ &= 
    \pi \Big( \mathcal{E}_0(B(x,x)) - \mathcal{E}_2({B(u,u)})\circ \mathcal{E}_0(x^2) -\Big(2\mathcal{E}_1(B(u,x)) \\ &-2\mathcal{E}_2(B(u,u))\circ \mathcal{E}_1(t\circ x)-\mathcal{E}_1(B(u,u)) \circ\mathcal{E}_1(x)\Big)\circ \mathcal{E}_0(x) \Big),
    \end{aligned} \right.
    \end{equation}
    for all $x\in \mathfrak{J}$. Observe that the Jordan factor representation $\pi$ was arbitrarily chosen, so by setting $\lambda = \mathcal{E}_2(B(u,u)),$ and defining $\mu : \mathfrak{J} \rightarrow \mathfrak{J}$ and $\nu : \mathfrak{J} \times \mathfrak{J}  \rightarrow \mathfrak{J}$ by \begin{align*}
      \mu(x) &= 2\mathcal{E}_1(B(u,x)) -2\mathcal{E}_2(B(u,u))\circ \mathcal{E}_1(t\circ x) -\mathcal{E}_1(B(u,u))\circ\mathcal{E}_1(x),\\ \hbox{ and }&
      \nu(x,x) =  \mathcal{E}_0(B(x,x))- \lambda \circ\mathcal{E}_0(x^2) - \mu(x)\circ\mathcal{E}_0(x),
      \end{align*} and combining \eqref{eq linear comb pix} and \eqref{eq expressions for alpha, beta and gamma}, we arrive to 
     \begin{equation}\label{eq B final}
	\pi(B(x,x)) = \pi(\lambda \circ x^2 + \mu(x)\circ x + 
        \nu(x,x)\circ\mathbf{1}),
     \end{equation}
    for every $x \in \mathfrak{J}$, and every Jordan factor representation $\pi$ of $\mathfrak{J}$. Since every JB$^*$-algebra admits a separating family of Jordan factor representations of this form (cf. \cite[Corollary 5.7]{AlfShulStor79GelfandNeumark}), it follows from \eqref{eq B final} that \begin{equation}\label{eq B final without pi}
	B(x,x) = \lambda \circ x^2 + \mu(x)\circ x + 
        \nu(x,x),
     \end{equation}
    for every $x \in \mathfrak{J}$. Note that $\pi(\lambda), \pi(\mu(x))$ and $\pi(\nu(x,x))$ belong to $Z(\mathfrak{J}_{\pi}) = \mathbb{C} \unit,$ for every Jordan factor representation $\pi: \mathfrak{J} \rightarrow \mathfrak{J}_{\pi}$, so by \Cref{Rem: central representation} we finally deduce that $\lambda, \mu(x),$ and $\nu(x,x)$ belong to $Z(\mathfrak{J})$ for each $x \in \mathfrak{J}$ as we desired.	
\end{proof}

\section{Associating traces on purely exceptional \texorpdfstring{JBW$^*$}{JBW*}-algebras}\label{sec: purely exceptional associating trace}

In this section we study symmetric bilinear maps with associating trace on exceptional JBW$^*$-algebras. Our goal is to show that these maps also have the standard form established when the domain is a JW$^*$-algebra without associative and spin part established in  \Cref{prop:JBWspecial}. As we commented in the introduction, the case of exceptional factors was treated by M. Bre\v{s}ar, D. Eremita, and the third author of this paper in \cite[\S 4]{BreEreVil}.  \smallskip

We begin by introducing some notation. Given a set $\Omega$ and a Jordan algebra $\mathfrak{J}$, we shall write $\mathcal{F} (\Omega, \mathfrak{J})$ for the space of all functions from $\Omega$ to $\mathfrak{J}$, and we equip it with the pointwise Jordan product. We can actually consider the Jordan product of an element in the space $C (\Omega, \mathfrak{J}),$ of all continuous functions from $\Omega$ to $\mathfrak{J})$, by an element in $\mathcal{F} (\Omega, \mathfrak{J})$. Given $w\in \mathfrak{J}$, we shall denote by $\widehat{w}$ the constant continuous function $\widehat{w} \in C(\Omega, \mathfrak{J}),$ $\widehat{w} (t) = w$ ($t\in \Omega$). Clearly, each elementary operator $\mathcal{E}\in \element(\mathfrak{J})$ extends to an elementary operator $\widetilde{\mathcal{E}}$ on $C(\Omega, \mathfrak{J})$ given by the corresponding extensions of the involved elements in $\mathfrak{J}$ as constant functions. The extended elementary operator satisfies  $\widetilde{\mathcal{E}}(u)(s) = {\mathcal{E}}(u(s)),$ and  $\widetilde{\mathcal{E}}(\alpha\circ u)(s) =\alpha (s)\circ \widetilde{\mathcal{E}}(u) (s),$ for all $u\in  C(\Omega,\mathfrak{J}),$ $\alpha\in \mathcal{F} (\Omega, \mathbb{C} \unit)$, and $s\in \Omega$.  \smallskip

The proof of \Cref{prop: uniq_sol} is essentially algebraic and the arguments remain valid to establish the following variant.

\begin{proposition}\label{prop: uniq_sol FOmegaJ} Let $\Omega$ be a compact Hausdorff space, $\mathfrak{J}$ a factor JBW$^*$-algebra containing an element $w \in \mathfrak{J}$ and admitting elementary operators $\mathcal{E}_0, \mathcal{E}_1, \mathcal{E}_2$ on $\mathfrak{J}$ such that $\mathcal{E}_i(w^j) = \delta_{ij}\mathbf{1},$ for all $ i,j \in \{0,1,2\}$. Let $B : C(\Omega,\mathfrak{J}) \times C(\Omega,\mathfrak{J}) \rightarrow C(\Omega,\mathfrak{J})$ be a symmetric bilinear mapping admitting the following weak standard representation: 
	$$ B(x,x) (t) = \left(\lambda \circ x^2  + \mu(x)\circ x + \nu(x,x)\right)(t),$$
	for all $x \in C(\Omega,\mathfrak{J})$, $t\in \Omega$, where $\lambda \in \mathcal{F}(\Omega,\mathbb{C} \unit)$, $\mu : C(\Omega,\mathfrak{J}) \rightarrow \mathcal{F}(\Omega,\mathbb{C}\unit)$ is a linear mapping, and $\nu: C(\Omega,\mathfrak{J}) \times C(\Omega,\mathfrak{J})  \rightarrow \mathcal{F}(\Omega,\mathbb{C} \unit)$ is a symmetric bilinear map. Then, this representation is unique and satisfies:
	\begin{enumerate}[$\bullet$]
		\item $\lambda = \widetilde{\mathcal{E}_2}(B(\widehat{w},\widehat{w}))\in C(\Omega,\mathbb{C}\unit)$,
		\item $\mu(x) = 2\widetilde{\mathcal{E}_1}(B(\widehat{w},x)) -2\widetilde{\mathcal{E}_2}(B(\widehat{w},\widehat{w}))\circ \widetilde{\mathcal{E}_1}(\widehat{w}\circ x)) -\widetilde{\mathcal{E}_1}(B(\widehat{w},\widehat{w}))\circ\widetilde{\mathcal{E}_1}(x)$ $\in$ $C(\Omega,\mathbb{C} \unit),$ for every\hyphenation{every} $x$ in $C(\Omega,\mathfrak{J})$, and 
		\item $\nu(x,x):= \widetilde{\mathcal{E}_0}(B(x,x))- \lambda \circ\widetilde{\mathcal{E}_0}(x^2) - \mu(x)\circ\widetilde{\mathcal{E}_0}(x)$ $\in$ $ C(\Omega,\mathbb{C}\unit),$ for every $x$ in $C(\Omega,\mathfrak{J})$, 
	\end{enumerate}
	where each $\widetilde{\mathcal{E}_i}$ is the natural elementary operators in $\element(C(\Omega,\mathfrak{J}))$ extending ${\mathcal{E}_i}$ to $C(\Omega,\mathfrak{J})$.
\end{proposition} 

Our main result towards determining all associating traces on exceptional JBW$^*$-algebras is stated in the next proposition. 

\begin{proposition}\label{prop:JBWex}
Let $\mathfrak{J}$ denote the JB$^*$- algebra $C(\Omega, H_3(\mathbb{O}))$, where $\Omega$ is a compact Hausdorff space. Let  $B : \mathfrak{J} \times \mathfrak{J} \rightarrow \mathfrak{J}$ be a symmetric bilinear map satisfying  
\begin{equation}\label{eq:asociator}
[B(x,x),\mathfrak{J} , x] = 0, \hspace{0.3cm}	  \hbox{ for all $x \in \mathfrak{J}$.} 	
\end{equation} Then there exist a unique element  $\lambda $ in $Z(\mathfrak{J}) = C(X,\mathbb{C})$, a unique linear mapping $\mu : \mathfrak{J} \rightarrow Z(\mathfrak{J}),$ and a unique symmetric bilinear map $\nu : \mathfrak{J} \times \mathfrak{J} \rightarrow Z(\mathfrak{J})$ such that 
\begin{equation*}
B(x,x) = \lambda \circ x^2 + \mu(x) \circ x + \nu(x,x), \ \hbox{ for all } x\in \mathfrak{J}.\end{equation*}
\end{proposition}

\begin{proof} Fix an arbitrary $t \in \Omega$. According to the notation of this section, given $a\in H_3(\mathbb{O})$, we shall write $\hat{a}$ for the continuous function in $\mathfrak{J}$ defined by $\hat{a}: \Omega \rightarrow H_3(\mathbb{O}),$ $s \mapsto \hat{a}(s) := a$. For each $u \in \mathfrak{J}$, we write $u_t := \widehat{u(t)}: \Omega \rightarrow H_3(\mathbb{O}),$ $s \mapsto u_t(s) := u(t)$.\smallskip
	
Note that for every $u \in \mathfrak{J}$ the bilinear mapping $B$ can be decomposed as follows: 
 \begin{equation}\label{eq first identity B(u,u) as three summands} B(u,u) = B(u-u_t, u-u_t) + 2B(u-u_t, u_t) - B(u_t,u_t).
 \end{equation} Let us see how each summand on the right hand side of the previous identity can be represented.\smallskip
 
\noindent $(1)$ We begin with the term $B(u_t,u_t)$. Note that for each $s \in \Omega$ we can define a symmetric bilinear mapping $B_{s}: H_3(\mathbb{O})\times H_3(\mathbb{O}) \to H_3(\mathbb{O}),$ $B_s (a,b) := B(\hat{a},\hat{b}) (s).$ We claim that $B_s$ is associating. Namely, if we take $a,b\in H_3(\mathbb{O})$ we can see from the hypotheses on $B$ that $$[B_s(a,a),b,a] = [B(\hat{a},\hat{a}),\hat{b}, \hat{a}] (s) = 0.$$ Thus,  by applying  \cite[Lemma 4.1]{BreEreVil}, $B_s(\cdot,\cdot)$ admits the following representation:
\begin{equation}\label{eq expression for Bs}
	B_s (a,a) = B(\hat{a},\hat{a})(s) = \lambda_s\circ  a^2 + \mu_s (a)\circ a + \nu_s(a, a),
\end{equation}
where $\lambda_s \in \mathbb{C} \ \unit $, $\mu_s:H_3(\mathbb{O})  \to \mathbb{C}\ \unit$ is a linear map, and $\nu_s: H_3(\mathbb{O})\times H_3(\mathbb{O}) \to \mathbb{C}\!\ \unit$ is a symmetric bilinear mapping, where $\unit$ stands for the unit of $H_3(\mathbb{O})$, and this procedure can be done for every $s\in \Omega$.\smallskip
      
Note that, by \Cref{cor:3_lib} there exists an element $w \in H_3(\mathbb{O})_{sa}$ and elementary operators $\mathcal{E}_0,\mathcal{E}_1,\mathcal{E}_2$ such that $\mathcal{E}_i(w^j) = \delta_{i j}\unit$ ($i,j \in \{0,1,2\}$). We can apply \Cref{prop: uniq_sol}, with the Jordan homomorphism $\pi: H_3(\mathbb{O})\to  H_3(\mathbb{O})$ as the identity operator, to deduce that
\begin{equation}\label{eq B_s}\left\{
	\begin{aligned}
		\lambda_s & = \mathcal{E}_2(B_s(w, w)),\\ \mu_s(a) &= 2\mathcal{E}_1(B_s(w, a)) -2\lambda_s \circ \mathcal{E}_1(w\circ a) -\mathcal{E}_1(B_s(w,w))\circ \mathcal{E}_1(a),\\ 
		\nu_s (a, a) &= \mathcal{E}_0(B_s (a, a)) - \lambda_s\circ  \mathcal{E}_0(a^2) - \mu_s(a)\circ  \mathcal{E}_0(a), 
		\end{aligned} \right.
\end{equation} for all $a\in H_3(\mathbb{O})$ and all $s\in \Omega$.\smallskip

Let us consider the natural extension of $w \in H_3(\mathbb{O})$ as the constant continuous function $\widehat{w} \in \mathfrak{J} =C(\Omega, H_3(\mathbb{O})),$\label{existence of tildew 030924} and the corresponding extension of the elementary operator $\mathcal{E}_i$ ($ i \in \{0,1,2\}$) as an elementary operator $\widetilde{\mathcal{E}_i}$ on $\mathfrak{J} =C(\Omega, H_3(\mathbb{O}))$ ($ i \in \{0,1,2\}$) given by the corresponding extensions of the involved elements as constant functions. By construction, $\widetilde{\mathcal{E}_i}(\widehat{w}^j)(s) =\mathcal{E}_i(\widehat{w}^j(s))= \mathcal{E}_i(w^j)= \delta_{i j}\mathbf{1},$ for every $s \in \Omega$. \smallskip

We also set \begin{equation}\label{eq def of lambda mu and nu}\left\{
	\begin{aligned}
		\lambda & := \widetilde{\mathcal{E}_2}(B(\widehat{w}, \widehat{w}))\in \mathfrak{J} = C(\Omega, H_3(\mathbb{O})),\\ 
		\mu (u) &:= 2\widetilde{\mathcal{E}_1}(B(\widehat{w}, u)) -2\lambda \circ \widetilde{\mathcal{E}_1}(\widehat{w}\circ u) -\widetilde{\mathcal{E}_1}(B(\widehat{w},\widehat{w}))\circ \widetilde{\mathcal{E}_1}(u)\in \mathfrak{J},\\ 
		\nu (u, u) &:= \widetilde{\mathcal{E}_0}(B (u, u)) - \lambda\circ  \widetilde{\mathcal{E}_0}(u^2) - \mu (u)\circ  \widetilde{\mathcal{E}_0}(u)\in \mathfrak{J}, 
	\end{aligned} \right.
\end{equation} for all $u\in\mathfrak{J} = C(\Omega, H_3(\mathbb{O}))$. Clearly, $\mu : \mathfrak{J}\to \mathfrak{J}$ is a linear mapping, while $\nu$ induces a  symmetric bilinear mapping on $\mathfrak{J}$. \smallskip

Now, by applying $\widetilde{\mathcal{E}_2} $ to $B(\widehat{w}, \widehat{w})$ and having in mind \eqref{eq B_s} we obtain 
     \begin{equation}\label{eq lambdas and lambda(s)}
	  \lambda (s) = \widetilde{\mathcal{E}_2}(B(\widehat{w}, \widehat{w}))(s) = \mathcal{E}_2(B(\widehat{w},\widehat{w})(s))= \mathcal{E}_2(B_s(w,w)) = \lambda_s\in \mathbb{C}\unit,
     \end{equation}
    for every $s\in \Omega$. Observe that $\lambda (s)\in \mathbb{C} \unit$ for  every $s\in \Omega$, and thus $\lambda\in Z(\mathfrak{J})$.  
    
Fix arbitrary elements $s,t\in \Omega$ and $u\in \mathfrak{J}$. By \eqref{eq B_s}, \eqref{eq def of lambda mu and nu}, and \eqref{eq lambdas and lambda(s)} we also get 
\begin{equation}\label{eq mus ut and mu ut (s)}\left\{ \begin{aligned}
		\mathbb{C} \unit \ni \mu_s(u(t)) &= {\mathcal{E}_1}\Big( 2 B_s({w}, u(t))\Big)   - 2\lambda(s)\circ {\mathcal{E}_1} \Big( w\circ u(t)\Big) \\
		&- {\mathcal{E}_1}\Big( B_s({w}, {w})\Big) \circ {\mathcal{E}_1}\Big( u(t)\Big)  \\
		&= 2 {\mathcal{E}_1}(B(\widehat{w}, u_t)  (s) ) - 2\lambda(s)\circ {\mathcal{E}_1}(\widehat{w}(s)\circ u_t(s))\\
		& - {\mathcal{E}_1}(B(\widehat{w}, \widehat{w})(s)) \circ {\mathcal{E}_1}(u_t(s)) \\
		&= 2 \widetilde{\mathcal{E}_1}(B(\widehat{w}, u_t)) (s) -2\lambda(s)\circ \widetilde{\mathcal{E}_1}(\widehat{w}\circ u_t) (s)  \\
		&- \widetilde{\mathcal{E}_1}(B(\widehat{w}, \widehat{w}))(s) \circ \widetilde{\mathcal{E}_1}(u_t)(s) \\
		=\Big(2  \widetilde{\mathcal{E}_1}&(B(\widehat{w}, u_t))-2\lambda\circ \widetilde{\mathcal{E}_1}(\widehat{w}\circ u_t) - \widetilde{\mathcal{E}_1}(B(\widehat{w}, \widehat{w}))\circ \widetilde{\mathcal{E}_1}(u_t)\Big)(s)\\
		& = \mu (u_t) (s),
	\end{aligned}\right.
\end{equation} and consequently, $\mu (u_t)\in Z(\mathfrak{J})$ for all $u\in \mathfrak{J}$, $t\in \Omega$. \smallskip

Fix $u\in \mathfrak{J}$ and $s,t\in \Omega$. By definition of $\nu (\cdot,\cdot)$ (cf. \eqref{eq def of lambda mu and nu}), and having in mind \eqref{eq B_s}, \eqref{eq def of lambda mu and nu}, \eqref{eq lambdas and lambda(s)}, and \eqref{eq mus ut and mu ut (s)}, we obtain
\begin{equation}\label{eq nus} \left\{ 
	\begin{aligned}\nu (u_t,u_t) (s) &=  \Big(\widetilde{\mathcal{E}_0}(B(u_t,u_t)) - \lambda\circ \widetilde{\mathcal{E}_0}(u_t^2) - \mu(u_t)\circ \widetilde{\mathcal{E}_0}(u_t)\Big)(s) \\ 
		&=  {\mathcal{E}_0}(B_s(u_t,u_t)) - \lambda (s) \circ {\mathcal{E}_0}(u(t)^2) - \mu (u_t) (s) \circ {\mathcal{E}_0}(u(t)) \\ 
		&=  {\mathcal{E}_0}(B_s(u_t,u_t)) - \lambda _{s} \circ {\mathcal{E}_0}(u(t)^2) - \mu_{s} (u(t)) \circ {\mathcal{E}_0}(u(t)) \\  
		&= \nu_s(u(t),u(t))\in \mathbb{C} \unit, \hbox{ for every $s \in X$.}
	\end{aligned} \right.
\end{equation} Consequently, $\nu (u_t,u_t)\in Z(\mathfrak{J})$ for all $u\in \mathfrak{J}$, $t\in \Omega.$\smallskip

The previous conclusions show, via \eqref{eq expression for Bs}, \eqref{eq lambdas and lambda(s)}, \eqref{eq mus ut and mu ut (s)}, and \eqref{eq nus}, that the term $B(u_t,u_t)$ can be written as follows
    \begin{equation}\label{eq B(ut,ut)}\left\{ \begin{aligned}
    		B(u_t,u_t)(s) &= B(\widehat{u(t)},\widehat{u(t)}) (s) = B_s (u(t),u(t)) \\
    		&= \lambda_s \circ u(t)^2 + \mu_s (u(t)) \circ u(t) + \nu_s (u(t),u(t)) \circ \unit \\
    		&= \Big(\lambda\circ  u_t^2 + \mu(u_t)\circ u_t + \nu(u_t,u_t)\circ \unit\Big)(s).
    	\end{aligned} \right.
    \end{equation}

\noindent $(2)$ We deal next with the term $B(u-u_t, u-u_t)$. It is not hard to check that for every $a, b \in H_3(\mathbb{O})$ and $s \in \Omega$, by \eqref{eq new remark 1206 Jordan} in \Cref{r identity Bresar}, we have 
\begin{align*}
	  [B(u-u_t,u-u_t)(s), a, b] &= [B(u-u_t,u-u_t), \widehat{a}, \widehat{b}](s)\\ 
	                            &=-2[B(u-u_t, \widehat{b}), \widehat{a}, u-u_t](s).
\end{align*} When $s = t$, we obtain that $B(u-u_t,u-u_t)(t) \in \mathbb{C} \unit = Z(H_3(\mathbb{O})).$ We can thus define a symmetric bilinear mapping $\gamma$ from $\mathfrak{J} \times \mathfrak{J}\to \mathcal{F} (\Omega,\mathbb{C}\unit)$ given by \begin{equation}\label{eq B(u-ut,u-ut)} \gamma(u,u)(t) =B(u-u_t, u-u_t) (t) \in \mathbb{C}\unit.
\end{equation}

\noindent $(3)$ It is now the turn to handle the term $B(u-u_t, u_t)$. For this purpose, for each $s\in \Omega$ and $u\in \mathfrak{J}$, we define a linear mapping $T_{s,u}: H_3(\mathbb{O}) \rightarrow H_3(\mathbb{O})$ defined by 
     \begin{equation*}
   	  T_{s,u} (a) :=B(u-u_t, \widehat{a})(s).
     \end{equation*} If we pick an arbitrary element $b \in H_3(\mathbb{O})$, it is easy to see, via a new application of \eqref{eq new remark 1206 Jordan}, that 
     \begin{align*}
     	[T_{s,u}(a), b, a] &= [B(u-u_t, \widehat{a})(s), \widehat{b}(s), \widehat{a}(s)]\\
     	             &=  [B(u-u_t, \widehat{a}), \widehat{b}, \widehat{a}](s)\\
     	             &= -\frac{1}{2} [B(\widehat{a}, \widehat{a}), \widehat{b}, u-u_t](s). 
     \end{align*} For $s= t$, we obtain that $[T_{t,u}(a), b, a] = 0$ for every $a\in H_3(\mathbb{O})$. Thus, by applying \cite[Theorem 2.1]{BreCabVil} we deduce the existence of $\alpha_{t,u} \in \mathbb{C}$ and a linear map $\beta_{t,u}: H_3(\mathbb{O}) \rightarrow \mathbb{C}$ such that
     \begin{equation*}
    	T_{t,u} (a) = \alpha_{t,u} a + \beta_{t,u} (a) \unit, \hbox{ for all } a\in H_3(\mathbb{O}).
     \end{equation*}
When taking ${a} = w,$ we obtain the following identity
     \begin{equation*}
     	B(u-u_t, \widehat{w})(t) = T_{t,u} (w) = \alpha_{t,u} w + \beta_{t,u} (w)\unit.
     \end{equation*} 
     
For each $u\in\mathfrak{J}$,  set $\alpha (u), \beta(u,u) : \Omega \to \mathbb{C}\unit $ given by \begin{align*}
	\alpha(u) (t)  &:= \widetilde{\mathcal{E}_1} (B(u-u_t, \widehat{w})) (t) = {\mathcal{E}_1} (B(u-u_t, \widehat{w}) (t)) = \alpha_{t,u} \unit, \hbox{ and } \\
	\beta(u,u)(t) & := \widetilde{\mathcal{E}_0} (B(u-u_t, u_t))(t) - \alpha(u)(t) \circ \mathcal{E}_0(u(t))= \beta_{t,u}(u(t)), \hbox{ respectively}.
\end{align*} It is easy to check that the assignment $\alpha: C(\Omega, H_3(\mathbb{O})) \to \mathcal{F} (\Omega, \mathbb{C} \unit)$, $u\mapsto \alpha(u)$ (respectively, $\beta: C(\Omega, H_3(\mathbb{O}))\times C(\Omega, H_3(\mathbb{O})) \to \mathcal{F} (\Omega, \mathbb{C}\unit),$ $(u,u)\mapsto \beta(u,u)$) defines a linear (respectively, symmetric bilinear) mapping. According to this notation we can write 
\begin{equation}\label{eq B(u-ut,ut) at t}      	B(u-u_t, u_t)(t) = T_{t,u} (u(t)) = \alpha(u) (t)\circ u(t) + \beta(u,u) (t).
\end{equation}

Back to \eqref{eq first identity B(u,u) as three summands}, \eqref{eq B(u-ut,u-ut)}, \eqref{eq B(u-ut,ut) at t} and \eqref{eq B(ut,ut)} we conclude that 
$$\begin{aligned}
B(u,u)(t) &= \lambda(t)\circ u(t)^2 + \Big(\mu(u_t) + \alpha(u)\Big)(t)\circ u(t) \\
&+ \Big(\nu(u_t,u_t) + \beta(u,u)+ \gamma(u,u) \Big)(t),
\end{aligned}$$ for all $u\in \mathfrak{J}$, $t\in \Omega$. \Cref{prop: uniq_sol FOmegaJ} asserts that $\mu (u) = \mu(u_t) + \alpha(u)$, and $\nu (u,u) = \nu(u_t,u_t) + \beta(u,u)+ \gamma(u,u)$, for all $u\in \mathfrak{J}$, which completes the proof of the existence of the desired decomposition. The uniqueness claimed in the conclusion follows from \Cref{prop: uniq_sol FOmegaJ}.
\end{proof}

Note that, by \cite[Theorem 3.9]{Shultz79}, every JBW-algebra $\mathfrak{J}$ admits a unique decomposition in the form $\mathfrak{J}= \mathfrak{J}_{sp}\oplus^{\infty}\mathfrak{J}_{ex}$, where $\mathfrak{J}_{sp}$ is special (i.e. a JW-algebra) and $\mathfrak{J}_{ex}$ is a purely exceptional JBW-algebra. Moreover, $\mathfrak{J}_{ex}$ is isomorphic to $C(\Omega, H_3(\mathbb{O}))_{sa}$, where $\Omega$ is a \emph{hyperStonean} compact Hausdorff space.  Therefore, each purely exceptional JBW$^*$-algebra is JBW$^*$-isomorphic to $C(\Omega, H_3(\mathbb{O}))$ for some compact Hausdorff space $\Omega$ (cf. \cite[Theorem 2.8]{Wright1977} and \cite[Theorem 3.2]{Edwards80}). The next result is a straightforward consequence of these comments and \Cref{prop:JBWex} above.

\begin{corollary}\label{cor excep bil}
    Let $\mathfrak{J}$ be a purely exceptional JBW$^*$-algebra. Let  $B : \mathfrak{J} \times \mathfrak{J} \rightarrow \mathfrak{J}$ be a symmetric bilinear map satisfying  $[B(x,x),\mathfrak{J} , x] = 0,$	for all $x \in \mathfrak{J}$. Then there exist a unique element $\lambda$ in $Z(\mathfrak{J})$, a unique linear mapping $\mu : \mathfrak{J} \rightarrow Z(\mathfrak{J})$, and a unique symmetric bilinear map  $\nu : \mathfrak{J} \times \mathfrak{J} \rightarrow Z(\mathfrak{J})$ such that $B$ admits the following representation: 
\begin{equation*}
B(x,x) = \lambda\circ  x^2 + \mu(x)\circ x + \nu(x,x), \ \hbox{ for all } x\in \mathfrak{J}.\end{equation*}
\end{corollary}

\section{Associating traces on type~\texorpdfstring{$I_2$}{I2} \texorpdfstring{JW$^*$}{JW*}-algebras}

The main result in \Cref{sec: bilinear with associating trace on JW*-algebras without spin} does not cover the case of associating traces on a JW$^*$-algebra admitting summands of type $I_1$ and $I_2$. This section is devoted to fill in the gap in the case of JW$^*$-algebras of type $I_2$. We recall that every JW$^*$-algebra of type $I_2$ admits a separating family of Jordan factor representations into spin factors (cf. \cite[Propositions~4.6.4, Proposition~5.3.12, and Theorem~6.1.8]{HOS}).\smallskip

We begin with a variant of \Cref{lemma:lin_exist}. Observe first that every non-associative factor JBW$^*$-algebra must be at least $3$-dimesional (cf. \cite[5.3.3]{HOS} and the fact that the dimension of every spin factor is at least $3$). So, if $\pi: \mathfrak{J}\to \mathfrak{J}_{\pi}$ is a Jordan factor representation, where $\mathfrak{J}$ is a JBW$^*$-algebra without type $I_1$ summands, then dim$(\pi(\mathfrak{J}))\geq 3$ (cf. \Cref{cor:3_lib including spin part}).

\begin{lemma}\label{lemma:lin_exist spin} Let $\mathfrak{J}$ be a JBW$^*$-algebra admitting no direct summands of type $I_1$, and let  $\pi: \mathfrak{J} \rightarrow \mathfrak{J}_{\pi}$ be a Jordan factor representation. Assume that $B : \mathfrak{J} \times \mathfrak{J} \rightarrow \mathfrak{J}_{\pi}$ is a symmetric bilinear mapping satisfying that for each $x \in \mathfrak{J}$ there exist $\alpha_x, \beta_x \in \mathbb{C}$ such that  $$B(x,x) = \alpha_x \pi(x) + \beta_x \mathbf{1}.$$  Then there exist a linear mapping $\mu: \mathfrak{J} \rightarrow \mathbb{C}$ and a symmetric bilinear map $\nu: \mathfrak{J} \times \mathfrak{J} \rightarrow \mathbb{C}$ such that 
    $$ B(x,x) = \mu(x) \pi(x) + \nu(x,x) \mathbf{1},$$
    for every $x \in \mathfrak{J}$.
\end{lemma}

\begin{proof} \Cref{cor:3_lib including spin part} implies the existence of an element $w \in \mathfrak{J}$ and elementary operators $\mathcal{E}_0, \mathcal{E}_1 \in \element(\mathfrak{J})$  such that $\mathcal{E}_i(w^j) = \delta_{ij}\mathbf{1},$ for all $i,j \in {0,1}$. By \Cref{lemma:elem_homo}(ii) we can find elementary operators $\widehat{\mathcal{E}}_{j}\in \element (\mathfrak{J}_{\pi})$ satisfying $\pi\mathcal{E}_i = \widehat{\mathcal{E}_i}\pi,$ for all $i \in \{0,1\}$. Consequently, $\widehat{\mathcal{E}}_{j} (\pi(w)^i) = \delta_{ij} \unit,$ $\forall i,j\in \{0,1\}$. \smallskip

As in the proof of \Cref{lemma:lin_exist}, the identity \begin{equation*}
		B(x+w, x+w) + B(x-w,x-w) = 2B(x,x)+ 2 B(w,w), \ \ x\in \mathfrak{J}, 
	\end{equation*}  implies that
     \begin{equation}\label{eq:lin_indep new}
     	\alpha_{x+w} = \alpha_x + \alpha_w,
     \end{equation}
     for every $x \in \mathfrak{J}$ such that $\pi(x)$ is not a linear combination of $\mathbf{1}$ and $\pi(w)$ --observe that the vectors $\mathbf{1}$ and $\pi(w)$ are linearly independent--. We define the mappings $\mu: \mathfrak{J} \rightarrow \mathfrak{J}_{\pi},$ and $\nu: \mathfrak{J} \times \mathfrak{J} \rightarrow \mathfrak{J}_{\pi}$ by
      \begin{align*}
      	\mu(x):=&2\widehat{\mathcal{E}}_1(B(x,w)) -\widehat{\mathcal{E}_1}(B(w,w))\circ \widehat{\mathcal{E}}_1(\pi(x)),\\
      	\nu(x,x):=& B(x,x) - \mu(x)\circ \pi(x).
      \end{align*} 
      By definition, $\mu$ is a linear map while $\nu$ is a bilinear map. It only remains to prove that both maps take values in $\mathbb{C}\mathbf{1}$. Indeed, by mimicking the proof of \Cref{lemma:lin_exist}
      \begin{enumerate}[$(1)$]
\item If $\pi(x) \notin \hbox{Span}\{\mathbf{1}, \pi(w)\}$,  $\mu(x) = 2 \widehat{\mathcal{E}}_1(B(x,w))- \widehat{\mathcal{E}_1}(B(w,w))\circ \widehat{\mathcal{E}}_1(\pi(x))  = \alpha_x\unit \in \mathbb{C}\mathbf{1}$, and $\nu(x,x)=B(x,x) - \mu(x)\circ \pi(x) =B(x,x) - \alpha_x \circ \pi(x) = \beta_x \unit.$
\item If $\pi(x) \in \hbox{Span} \{\mathbf{1}, \pi(w)\}$, by the comments preceding this lemma, there exists $y\in \mathfrak{J}$ such that $\pi(y)  \notin \hbox{Span} \{\mathbf{1}, \pi(w)\}$, and we can thus conclude that for every $\varepsilon > 0,$ the element $\pi(x +\varepsilon y)$ is not in $\hbox{Span} \{\mathbf{1}, \pi(w)\}$. By applying \eqref{eq:lin_indep new} and the conclusion in $(1)$ we obtain
      	\begin{equation*}
      		B(x+\varepsilon y,x+\varepsilon y ) - \mu(x+\varepsilon y)\pi(x+\varepsilon y) \in \mathbb{C}\mathbf{1}.
      	\end{equation*}
       Note that, by what we have just proved in $(1)$, $\mu(x +\varepsilon y) = \alpha_{x + \varepsilon y} \unit.$ We also know that $\mu$ is linear. Thus, it follows that 
       $$\begin{aligned}
       	B(x,x) + 2B(x,\varepsilon y) + \varepsilon^2B(y,y) &- (\mu(x) + \varepsilon\mu(y))\circ\pi(x+\varepsilon y) \in \mathbb{C}\mathbf{1},\\
        \mu (x + \varepsilon y) = & \mu(x) + \varepsilon \mu (y)\in \mathbb{C}\unit.
       \end{aligned}$$  Taking $\varepsilon \rightarrow 0$ we get 
      	$\nu(x,x)= B(x,x) - \mu(x)\circ \pi(x) \in \mathbb{C}\mathbf{1},$  and $\mu (x)\in \mathbb{C}\unit,$
       for every $x \in \mathfrak{J}$.
      \end{enumerate}\end{proof}

The next result determines the general form of all associating traces on a JBW$^*$-algebra of type $I_2$. Briefly speaking, we shall see that the set of all  associating traces on JBW$^*$-algebras of type $I_2$ constitute a more restrictive set than what we have in other cases.  

\begin{theorem}\label{t bilinear associating trace on spin}
Let $\mathfrak{J}$ be a JBW$^*$-algebra of type $I_2$. Let  $B : \mathfrak{J} \times \mathfrak{J} \rightarrow \mathfrak{J}$ be a symmetric bilinear map satisfying  
$$
[B(x,x),\mathfrak{J} , x] = 0, \hspace{0.3cm}	  \hbox{ for all $x \in \mathfrak{J}$.}$$ Then $B$ admits the following  unique representation: 
\begin{equation*}
B(x,x) =  \mu(x)\circ x + \nu(x,x), \ \hbox{ for all } x\in \mathfrak{J},\end{equation*}
where $\mu : \mathfrak{J} \rightarrow Z(\mathfrak{J})$ is linear, and $\nu : \mathfrak{J} \times \mathfrak{J} \rightarrow Z(\mathfrak{J}) $ is bilinear and symmetric.
\end{theorem}

\begin{proof} Fix an arbitrary Jordan factor representation $\pi : \mathfrak{J}\to \mathfrak{J}_{\pi}$, where $\mathfrak{J}_{\pi}$ must be a spin factor by the hypothesis on $\mathfrak{J}$. Take an arbitrary $x\in \mathfrak{J}$. Since, by assumptions, $B(x,x)$ and $x$ operator commute, the same conclusion holds for $\pi B(x,x)$ and $\pi(x)$. Therefore $\pi B(x,x)$ must be a linear combination of $\pi(x)$ and $\pi (\unit)$ (cf. \Cref{r associator-commutator of an element in a spin}). By applying  \Cref{lemma:lin_exist spin} we prove the existence of a linear mapping $\mu_{\pi}: \mathfrak{J} \rightarrow \mathbb{C}$ and a bilinear mapping $\nu_{\pi}: \mathfrak{J}\times \mathfrak{J}\rightarrow \mathbb{C}$ such that 
\begin{equation}\label{eq 2907 A}
    \pi(B(x,x)) = \mu_{\pi}(x)\circ \pi(x) + \nu_{\pi}(x,x).
\end{equation}

By \Cref{cor:3_lib including spin part} and \Cref{lemma:elem_homo}(ii) we find an element $w \in \mathfrak{J}$ and elementary operators ${\mathcal{E}}_{j}\in \element (\mathfrak{J}),$ $\widehat{\mathcal{E}}_{j}\in \element (\mathfrak{J}_{\pi})$ satisfying $\mathcal{E}_i(w^j)=\delta_{ij}\mathbf{1},$ for all $i,j$ in $\{0,1\}$, and $\pi\mathcal{E}_i = \widehat{\mathcal{E}_i}\pi,$ for all $i \in \{0,1\}$. Due to \Cref{prop: uniq_sol spin}, the previous representation of $\pi B$ is unique and satisfies 
\begin{equation}\label{eq 2907 B}
    \begin{aligned}
        \mu_{\pi}(x) &= 2\widehat{\mathcal{E}_1}(\pi(B(w,x)))  -\widehat{\mathcal{E}_1}(\pi(B(w,w)))\circ\widehat{\mathcal{E}_1}(\pi(x)), \\
        \nu_{\pi}(x,x) & = \widehat{\mathcal{E}_0}(\pi(B(x,x))) - \mu_{\pi}(x)\circ\widehat{\mathcal{E}_0}(\pi(x)),
    \end{aligned}
\end{equation}
for every $x \in \mathfrak{J}$. Observe that the Jordan representation $\pi$ was arbitrarily chosen, so defining
$\mu: \mathfrak{J}\rightarrow \mathfrak{J}$ and $\nu: \mathfrak{J}\times \mathfrak{J} \rightarrow\mathfrak{J}$ by 
\begin{equation} \label{eq 2907 defi}
    \left \{\begin{aligned} 
        \mu(x) &= 2\mathcal{E}_1(B(w,x)) -\mathcal{E}_1(B(w,w))\circ\mathcal{E}_1(x),\\
      \nu(x,x) &=  \mathcal{E}_0(B(x,x))- \mu(x)\circ\mathcal{E}_0(x),
    \end{aligned}\right.
\end{equation}
and combining \eqref{eq 2907 A}, \eqref{eq 2907 B} and \eqref{eq 2907 defi} we arrive to 
\begin{equation}\label{eq 2907 C}
    \pi(B(x,x)) = \pi(\mu(x)\circ x + \nu(x,x))
\end{equation}
for every $x \in \mathfrak{J}$, and every Jordan factor representation $\pi$ of $\mathfrak{J}$. Since every JB$^*$-algebra admits a separating family of Jordan factor representations of this form (cf. \cite[Corollary 5.7]{AlfShulStor79GelfandNeumark} and the comments at the beginning of this section), it follows from \eqref{eq 2907 C} that \begin{equation}\label{eq spin C final without pi}
	B(x,x) = \mu(x)\circ x + \nu(x,x)
     \end{equation}
    for every $x \in \mathfrak{J}$. Note that $\pi(\mu(x))$ and $\pi(\nu(x,x))$ belong to $Z(\mathfrak{J}_{\pi}) = \mathbb{C} \unit$ for every Jordan factor representation $\pi: \mathfrak{J} \rightarrow \mathfrak{J}_{\pi}$, so it follows from \Cref{Rem: central representation} that $\mu(x)$ and $\nu(x,x)$ belong to $Z(\mathfrak{J}),$ for each $x \in \mathfrak{J},$ as we desired.	
\end{proof}

\section{Associating linear maps on \texorpdfstring{JBW$^*$}{JBW*}-algebras}\label{sec: associating linear maps}

It is perhaps of independent interest to determine the structure of those linear maps $T$ on a JBW$^*$-algebra $\mathfrak{J}$ enjoying the property of being associating, i.e. $T(x)$ and $x$ operator commute for all $x$ in $\mathfrak{J}$. By obvious reasons we must avoid JBW$^*$-algebras admitting associative summands --observe that every linear mapping between associative JB$^*$-algebras is associating. In this note, these type of linear maps will also play a special role in the next section.   

\begin{lemma}\label{l linear maps from scalars to maps} Let $\mathfrak{J}, \mathfrak{M}$ be JBW$^*$-algebras where $\mathfrak{J}$ has no direct summands of type $I_1$. Let $\pi: \mathfrak{J} \rightarrow \mathfrak{M}$ be a Jordan homomorphism whose image is at least 3-dimensional and weak$^*$-dense in $\mathfrak{M}$, and let $T : \mathfrak{J} \rightarrow \mathfrak{M}$ be a linear map. Suppose that  for each $x \in \mathfrak{J}$ there exist $\alpha_x, \beta_x \in \mathbb{C}$ such that  $$T(x) = \alpha_x \pi(x) + \beta_x \mathbf{1}.$$  Then there exist a complex number $\lambda \in \mathbb{C}$ and a linear map $\mu: \mathfrak{J} \rightarrow \mathbb{C}$ such that 
$$ T(x) = \lambda \pi(x) + \mu(x)\mathbf{1},$$
 for every $x \in \mathfrak{J}$.
\end{lemma}

\begin{proof} By \Cref{cor:3_lib including spin part} and \Cref{lemma:elem_homo}(ii) we can find an element $u \in \mathfrak{J}$ and elementary operators ${\mathcal{E}}_{j}\in \element (\mathfrak{J}),$ $\widehat{\mathcal{E}}_{j}\in \element (\mathfrak{M})$ satisfying $\mathcal{E}_i(u^j)=\delta_{ij}\mathbf{1},$ for all $i,j$ in $\{0,1\}$, and $\pi\mathcal{E}_i = \widehat{\mathcal{E}_i}\pi,$ for all $i \in \{0,1\}$. In particular, $\widehat{\mathcal{E}}_{j} (\pi(u)^i) = \delta_{ij} \unit,$ for all $i,j\in \{0,1\}$.\smallskip

Let us fix an arbitrary $x \in \mathfrak{J}$. The identities $$T(x +u) -  T(x-u) = 2T(u), \hbox{ and }  T(x +u) +  T(x-u)  = 2T(x) $$
    hold by hypothesis. Thus we have 
    $$\begin{aligned}
        (\lambda_{x+u} - \lambda_{x-u})\pi(x) & + (\lambda_{x+u} + \lambda_{x-u} - 2\lambda_u)\pi(u) + (\mu_{x+u} - \mu_{x-u} - 2\mu_u)\unit = 0,\\
        (\lambda_{x+u} - \lambda_{x-u})\pi(u) & + (\lambda_{x+u} + \lambda_{x-u} - 2\lambda_x)\pi(x) + (\mu_{x+u} + \mu_{x-u} - 2\mu_x)\unit = 0,
    \end{aligned}$$
for every $x \in \mathfrak{J}$. If the vectors in the set $\{\unit, \pi(x), \pi(u)\}$ are linearly independent, we deduce from the previous identities that $\lambda_{x+u} = \lambda_{x-u} = \lambda_u$ and $\mu_{x+u} = \mu_x + \mu_u$. Furthermore, since 
    $$ T(x+u) = \lambda_u(\pi(x) + \pi(u)) + (\mu_x + \mu_u)\unit,$$ we derive that $\lambda_x = \lambda_u$.
Let us define $\lambda: = \lambda_u = \widehat{\mathcal{E}}_1(T(u))$. We shall prove that $T(x) - \lambda \pi(x) \in \mathbb{C}\unit$ for all $x \in \mathfrak{J}$.\smallskip

Two possible cases can be considered here:\smallskip

\noindent [$(1)$] \underline{\textit{Case 1}}: $\{\unit, \pi(x), \pi(u)\}$ are linearly independent. It follows from the previous conclusions and the hypotheses that $T(x) = \lambda \pi(x) + \mu_x\unit$, and thus $T(x) - \lambda \pi(x)$ lies in $\mathbb{C}\unit$.\smallskip

\noindent [$(2)$] \underline{\textit{Case 2}}: $\{\unit, \pi(x), \pi(u)\}$ are linearly dependent. By hypotheses, there exists $z\in \mathfrak{J}$ such that the vectors $\pi (z), \unit = \pi (\unit),$ and $ \pi(u)$ are linearly independent. Therefore for any $\varepsilon > 0$ , the vectors $\{\unit, \pi(u), \pi(x + \varepsilon z)\}$ are linearly independent, and thus, by applying the previous \textit{Case 1} we get $T(x+\varepsilon z) - \lambda \pi(x + \varepsilon z) \in \mathbb{C}\unit$, or equivalently, $T(x) + \varepsilon T(z) - \lambda \pi(x) - \varepsilon \lambda \pi(u) \in \mathbb{C}\unit$ for every $\varepsilon>0$. Taking $\varepsilon \rightarrow 0$ we deduce that $T(x) - \lambda \pi(x) \in \mathbb{C}\unit$.\smallskip

Summarising, $T(x) - \lambda\pi(x) \in \mathbb{C}\unit$ for every $x \in \mathfrak{J}$. Finally, observe that the map $\mu: \mathfrak{J} \rightarrow \mathfrak{M}$ given by $\mu(x):= T(x) - \lambda \pi(x)$ is a linear map such that $\mu(x) \in \mathbb{C}\unit$, and $T(x) = \lambda \pi(x)+ \mu(x)\unit$ for every $x \in \mathfrak{J}$.
\end{proof}

\begin{lemma}\label{l auxiliary lemma for ellinfty sums} Let $p$ be a central projection in a unital JB$^*$-algebra $\mathfrak{J}$. Suppose $T : \mathfrak{J} \rightarrow \mathfrak{J}$ is an associating linear map. Then $p \circ T((\unit-p)\circ x)$ and $(\unit-p) \circ T(p\circ x)$ are central elements for all $x\in \mathfrak{J}.$ 
\end{lemma}

\begin{proof} To simplify the notation, let us denote $q = \unit-p.$ By hypotheses, $[T(x), \mathfrak{J}, x] =0,$ for all $x\in \mathfrak{J}$. So, by linearizing we get $$[T(x), a, y ] =  [x, a, T(y)],$$ for all $x,y,a\in \mathfrak{J}.$ Having this identity in mind and the assumption on $p$ we deduce that 
$$[q\circ T(p \circ x), \mathfrak{J}, y] = q\circ [T(p \circ x), \mathfrak{J}, y] =  q\circ [p \circ x, \mathfrak{J}, T(y)] =  [ q\circ (p \circ x), \mathfrak{J}, T(y)] =0,$$ for all $x,y\in \mathfrak{J},$ witnessing that 
$q\circ T(p \circ x) \in Z(\mathfrak{J}).$ The rest is clear since $q$ also is central.     
\end{proof}

We begin our study with the case of associating linear maps on special JW$^*$-algebras without type $I_1$ part. 

\begin{theorem}\label{t linear associating maps on special JBW*-algebras without type I1} Let $\mathfrak{J}$ be a JW$^*$-algebra admitting no direct summands of type $I_1$,  and let $T : \mathfrak{J} \rightarrow \mathfrak{J}$ be a linear mapping. Suppose additionally that $[T(x), \mathfrak{J}, x] = 0$, for all $x \in \mathfrak{J}$. Then $T$ can be uniquely written in the form: 
    \begin{equation*}
    	T(x) = \lambda\circ x + \mu(x), \ \ \hbox{ for all }x\in \mathfrak{J},
    \end{equation*}
    where $\lambda \in Z (\mathfrak{J})$, and $\mu: \mathfrak{J} \rightarrow Z (\mathfrak{J})$ is linear. 
\end{theorem}

\begin{proof} Let us take an arbitrary Jordan factor representation $\pi: \mathfrak{J} \rightarrow \mathfrak{J}_{\pi}$, where $\mathfrak{J}_{\pi}$ is a factor JW$^*$-algebra which is not of type $I_1$. \smallskip

Since $\pi$ is a Jordan homomorphism, it follows from the hypotheses that the identity 
    $$ [\pi(T(x)), \pi(\mathfrak{J}), \pi(x)] =0,$$  holds for every $x \in \mathfrak{J}$. The weak$^*$-density of $\pi(\mathfrak{J})$ and the separate weak$^*$-continuity of the Jordan product of $\mathfrak{J}_{\pi}$ lead to  
    \begin{equation}\label{eq:lin_asoc}
        [\pi(T(x)), \mathfrak{J}_{\pi}, \pi(x)] =0,
    \end{equation}
    for every $x \in \mathfrak{J}$.\smallskip

Assume first that $\mathfrak{J}_{\pi}$ is a spin factor or a factor JW$^*$-algebra of type $I_2$ (cf. \cite[Theorem 6.1.8]{HOS}), having in mind that $\pi T(x)$ and $\pi (x)$ operator commute in $\mathfrak{J}_{\pi}$ (cf. \eqref{eq:lin_asoc}), \Cref{r associator-commutator of an element in a spin} implies the existence of $ \lambda_{1,x} \lambda_{0,x} \in \mathbb{C}$ such that 
\begin{equation}\label{eq: lin_conm_spin}
    T(x) = \lambda_{1,x} \pi(x) + \lambda_{0,x}\unit
\end{equation}
for all $x \in \mathfrak{J}$. \smallskip

If $\mathfrak{J}_{\pi}$ is not a spin factor, by the same arguments employed in \Cref{prop:JBWspecial}, $\mathfrak{J}_{\pi}$ cannot be a JW$^*$-algebra factor of type $I_1$ or $I_2$, and its universal von Neumann algebra, $W^* (\mathfrak{J}_{\pi}),$ satisfies that $\mathfrak{J}_{\pi} = H(W^* (\mathfrak{J}_{\pi}),\tau)$, where $\tau$ is a period-2 $^*$-anti-automorphism on $W^* (\mathfrak{J}_{\pi})$ and the latter either is a factor von Neumann algebra, or it can be written in the form $W^* (\mathfrak{J}_{\pi}) = N\oplus^{\infty}\tau(N)$, for some weak$^*$-closed ideal $N$ of $W^* (\mathfrak{J}_{\pi})$ which is a factor von Neumann algebra (cf. \Cref{l HWtau factor implies W factor} and the comments prior to it or \cite[Theorem 7.1.9, Remark 7.2.8 and Proposition 7.3.3]{HOS}). As in previous cases, the product of $W^* (\mathfrak{J}_{\pi})$ will be denoted by mere juxtaposition.\smallskip

\Cref{P Topping for non-self adjoint} proves that \eqref{eq:lin_asoc} is equivalent to 
    \begin{equation}\label{eq: lin_comm}
        [\pi(T(x)), \pi(x)]= \pi(T(x))\pi(x) - \pi(x)\pi(T(x)) = 0, \hbox{ for all } x \in \mathfrak{J},
    \end{equation} in terms of the associative product of $ \mathfrak{J}_{\pi} = W^*(\mathfrak{J}_{\pi})$. A linearization argument gives $$[\pi(T(x+y), \pi(x + y)] = [\pi(T(x)),\pi(y)] + [\pi(T(y)),\pi(x)] = 0,$$  for all $x, y \in \mathfrak{J}$.
    Therefore, 
    \begin{equation}\label{eq: id_lin2407}
        [\pi(T(x)), \pi(y)] = - [\pi(T(y)), \pi(x)]
    \end{equation}
    for every $x,y \in \mathfrak{J}$.\smallskip

Now, by combining \eqref{eq: lin_comm} and \eqref{eq: id_lin2407} it can be easily checked that the identity 
    $$ \begin{aligned}
        \relax [\pi(T(x^2)),\pi(y)] & = -[\pi(T(y)),\pi(x^2)] \\ &=   -([\pi(T(y)),\pi(x)]\pi(x) + \pi(x)[\pi(T(y)),\pi(x)])\\
                               & =  [\pi(T(x)),\pi(y)]\pi(x) + \pi(x)[\pi(T(x)),\pi(y)]  \\
                               & = \pi(T(x))\pi(y)\pi(x) - \pi(y)\pi(T(x))\pi(x) \\
                               & + \pi(x)\pi(T(x))\pi(y) - \pi(x)\pi(y)\pi(T(x)),
\end{aligned}$$ holds for every $x, y \in \mathfrak{J}$. Obviously, the previous identity holds when $\pi(y)$ is replaced with an arbitrary $z \in \mathfrak{J}_{\pi}$, and hence  
    \begin{equation}\label{eq: id_linear}
     \begin{aligned}
        &\pi(T(x^2))z\pi(\unit) - \pi(\unit)z\pi(T(x^2)) + z\pi(T(x))\pi(x) \\  & + \pi(x)z\pi(T(x)) - \pi(T(x))z\pi(x) - \pi(x)\pi(T(x))z = 0,
    \end{aligned}
    \end{equation}
    for all $x \in \mathfrak{J}$ and any $z \in \mathfrak{J}_{\pi}$.
Two possibilities arise here:\smallskip

    \begin{enumerate}[$(1)$]
        \item \underline{\emph{Case 1}}: $\{\unit = \pi(\unit), \pi(x)\}$ are $\mathbb{C}$-linearly independent. Having in mind that $\mathfrak{J}_{\pi} = H(W^*(\mathfrak{J}_{\pi}), \tau)$ is a factor, we can combine \Cref{l HWtau factor implies W factor} and \Cref{independencia2} with the identity in \eqref{eq: id_linear}, to conclude that 	
        $$\pi(T(x)) = \lambda_{1,x} \pi(x) + \lambda_{0,x}\unit ,$$
        where $\lambda_{i,x} \in \mathbb{C}$ for $i \in \{0,1\}$.
         \item \underline{\emph{Case 2}}: $\{\mathbf{1}, \pi(x)\} $ are $\mathbb{C}$-linearly dependent. In this case there exist $\alpha\in \mathbb{C}$ such that 
     	$$ \pi(x) = \alpha \pi(1).$$ 
         Then it is clear by \eqref{eq: id_lin2407} that  $$[\pi(T(x)), \pi(y)] = -[\pi(T(y)), \pi(x)] = 0,  \hbox{ for every } y \in \mathfrak{J}.$$ Moreover, due to the weak$^*$-density of $\pi(\mathfrak{J})$ in $\mathfrak{J}_{\pi}$, and the separate weak$^*$-continuity of the product of $W^*(\mathfrak{J_{\pi}})$, the previous identity also holds when $\pi(y)$ is replaced by an arbitrary  $z \in \mathfrak{J}_{\pi}$, and hence $T(x) \in Z(\mathfrak{J}_{\pi}) = \mathbb{C}\unit$. Therefore, there exists $\lambda_{0,x} \in \mathbb{C}$ such that
       \begin{equation*}
              \pi(T(x)) = \lambda_{0,x}\unit.
       \end{equation*}
    \end{enumerate}
    
Summarising, all the previous cases lead to the conclusion that for each $x \in \mathfrak{J}$ there exist complex numbers $\lambda_{0,x}, \lambda_{1,x} $ satisfying
      \begin{equation}\label{eq special B linear}
      	\pi(T(x)) =   \lambda_{1,x}\pi(x) + \lambda_{0,x} \unit.
      \end{equation}     
      
We are now in a position to apply  \Cref{l linear maps from scalars to maps}  to deduce the existence of $\alpha \in \mathbb{C}$ and a linear map $\beta: \mathfrak{J} \rightarrow \mathbb{C}$ satisfying  
      \begin{equation} \label{eq special C}
          \pi(T(x)) =   \alpha \pi(x) + \beta(x) \unit.
      \end{equation}
      
Now observe that due to \Cref{cor:3_lib including spin part} and \Cref{lemma:elem_homo} there exist an element $u \in \mathfrak{J}_{sa}$ and elementary operators $\mathcal{E}_i\in\element (\mathfrak{J}),$ $\widehat{\mathcal{E}_i}\element (\mathfrak{J}_{\pi})$ for all $i \in \{0,1\}$ such that $\widehat{\mathcal{E}_i} \pi = \pi \mathcal{E}_i$, and  $\widehat{\mathcal{E}_i}(\pi(u)^j) = \delta_{ij}\unit,$ for any $ i,j \in \{0,1\}$. Arguing as in the proof of \Cref{prop: uniq_sol}, we deduce from \eqref{eq special C} that $\alpha$ and $ \beta(x)$ can be uniquely written as follows:
\begin{equation}\label{eq expressions for alpha, beta linear} \left\{\begin{aligned}
    \alpha \unit &= \widehat{\mathcal{E}_1}(\pi(T(u))) = \pi(\mathcal{E}_1(T(u)),\\
    \beta(x) \unit &= \widehat{\mathcal{E}_0}(\pi(T(x))) - \widehat{\mathcal{E}_1}(\pi(T(u)))\circ \widehat{\mathcal{E}_0}(\pi(x))\\ & = \pi\Big (\mathcal{E}_0(T(x)) -\mathcal{E}_1(T(u))\circ \mathcal{E}_0(x) \Big),\\
    \end{aligned} \right.
    \end{equation} for all $x\in \mathfrak{J}$. Observe that the Jordan factor representation $\pi$ was arbitrarily chosen, so by setting $\lambda = \mathcal{E}_1(T(u))\in\mathfrak{J}$ and defining a linear mapping $\mu : \mathfrak{J} \rightarrow \mathfrak{J}$ by \begin{equation*}
      \mu(x) = \mathcal{E}_0(T(x)) -\mathcal{E}_1(T(u))\circ \mathcal{E}_0(x),
      \end{equation*} we arrive, via \eqref{eq special C} and \eqref{eq expressions for alpha, beta linear}, to 
     \begin{equation}\label{eq special C with pi}
	\pi(T(x)) = \pi(\lambda \circ x + \mu(x)),
     \end{equation}
    for every $x \in \mathfrak{J}$, and every Jordan factor representation $\pi$ of $\mathfrak{J}$. Since every JB$^*$-algebra admits a separating family of Jordan factor representations of this form (cf. \cite[Corollary 5.7]{AlfShulStor79GelfandNeumark}), it follows from \eqref{eq special C with pi} that \begin{equation}\label{eq C special final without pi}
	T(x) = \lambda \circ x + \mu(x),
     \end{equation}
    for every $x \in \mathfrak{J}$. Note that $\pi(\lambda)$ and $ \pi(\mu(x))$ belong to $Z(\mathfrak{J}_{\pi}) = \mathbb{C} \unit$ for every Jordan factor representation $\pi: \mathfrak{J} \rightarrow \mathfrak{J}_{\pi}$, so \Cref{Rem: central representation} implies that $\lambda$ and $ \mu(x)$ belong to $Z(\mathfrak{J})$ for each $x \in \mathfrak{J}$ as we desired.\smallskip

Observe finally, that the uniqueness of $\lambda\in Z(\mathfrak{J})$ and the linear mapping $\mu$ is essentially guaranteed by the identities $\lambda = \mathcal{E}_1(T(u)),$ and $\mu(x) = \mathcal{E}_0(T(x)) -\mathcal{E}_1(T(u))\circ \mathcal{E}_0(x).$    
\end{proof}

Linear associating maps on an exceptional JBW$^*$-algebra are treated in our next result. 

\begin{theorem}\label{t associating on exceptional} Let $\mathfrak{J}$ denote the  JB$^*$-algebra $C(\Omega, H_3(\mathbb{O}))$, where $\Omega$ is a compact Hausdorff space. Suppose $T: \mathfrak{J}\rightarrow \mathfrak{J}$ is a linear mapping satisfying   $$[T(x), \mathfrak{J}, x] = 0, \hbox{ for all } x \in \mathfrak{J}.$$ Then $T$ can be uniquely written in the form: 
    \begin{equation*}
    	T(x) = \lambda\circ x + \mu(x),  \hbox{ for all } x\in \mathfrak{J},
    \end{equation*}
    where $\lambda \in Z (\mathfrak{J})$, $\mu: \mathfrak{J} \rightarrow Z (\mathfrak{J})$ is linear. 
\end{theorem}

\begin{proof} We keep the notation employed in \Cref{sec: purely exceptional associating trace}. That is, for each $a\in H_3(\mathbb{O})$, we write $\hat{a}$ for the continuous function on $\mathfrak{J}$ with constant value $a$, while for $u \in \mathfrak{J}$ and $t\in \Omega$, we write $u_t := \widehat{u(t)}\in \mathfrak{J}$.\smallskip

The uniqueness can be justified as in the final observation in the proof of the previous \Cref{t linear associating maps on special JBW*-algebras without type I1}.\smallskip

Fix now an arbitrary $t \in \Omega$. Note that  the linear map $T$ can be decomposed as follows: 
    \begin{equation}\label{eq lin_excep_descomp}
        T(u) = T(u-u_t) + T(u_t).
    \end{equation}
Let us deal with each summand on the right hand side of the previous identity independently.\\

\noindent $(1)$ We begin with the term $T(u_t)$. Note that for each $s\in \Omega$ we can define a linear mapping $T_s: H_3(\mathbb{O}) \rightarrow  H_3(\mathbb{O}),$ $T_s(a):= T(\widehat{a})(s)$. We claim that $T_s$ is associating. Indeed, if we take $a,b \in H_3(\mathbb{O})$ we can see from the hypotheses on $T$ that 
$$ [T_s(a),b,a] = [T(\widehat{a}),\widehat{b},\widehat{a}](s) = 0.$$ Thus, by applying \cite[Theorem 2.1]{BreCabVil} we deduce that $T_s$ admits the following representation: 
    \begin{equation}\label{eq lin_excep A}
        T_s(a) = \lambda_s \circ a + \mu_s(a) = \lambda_s \circ a + \mu_s(a)\circ \unit,
    \end{equation}
where $\lambda_s \in Z(H_3(\mathbb{O}))=\mathbb{C}\unit$, and $\mu_s: H_3(\mathbb{O}) \rightarrow Z(H_3(\mathbb{O}))=\mathbb{C}\unit$ is a linear map. Note that this procedure can be done for every $s \in \Omega$. \smallskip

Note now that, by \Cref{cor:3_lib including spin part}, we can find an element $w \in H_3(\Omega)_{sa}$ and elementary operators $\mathcal{E}_0, \mathcal{E}_1\in \element(H_3(\mathbb{O}))$ such that $\mathcal{E}_i(w^j) = \delta_{ij}\unit$  $(i,j \in \{0,1\}).$ \Cref{prop: uniq_sol}, applied to $T_s$ and the Jordan homomorphism $\pi :H_3(\mathbb{O})  \rightarrow H_3(\mathbb{O})  $ given by the identity operator, assures that
\begin{equation}\label{eq: lin_uniq_A}\left\{
	\begin{aligned}
		\lambda_s & = \mathcal{E}_1(T_s(w))\in Z(H_3(\mathbb{O})) = \mathbb{C} \unit,\\ 
        \mu_s(a) &= \mathcal{E}_0(T_s(a)) -\lambda_s \circ \mathcal{E}_0(a) \in Z(H_3(\mathbb{O})) = \mathbb{C} \unit,
		\end{aligned} \right.
\end{equation} for all $a\in H_3(\mathbb{O})$ and all $s\in \Omega$.\smallskip 

Let us consider the natural extension of $w \in H_3(\mathbb{O})_{sa}$ as the constant continuous function $\widehat{w}\in \mathfrak{J}_{sa}$, and the corresponding extension of the elementary operator $\mathcal{E}_i$  $(i \in \{0,1\})$ as an elementary operator $\widetilde{\mathcal{E}_i}$ on $\mathfrak{J} =C(\Omega, H_3(\mathbb{O}))$ ($ i \in \{0,1\}$) given by the corresponding extensions of the involved elements as constant functions. It follows from the definition that $\widetilde{\mathcal{E}_i}(\widehat{w}^j)(s) =\mathcal{E}_i(\widehat{w}^j(s))= \mathcal{E}_i(w^j)= \delta_{i j}\mathbf{1},$ for every $s \in \Omega$. \smallskip

We also set \begin{equation}\label{eq def of lambda mu lin_excep}\left\{
	\begin{aligned}
		\lambda & := \widetilde{\mathcal{E}_1}(T(\widehat{w}))\in \mathfrak{J} = C(\Omega, H_3(\mathbb{O})),\\ 
		\mu (u) &:= \widetilde{\mathcal{E}_0}(T(u)) -\lambda \circ \widetilde{\mathcal{E}_0}(u) \in \mathfrak{J},
	\end{aligned} \right.
\end{equation} for all $u\in\mathfrak{J} = C(\Omega, H_3(\mathbb{O}))$. Clearly, $\mu : \mathfrak{J}\to \mathfrak{J}$ is a linear mapping on $\mathfrak{J}$. \smallskip

Now by applying $\widetilde{\mathcal{E}}_1$ to $T(\widehat{w})$ and having in mind \eqref{eq lin_excep A} we get 
\begin{equation}\label{eq lambda excep A}
    \lambda(s) = \widetilde{\mathcal{E}}_1(T(\widehat{w}))(s) = \mathcal{E}_1(T(\widehat{w})(s)) = \mathcal{E}_1(T_s(w)) = \lambda_s \in \mathbb{C}\unit,
\end{equation}
for every $s \in \Omega$. Observe that $\lambda(s) \in \mathbb{C}\unit$ for every $s \in \Omega$, and thus $\lambda \in Z(\mathfrak{J})$. \smallskip

Fix now two arbitrary elements $s,t \in \Omega$ and $u \in \mathfrak{J}$. By combining \eqref{eq lin_excep A}, \eqref{eq def of lambda mu lin_excep} and \eqref{eq lambda excep A} we derive that  
\begin{equation}\label{eq mu excep A}\left\{ \begin{aligned}
		\mathbb{C} \unit \ni \mu_s(u(t)) &= {\mathcal{E}_0}\Big( T_s( u(t))\Big)   - \lambda(s)\circ {\mathcal{E}_0} \Big(u(t)\Big) \\
		&= {\mathcal{E}_0}(T(u_t)  (s) ) - \lambda(s)\circ {\mathcal{E}_0}( u_t(s))\\
		&=\widetilde{\mathcal{E}_0}(T(u_t)) (s) -\lambda(s)\circ \widetilde{\mathcal{E}_0}(u_t) (s)  \\
		&= \Big(\widetilde{\mathcal{E}_0}(T(u_t))-\lambda\circ \widetilde{\mathcal{E}_0}(u_t)\Big)(s)\\
		& = \mu (u_t) (s),
	\end{aligned}\right.
\end{equation} and consequently, $\mu (u_t)\in Z(\mathfrak{J})$ for all $u\in \mathfrak{J}$, $t\in \Omega$. \smallskip

This shows, via \eqref{eq lin_excep A}, \eqref{eq lambda excep A} and \eqref{eq mu excep A} that the term $T(u_t)$ can be written as follows 
\begin{equation}\label{eq T(ut)}\left\{ \begin{aligned}
    		T(u_t)(s) &= T(\widehat{u(t)}) (s) = T_s (u(t)) \\
    		&= \lambda_s \circ u(t) + \mu_s (u(t))  \unit \\
    		&= \Big(\lambda\circ  u_t + \mu(u_t)\Big)(s).
    	\end{aligned} \right.
    \end{equation}

\noindent $(2)$ We deal next with the term $T(u-u_t)$ in \eqref{eq lin_excep_descomp}. It is not hard to check by the hypotheses on $T$ that for every $a, b \in H_3(\mathbb{O})$ we have 
\begin{align*}
	  [T(u-u_t)(s), a, b] &= [T(u-u_t), \widehat{a}, \widehat{b}](s) =-[T(\widehat{b}), \widehat{a}, u-u_t](s).
\end{align*} When $s = t$, we have $T(u-u_t)(t) \in \mathbb{C} \unit = Z(H_3(\mathbb{O})).$ We can thus define a linear mapping $\beta: \mathfrak{J} \to \mathcal{F} (\Omega,\mathbb{C}\unit)$ given by \begin{equation}\label{eq T(u-ut)} \beta(u)(t) =T(u-u_t) (t) \in \mathbb{C}\unit.
\end{equation}
Back to \eqref{eq lin_excep_descomp}, \eqref{eq T(u-ut)} and \eqref{eq T(ut)} we conclude that 
\begin{equation}\label{eq T(u) (t) final linear exceptional}
T(u)(t) = \lambda(t)\circ u(t) + \Big(\mu(u_t) + \beta(u)\Big)(t),    
\end{equation} for all $u\in \mathfrak{J}$, $t\in \Omega$.  By 
\eqref{eq lambda excep A}, \eqref{eq mu excep A}, and \eqref{eq lin_excep A} we have 
$$\begin{aligned}
\widetilde{\mathcal{E}_0} (T(u))(t) &= \mathcal{E}_0 (T(u)(t)) = \mathcal{E}_0 \Big(\lambda(t)\circ u(t) + \mu(u_t) (t)+ \beta(u)(t) \Big) \\ 
&= \mathcal{E}_0 \Big(\lambda_t \circ u(t)+ \mu_t(u(t)) \Big) + \beta(u)(t) \\
&= \mathcal{E}_0 \Big(T_t (u(t)) \Big) + \beta(u)(t) = \mathcal{E}_0 \Big(T (u_{t}) (t) \Big) + \beta(u)(t) \\
&= \widetilde{\mathcal{E}_0} \Big(T (u_{t}) \Big)  (t) + \beta(u)(t),
\end{aligned} $$ for all $u\in\mathfrak{J}$, $t\in \Omega$, which proves that $\beta (u) = \widetilde{\mathcal{E}_0} (T(u)) - \widetilde{\mathcal{E}_0} \Big(T (u_{t}) \Big)\in Z(\mathfrak{J}).$ Furthermore, back to \eqref{eq def of lambda mu lin_excep} and \eqref{eq T(u) (t) final linear exceptional}, we derive that $$\mu (u) =  \widetilde{\mathcal{E}_0} \Big(T (u) \Big) -\lambda \circ \widetilde{\mathcal{E}_0} (T (u)) = \mu (u_t) + \beta(u), $$ for all $u\in \mathfrak{J},$ which concludes the proof. 
\end{proof}

By the same observations made before \Cref{cor excep bil}, the next result is a straightforward consequence of the representation theory for exceptional JBW$^*$-algebras and the previous  \Cref{t associating on exceptional}.

\begin{corollary}\label{cor associating on exceptional B} Let $\mathfrak{J}$ be a purely exceptional JBW$^*$-algebra. Suppose $T: \mathfrak{J}\rightarrow \mathfrak{J}$ is a linear map satisfying $$[T(x), \mathfrak{J}, x] = 0, \hbox{ for all } x \in \mathfrak{J}.$$ Then $T$ can be uniquely written in the form: 
    \begin{equation*}
    	T(x) = \lambda\circ x + \mu(x), \hbox{ for all } x\in \mathfrak{J},
    \end{equation*}
    where $\lambda \in Z (\mathfrak{J})$, and $\mu: \mathfrak{J} \rightarrow Z (\mathfrak{J})$ is linear. 
\end{corollary}

The main result of this section reads as follows:

\begin{theorem}\label{t linear associating maps on JBW*-algebras without type I1} Let $\mathfrak{J}$ be a JBW$^*$-algebra admitting no direct summands of type $I_1$,  and let $T : \mathfrak{J} \rightarrow \mathfrak{J}$ be a linear mapping. Suppose, additionally, that $[T(x), \mathfrak{J}, x] = 0$, for all $x \in \mathfrak{J}$. Then $T$ can be uniquely written in the form: 
    \begin{equation*}
    	T(x) = \lambda\circ x + \mu(x), \hbox{ for all } x\in \mathfrak{J},
    \end{equation*}
    where $\lambda \in Z (\mathfrak{J})$, and $\mu: \mathfrak{J} \rightarrow Z (\mathfrak{J})$ is linear. 
\end{theorem}

\begin{proof} The uniqueness follows as in the previous two theorems. \smallskip

Standard structure theory guarantees that $\mathfrak{J}$ decomposes as the orthogonal sum $\mathfrak{J} = \mathfrak{J}_1 \oplus^{\infty} \mathfrak{J}_2$ of two weak$^*$-closed ideals $\mathfrak{J}_1$ and $\mathfrak{J}_2$, where $\mathfrak{J}_1$ is a (possibly zero) JW$^*$-algebra admitting no direct summands of type $I_1$, and $\mathfrak{J}_2$ is a (possibly zero) purely exceptional JBW$^*$-algebra \cite[Theorems 5.1.5 and 5.3.5]{HOS}. Let $p$ denote the projection in $\mathfrak{J}$ corresponding to the unit of $\mathfrak{J}_1$. Clearly, $p$ is a central projection in $\mathfrak{J}$, and $q= \unit -p$ is the unit of $\mathfrak{J}_2$ (i.e., $\mathfrak{J}_1 = p\circ \mathfrak{J}$ and $\mathfrak{J}_2 = q\circ \mathfrak{J}$).\smallskip
    
\Cref{l auxiliary lemma for ellinfty sums} implies that the linear mappings $T_{1,2} :=M_p T M_{q}$ and $T_{2,1}:=M_q T M_{p}$ take values in $p\circ Z(\mathfrak{J}) = Z(p\circ \mathfrak{J}) = Z(\mathfrak{J}_1)$ and $ q\circ Z(\mathfrak{J})= Z(q\circ \mathfrak{J}) =Z(\mathfrak{J}_2)$, respectively. \smallskip

By applying that $p$ and $q$ are central projections, it can be easily checked that the linear mappings $T_1 = M_p T|_{\mathfrak{J}_1} : \mathfrak{J}_1\to \mathfrak{J}_1$ and $T_2 = M_q T|_{\mathfrak{J}_2} : \mathfrak{J}_2\to \mathfrak{J}_2$ are associating. We deduce from \Cref{t linear associating maps on special JBW*-algebras without type I1} and \Cref{cor associating on exceptional B} the existence of $\lambda_i \in Z(\mathfrak{J}_i)$ and $\mu_i : \mathfrak{J}_i \to Z(\mathfrak{J}_i)$ linear satisfying 
$$T_i (x) = \lambda_i \circ x + \mu_i (x)\ \ (x\in \mathfrak{J}_i,\  i=1,2).$$ Finally, we can write 
$$\begin{aligned}
T(x) &= p\circ T(p\circ x) + q\circ T(q\circ x) + q\circ T(p\circ x) + p\circ T(q\circ x) \\
&= T_1 (p\circ x) + T_2 (q\circ x) + T_{2,1} (x) + T_{1,2}(x) \\
&= \lambda_1 \circ (p\circ x) + \mu_1 (p\circ x) + \lambda_2 \circ (q\circ x) + \mu_2 (q\circ x) + T_{2,1} (x) + T_{1,2}(x) \\
&= (\lambda_1 + \lambda_2) \circ x + \mu_1 (p\circ x) + \mu_2 (q\circ x) + T_{2,1} (x) + T_{1,2}(x),
\end{aligned}$$ for all $x\in \mathfrak{J}.$ Taking $\lambda = \lambda_1 + \lambda_2\in Z(\mathfrak{J})$, and $\mu : \mathfrak{J}\to Z(\mathfrak{J}) = Z(\mathfrak{J}_1)\oplus^{\infty} Z(\mathfrak{J}_2),$ $\mu(x) := \mu_1 (p\circ x) + \mu_2 (q\circ x) + T_{2,1} (x) + T_{1,2}(x),$ we get the desired decomposition for $T$.
 \end{proof}
 
Determining all associating linear maps on JBW-algebras is a problem of its own importance. We shall see next how we can obtain a general characterization of all associating linear maps on JBW-algebras  from \Cref{t linear associating maps on JBW*-algebras without type I1} above.

\begin{theorem}\label{t linear associating maps on JBW-algebras without type I1} Let $\mathfrak{J}$ be a JBW$^*$-algebra admitting no direct summands of type $I_1$,  and let $T : \mathfrak{J}_{sa} \rightarrow \mathfrak{J}_{sa}$ be a linear mapping. Suppose, additionally, that $$[T(x), \mathfrak{J}_{sa}, x] = 0, \hbox{ for all } x \in \mathfrak{J}_{sa}.$$ Then $T$ can be uniquely written in the form
	\begin{equation*}
		T(x) = \lambda\circ x + \mu(x), \hbox{ for all } x\in \mathfrak{J}_{sa},
	\end{equation*}
	where $\lambda \in Z (\mathfrak{J}_{sa})$, and $\mu: \mathfrak{J}_{sa} \rightarrow Z (\mathfrak{J}_{sa})$ is linear. 
\end{theorem} 

\begin{proof} We can extend $T$ to a complex linear mapping on $\mathfrak{J}$ given by $T_{\mathbb{C}} : \mathfrak{J} \to \mathfrak{J}$, $T_{\mathbb{C}} (a + i b) = T(a) + i T(b)$. We claim that $T_{\mathbb{C}}$ is associating. Namely, since $T$ is associating, we have $$\begin{aligned}
		0= [T(a+ b), \mathfrak{J}_{sa}, a+b] &= [T(a), \mathfrak{J}_{sa}, a] + [T(a), \mathfrak{J}_{sa}, b] \\
		&+[T(b), \mathfrak{J}_{sa}, a] + [T(b), \mathfrak{J}_{sa}, b]\\
		&=  [T(a), \mathfrak{J}_{sa}, b] +[T(b), \mathfrak{J}_{sa}, a],
	\end{aligned}$$  for all $a,b\in \mathfrak{J}_{sa}.$ Therefore, the previous identity shows that 
$$\begin{aligned}
 [T_{\mathbb{C}}(a+ i b), \mathfrak{J}_{sa}, a+i b] &= [T(a), \mathfrak{J}_{sa}, a] - [T(b), \mathfrak{J}_{sa}, b] \\
 &+ i [T(a), \mathfrak{J}_{sa}, b] + i [T( b), \mathfrak{J}_{sa}, a] =0 ,
\end{aligned}$$  for all $a+i b\in \mathfrak{J}$, and thus  $ [T_{\mathbb{C}}(a+ i b), \mathfrak{J}, a+i b] = 0,$ which proves the claim. \smallskip

Now, \Cref{t linear associating maps on JBW-algebras without type I1} applied to $T_{\mathbb{C}}$ assures the existence of a unique $\gamma \in Z (\mathfrak{J})$, and a unique linear mapping $\nu: \mathfrak{J} \rightarrow Z (\mathfrak{J})$ such that $$T_{\mathbb{C}} (a+ i b) = \gamma\circ (a + i b) + \nu (a+ i b ), \hbox{ for all } a+i b \in \mathfrak{J}.$$

Observe that, by construction, $T_{\mathbb{C}} (a + i b)^* = T_{\mathbb{C}} \left( (a+ i b)^* \right)$ for all $a+ i b \in \mathfrak{J}.$ It then follows that $$T_{\mathbb{C}} (a+ i b) = T_{\mathbb{C}} ((a +i b)^*)^* = \gamma^*\circ (a + i b) + \nu((a+ i b)^* )^*, \hbox{ for all } a+i b \in \mathfrak{J}.$$ The uniqueness of the decomposition implies that $\gamma^* = \gamma$ lies in $\mathfrak{J}_{sa}$ and $\nu (a + i b )^* = \nu ((a+i b)^*)$, for all $a+i b \in \mathfrak{J}.$	The mapping $\mu := \nu|_{\mathfrak{J}}$ takes values in $Z(\mathfrak{J}_{sa}) = Z(\mathfrak{J})_{sa}$, and taking $\lambda = \gamma,$ we get  $$T(x) = \lambda\circ x + \mu(x), \hbox{ for all } x\in \mathfrak{J}_{sa}.$$
\end{proof}

\section{Conclusions on associating traces on general \texorpdfstring{JBW$^*$}{JBW*}-algebras}\label{sec: bilinear associating maps on general JBW*-algebras}

As the reader can already guess by the title, the goal of this section is to gather all previous results to determine the form of all symmetric bilinear mappings on a general JBW$^*$-algebra whose trace is associating.\smallskip 

Our first technical result is a bilinear version of the previous \Cref{l auxiliary lemma for ellinfty sums}.

\begin{lemma}\label{l bilinear associting trace and central projections bilinear version of auxiliary lemma} Let $p$ be a central projection in a unital JB$^*$-algebra $\mathfrak{J}$, and let $B:\mathfrak{J}\times\mathfrak{J}\rightarrow \mathfrak{J}$ be a symmetric bilinear map satisfying $[B(x,x),\mathfrak{J},x]=0$ for all $x \in \mathfrak{J}$. Then the mappings $$x\mapsto p\circ B((\unit-p)\circ x,(\unit-p)\circ x), (\unit-p) \circ B(p\circ x,p\circ x) \  \ \ (x\in \mathfrak{J}),$$ are center-valued. 
\end{lemma}

\begin{proof} Set $q = \unit-p$, which is also central. Let us take $x, z, y \in \mathfrak{J}$ then by \eqref{eq new remark 1206 Jordan} in \Cref{r identity Bresar} we have 
        \begin{align*}
            [p\circ B(q\circ x,q\circ x), y,z] & = p\circ [B(q\circ x,q\circ x), y,z]\\ & = p\circ (-2 [B(q\circ x,z), y, q\circ x])= 0,
        \end{align*} where the last equality holds because $p$ is central. The remaining statement follows by similar arguments. 
\end{proof}

The next technical conclusion is perhaps interesting by itself. 

\begin{proposition}\label{p mixed products} Let  
$\mathfrak{J}$ be a JBW$^*$-algebra	without direct summands of type $I_1$, and let $p$ be a central projection in $\mathfrak{J}$. Suppose that $B: \mathfrak{J} \times \mathfrak{J} \rightarrow \mathfrak{J}$ is a symmetric bilinear mapping satisfying  $[B(x,x),\mathfrak{J},x]=0,$ for all $x \in \mathfrak{J}$. Then there exist a unique linear mapping $\mu: p\circ \mathfrak{J}\to p\circ Z(\mathfrak{J})=   Z(p\circ \mathfrak{J}),$ and a unique bilinear mapping $\nu : p\circ \mathfrak{J}\times q\circ \mathfrak{J}\to p\circ  Z(\mathfrak{J}) =   Z(p\circ \mathfrak{J})$ satisfying 
$$p\circ B(p\circ x,q\circ x) = \mu (p\circ 
 x)\circ  x + \nu(p\circ x, q\circ x),$$ for all $x\in \mathfrak{J}.$
\end{proposition}

\begin{proof} Set $q = \unit-p$, which is also central. Let us fix and arbitrary $x \in \mathfrak{J}$ and define a linear mapping $T_{p,x}: p\circ \mathfrak{J} \rightarrow p\circ \mathfrak{J} $ given by $a \mapsto T_{p,x} (a)=p\circ B(a, q\circ x)$. The assumptions combined with \eqref{eq new remark 1206 Jordan} in \Cref{r identity Bresar}, give
        \begin{align*}
            [T_{p,x}(a), b,a] &= [p\circ B(a, q\circ x), b, a] = p\circ [B(a,q\circ x), b, a]\\ &= -\frac{1}{2}p\circ [B(a,a), b , q\circ x] = 0,
        \end{align*} for all $b\in\mathfrak{J}$. Consequently, $T_{p,x}$ is an associating linear map on $p\circ \mathfrak{J}$, and the latter cannot contain summands of type $I_1$. Therefore, by applying \Cref{t linear associating maps on JBW*-algebras without type I1}, there exist a unique element $\lambda_{p,x} \in Z(p\circ \mathfrak{J})$ and a unique linear mapping $\mu_{p, x}: p\circ \mathfrak{J} \rightarrow Z(p\circ \mathfrak{J})$ such that 
        \begin{equation}\label{eq 2907 1344}
            T_{p,x} (p\circ y) = \lambda_{p,x} \circ (p \circ y) +\mu_{p,x}(p\circ y),
        \end{equation}
        for all $y \in \mathfrak{J}$.\smallskip 
        
Now, let us define a linear mapping $\mu: p\circ \mathfrak{J} \rightarrow p\circ \mathfrak{J},$ and a symmetric bilinear map $\nu: p\circ \mathfrak{J} \times p\circ \mathfrak{J} \rightarrow p\circ \mathfrak{J}$ given by
        \begin{equation}\label{eq definition of mu and nu in section 7 2907}
            \left\{ 
            \begin{aligned}
                \mu(y) &= \mathcal{E}_1(p\circ B(p\circ u, q\circ y)) =p\circ  \mathcal{E}_1(B(p\circ u, q\circ y)),\\
                \nu(y,y) &= \mathcal{E}_0(p\circ B(p\circ y, q\circ y)) - \mu(y)\circ \mathcal{E}_0(p\circ y)\\
                &= p\circ \mathcal{E}_0( B(p\circ y, q\circ y)) - \mu(y)\circ \mathcal{E}_0(p\circ y),
            \end{aligned}
            \right.
        \end{equation}
        where $\mathcal{E}_i\in \element (\mathfrak{J})$ ($i=0,1$) and $u \in \mathfrak{J}$ are the elementary operators on $\mathfrak{J}$ and the element whose existence is guaranteed by \Cref{cor:3_lib including spin part}. Now observe that, by applying $\mathcal{E}_1$ to $T(p\circ u)$ in \eqref{eq 2907 1344}, we obtain 
        \begin{equation}\label{eq mu(x) and mux,p(x)} \begin{aligned}
        \mu (x) = \mathcal{E}_1(p\circ B(p\circ u, q\circ x)) &= \mathcal{E}_1(T(p\circ u)) = \mathcal{E}_1( \lambda_{p,x} \circ (p \circ u) ) \\
        &= p \circ \mathcal{E}_1( \lambda_{p,x} \circ u ) = \lambda_{p,x} \in Z(p\circ \mathfrak{J}).
        \end{aligned}
        \end{equation}
        Thus, the arbitrariness of $x\in \mathfrak{J}$ implies that $\mu$ is a linear mapping whose image is in $Z(p\circ \mathfrak{J})$.\smallskip

On the other hand, having in mind \eqref{eq mu(x) and mux,p(x)} and \eqref{eq 2907 1344} we derive: 
$$\begin{aligned}
\nu (x,x) &= \mathcal{E}_0\Big( p\circ B(p\circ x, q\circ x)) - \mu(x)\circ (p\circ x)\Big) \\
&=  \mathcal{E}_0\Big( T_{p,x} (p\circ x) - \lambda_{x,u} \circ (p\circ x)\Big)  = \mu_{p,x} (p\circ x)   \in Z(p\circ \mathfrak{J}),
\end{aligned}$$ which shows that $\nu$ is center-valued. By construction 
$$\begin{aligned}
p\circ B(p\circ x, q\circ x) &= T_{p,x} (p\circ x) = \lambda_{p,x} \circ (p \circ x) +\mu_{p,x}(p\circ x)\\
&= \mu(x)\circ (p\circ x) + \nu(x,x) = \mu(x)\circ x + \nu(x,x),
\end{aligned}$$ for all $x \in \mathfrak{J}$.\smallskip

The uniqueness of $\mu$ and $\nu$ has been implicitly proved by the uniqueness of the expressions in \eqref{eq definition of mu and nu in section 7 2907} in terms of $\mathcal{E}_i$ and the mapping $p\circ B(p\circ x, q\circ x)$ ($x\in \mathfrak{J}$).
\end{proof}

For the sake of simplicity, before stating our next technical tool we introduce some notation. We shall say that a JB$^*$-algebra $\mathfrak{J}$ satisfies the \emph{standard factorization property for associating traces} (SFP in short) if for every symmetric bilinear mapping $B: \mathfrak{J} \times \mathfrak{J}\to \mathfrak{J}$ with associating trace there exist a unique $\lambda \in Z (\mathfrak{J})$, a unique linear mapping $\mu: \mathfrak{J}\rightarrow Z (\mathfrak{J}),$ and a unique bilinear mapping $\nu : \mathfrak{J}\times \mathfrak{J} \rightarrow Z (\mathfrak{J})$ such that $$B(x,x) = \lambda \circ x^2 + \mu(x) \circ x + \nu(x,x).$$ for all $x\in \mathfrak{J}$. 

\begin{proposition}\label{p direct sums of factors with the SFP} Let $p$ be a central projection in a JBW$^*$-algebra $\mathfrak{J}$ without direct summands of type $I_1$. Suppose that $p\circ \mathfrak{J}$ and $(\unit-p)\circ \mathfrak{J}$ satisfy the SFP. Then $\mathfrak{J}$ satisfies the SFP.
\end{proposition}

\begin{proof} Let $B: \mathfrak{J} \times \mathfrak{J} \rightarrow \mathfrak{J}$ be a symmetric bilinear map whose trace is associating, that is, $[B(x,x),\mathfrak{J},x]=0,$ for every $x \in \mathfrak{J}$. Let us decompose $B(x,x)$ in the from \begin{equation}\label{eq general form of B(x,x) in proposition ell sums 7.3}\left\{
   \begin{aligned}
        B(x,x) 
                =& p\circ B(p\circ x, p\circ x) + 2p\circ B(p\circ x, q\circ x) + p\circ B(q\circ x,q\circ x)  \\
               & + q\circ B(p\circ x, p\circ x) + 2q\circ B(p\circ x, q\circ x) + q\circ B(q\circ x,q\circ x), 
    \end{aligned} \right.
\end{equation} with $x\in \mathfrak{J}$, and write $B_p : p\circ \mathfrak{J} \times p\circ \mathfrak{J} \rightarrow p\circ \mathfrak{J}$, $B_q : q\circ \mathfrak{J} \times q\circ \mathfrak{J} \rightarrow q\circ \mathfrak{J}$ for the symmetric bilinear mappings given by $$B_p(p\circ x,p\circ x) = p\circ B(p\circ x, p\circ x),\hbox{ and } B_q(q\circ x,q\circ x) = q\circ B(q\circ x, q\circ x),$$ respectively.\smallskip 

It is not hard to see, from the fact that $p$ and $q$ are central projections and the assumptions on $B$, that the traces of the bilinear maps $B_p$ and $B_q$ are associating. Since $p\circ \mathfrak{J}$ and $(\unit-p)\circ \mathfrak{J}$ satisfy the SFP, there exist unique $\lambda_p \in Z (p\circ \mathfrak{J})$, $\lambda_q \in Z (q\circ \mathfrak{J})$, linear mappings $\mu_p: p\circ \mathfrak{J}\rightarrow Z (p\circ \mathfrak{J}),$ $\mu_q: q\circ \mathfrak{J}\rightarrow Z (q\circ \mathfrak{J}),$ and bilinear mappings $\nu_{p} : p\circ \mathfrak{J}\times p\circ \mathfrak{J} \rightarrow Z (p\circ \mathfrak{J})$ and $\nu_{q} : q\circ \mathfrak{J}\times q\circ \mathfrak{J} \rightarrow Z (q\circ \mathfrak{J})$ such that 
\begin{equation}\label{eq expression for Bp} B_p(p\circ x,p\circ x) = \lambda_p \circ x^2 + \mu_p(p\circ x) \circ x + \nu_p (p \circ x,p\circ x),
\end{equation} and 
\begin{equation}\label{eq expression for Bq} 
B_q(q\circ x,q\circ x) = \lambda_q \circ x^2 + \mu_q(q\circ x) \circ x + \nu_q (q \circ x,q\circ x), 
\end{equation} for all $x\in \mathfrak{J}$. \smallskip

We shall also employ the following maps: $$\nu_{p,q}: q\circ \mathfrak{J} \times q\circ \mathfrak{J} \longrightarrow p\circ \mathfrak{J}, \ \ \ \ \ \ \ \ \  \nu_{q,p}: p\circ \mathfrak{J} \times p\circ \mathfrak{J} \longrightarrow q\circ \mathfrak{J} $$
    $$\nu_{p,q} (q\circ x, q\circ x)= p\circ B(q\circ x,q\circ x), \ 
    \nu_{q,p} (p\circ x, p\circ x)= q\circ B(p\circ x,p\circ x).$$ \Cref{l bilinear associting trace and central projections bilinear version of auxiliary lemma} assures that the above two mappings are center-valued, that is,  \begin{equation}\label{eq nus are center-valued} 
    \nu_{p,q}\left( q\circ \mathfrak{J} \times q\circ \mathfrak{J}\right) \subseteq  Z(p\circ \mathfrak{J}) = p\circ   Z(\mathfrak{J}), \hbox{ and  } \nu_{q,p}:\left( p\circ \mathfrak{J} \times p\circ \mathfrak{J}\right)  \subseteq  Z(q\circ \mathfrak{J}).
    \end{equation}
    
Consider now the mappings $$\delta_{p,p,q}: p\circ \mathfrak{J} \times q\circ \mathfrak{J} \longrightarrow p\circ \mathfrak{J}, \ \ \ \ \ \ \ \ \  \delta_{q,p,q}: p\circ \mathfrak{J} \times q\circ \mathfrak{J} \longrightarrow q\circ \mathfrak{J} $$
    $$\delta_{p,p,q} (p\circ x, q\circ x)= 2p\circ B(p\circ x, q\circ x), \ 
    \delta_{q,p,q} (p\circ x, q\circ x)= 2q\circ B(p\circ x, q\circ x).$$ Having in mind that $\mathfrak{J}$ admits no type $I_1$ central summands, we are in a position to apply \Cref{p mixed products} to deduce the existence of unique linear mappings $\mu_{p,p,q}: p\circ \mathfrak{J}\to p\circ Z(\mathfrak{J})=   Z(p\circ \mathfrak{J}),$ $\mu_{q,p,q}: q\circ \mathfrak{J}\to q\circ Z(\mathfrak{J})=   Z(q\circ \mathfrak{J})$, and unique bilinear mappings $\nu_{p,p,q} : p\circ \mathfrak{J}\times q\circ\mathfrak{J}\to p\circ  Z(\mathfrak{J}) =   Z(p\circ \mathfrak{J})$, and $\nu_{q,p,q} : p\circ \mathfrak{J}\times q\circ\mathfrak{J}\to q\circ  Z(\mathfrak{J}) =   Z(q\circ \mathfrak{J})$ satisfying
\begin{equation}\label{eq expression deltappq} 
\delta_{p,p,q}(p\circ x, q\circ x) = p\circ B(p\circ x,q\circ x) = \mu_{p,p,q} (p\circ 
 x)\circ  x + \nu_{p,p,q} (p\circ x, q\circ x),
\end{equation} and 
\begin{equation}\label{eq expression deltaqpq} 
\delta_{q,p,q}(p\circ x, q\circ x) = q\circ B(p\circ x,q\circ x) = \mu_{q,p,q} (q\circ 
 x)\circ  x + \nu_{q,p,q} (p\circ x, q\circ x),
\end{equation} for all $x\in \mathfrak{J}.$ Mixing \eqref{eq general form of B(x,x) in proposition ell sums 7.3}, \eqref{eq expression for Bp}, \eqref{eq expression for Bq}, \eqref{eq expression deltappq}, and \eqref{eq expression deltaqpq} we derive that $$\begin{aligned}
B(x,x) &= \lambda_p \circ x^2 +  \lambda_q \circ x^2 + \mu_p(p\circ x) \circ x   + \mu_q(q\circ x) \circ x + \mu_{p,p,q} (p\circ 
 x)\circ  x  \\  &+ \mu_{q,p,q} (q\circ 
 x)\circ  x + \nu_p (p \circ x,p\circ x) + \nu_q (q \circ x,q\circ x) + \nu_{p,q} (q\circ x, q\circ x) \\
 &+ \nu_{q,p} (p\circ x, p\circ x)  + \nu_{p,p,q} (p\circ x, q\circ x) + \nu_{q,p,q} (p\circ x, q\circ x) \\
 =& (\lambda_p +  \lambda_q) \circ x^2 + (\mu_p(p\circ x) + \mu_q(q\circ x) + \mu_{p,p,q} (p\circ 
 x) + \mu_{q,p,q} (q\circ 
 x)) \circ  x \\
&+ \nu_p (p \circ x,p\circ x) + \nu_q (q \circ x,q\circ x) + \nu_{p,q} (q\circ x, q\circ x) \\
 &+ \nu_{q,p} (p\circ x, p\circ x)  + \nu_{p,p,q} (p\circ x, q\circ x) + \nu_{q,p,q} (p\circ x, q\circ x),
\end{aligned}$$ for all $x\in \mathfrak{J}$. 
Setting $\lambda = \lambda_p +  \lambda_q\in Z(\mathfrak{J}),$ $\mu: \mathfrak{J}\to Z(\mathfrak{J}),$ $$\mu (x):= \mu_p(p\circ x) + \mu_q(q\circ x) + \mu_{p,p,q} (p\circ 
 x) + \mu_{q,p,q} (q\circ 
 x),$$ and $\nu : \mathfrak{J}\times \mathfrak{J}\to Z(\mathfrak{J}),$ $$\begin{aligned}
  \nu (x,x) &:= \nu_p (p \circ x,p\circ x) + \nu_q (q \circ x,q\circ x) + \nu_{p,q} (q\circ x, q\circ x) \\
  &+ \nu_{q,p} (p\circ x, p\circ x)  + \nu_{p,p,q} (p\circ x, q\circ x) + \nu_{q,p,q} (p\circ x, q\circ x),   
 \end{aligned}$$ we establish the existence of the desired factorization for $B(x,x)$ in the form $$B(x,x) = \lambda x^2 + \mu(x) \circ x + \nu (x,x) \ \ (x\in \mathfrak{J}).$$

The uniqueness deserves to be commented at least briefly. It essentially follows from the uniqueness in the SFP and in \Cref{p mixed products}. Namely, suppose we can also write $$B(x,x) = \widetilde{\lambda} x^2 + \widetilde{\mu}(x) \circ x + \widetilde{\nu} (x,x) \ \ (x\in \mathfrak{J}),$$ where $\widetilde{\lambda},$ $\widetilde{\mu},$ and $\widetilde{\nu}$ satisfy the required properties. Clearly 
$$ \begin{aligned}
(p \circ \widetilde{\lambda})& (p\circ x)^2 + (p\circ \widetilde{\mu}(p\circ x)) \circ (p\circ x) + p\circ \widetilde{\nu} (p\circ x,p\circ x) = p\circ B(p\circ x, p\circ x) \\
&= (p \circ {\lambda}) (p\circ x)^2 + (p\circ {\mu}(p\circ x)) \circ (p\circ x) + p\circ {\nu} (p\circ x,p\circ x) 
\end{aligned},$$ for all $x\in \mathfrak{J}$, and thus, the SFP of $p\circ \mathfrak{J}$ assures that \begin{equation}
\left\{\begin{aligned}\label{eq first identities towards uniqueness 1 3007}
p \circ \widetilde{\lambda} &= p \circ {\lambda}, \  p\circ \widetilde{\mu}(p\circ x) =p\circ {\mu}(p\circ x),\\ 
\hbox{  and } & p\circ \widetilde{\nu} (p\circ x,p\circ x) = p\circ {\nu} (p\circ x,p\circ x),    
\end{aligned} \right.
\end{equation} for all $x\in \mathfrak{J}.$
Similarly, \begin{equation}\label{eq first identities towards uniqueness 2 3007}
\left\{\begin{aligned}
q \circ \widetilde{\lambda} &= q \circ {\lambda}, \ \ q\circ \widetilde{\mu}(q\circ x)= q\circ {\mu}(q\circ x),\\
\hbox{ and } & q\circ \widetilde{\nu} (q\circ x,q\circ x) = q\circ {\nu} (q\circ x,q\circ x),
\end{aligned} \right.    
\end{equation} for all $x\in \mathfrak{J}.$ Consquently, $\lambda = \lambda\circ (p+q) = \widetilde{\lambda} \circ (p+q)  =  \widetilde{\lambda}.$ \smallskip

On the other hand, 
$$\begin{aligned}
p\circ B(p\circ x,q\circ x) &=\frac12 p\circ\Big(
B(x, x)- B(p\circ x,p\circ x) - B(q\circ x,q\circ x)  \Big) \\ 
&= \frac12 \Big( (p\circ \lambda)\circ (p \circ x)^2+ (p\circ \mu (x))\circ (p\circ x) + p\circ \nu (x,x) \\
&- (p\circ \lambda)\circ (p \circ x)^2 -  (p\circ \mu (p\circ x))\circ (p\circ x) - p\circ \nu (p\circ x,p \circ x) \\
&- p\circ \nu (q\circ x,q \circ x)  \Big)\\
&=\frac12 \Big( (p\circ \mu (q\circ x))\circ (p\circ x) + 2 p\circ \nu (p\circ x,q \circ x)\Big),
\end{aligned}$$ and similarly 
$$p\circ B(p\circ x,q\circ x)= \frac12 \Big( (p\circ \widetilde{\mu} (q\circ x))\circ (p\circ x) +2 p\circ \widetilde{\nu} (p\circ x,q \circ x)\Big),$$ for all $x\in \mathfrak{J}.$ The uniqueness of the decomposition in \Cref{p mixed products} implies that \begin{equation}\label{eq first identities towards uniqueness 3 3007}
(p\circ \mu (q\circ x)) = (p\circ \widetilde{\mu} (q\circ x)), \hbox{ and } p\circ \widetilde{\nu} (p\circ x,q \circ x) = p\circ {\nu} (p\circ x,q \circ x),    
\end{equation}
 for all $x\in \mathfrak{J}.$ \smallskip

Since the roles of $p$ and $q$ are clearly symmetric, we also get \begin{equation}\label{eq first identities towards uniqueness 4 3007} 
(q\circ \mu (p\circ x)) = (q\circ \widetilde{\mu} (p\circ x)), \hbox{ and } q\circ \widetilde{\nu} (q\circ x,p \circ x) = q\circ {\nu} (q\circ x,p \circ x),    
\end{equation}
for all $x\in \mathfrak{J}.$ We can easily deduce from \eqref{eq first identities towards uniqueness 1 3007}, \eqref{eq first identities towards uniqueness 2 3007}, \eqref{eq first identities towards uniqueness 3 3007}, and \eqref{eq first identities towards uniqueness 4 3007} that $\mu = \widetilde{\mu}$.\smallskip

Finally, the identities 
$$p\circ \widetilde{\nu} (q\circ x,q \circ x) = p\circ B(q\circ x, q\circ x) = p\circ {\nu} (q\circ x,q \circ x),$$ and 
$$q\circ \widetilde{\nu} (p\circ x,p \circ x) = q\circ B(p\circ x, p\circ x) = q\circ {\nu} (p\circ x,p \circ x), \ (x\in \mathfrak{J}),$$ combined with the remaining information from \eqref{eq first identities towards uniqueness 1 3007}, \eqref{eq first identities towards uniqueness 2 3007}, \eqref{eq first identities towards uniqueness 3 3007}, and \eqref{eq first identities towards uniqueness 4 3007} give  
$\nu = \widetilde{\nu}$.
\end{proof}

The following extension of the previous result follows from a simple induction argument. 

\begin{proposition}\label{p direct sums of n factors with the SFP} Let $p_1,\ldots,p_n$ ($n\in \mathbb{N}$) be central projections in a JBW$^*$-algebra $\mathfrak{J}$ without direct summands of type $I_1$. Suppose that $p_1+\ldots+p_n=\unit $ and  $p_j \circ \mathfrak{J}$ satisfies the SFP for all $j=1,\ldots,n.$ Then $\mathfrak{J}$ satisfies the SFP.
\end{proposition}

\begin{proof} In order to give a brief idea of the induction argument, suppose that the statement is true for $n$, and let us consider $n+1$ central projections $p_1,\ldots,p_n,p_{n+1}$ such that  $p_1+\ldots + p_n +p_{n+1} = \unit$ and $p_j \circ \mathfrak{J}$ satisfies the SFP for all $j.$ We deduce, from the induction hypothesis, that taking $p= p_1+\ldots +p_n$ we get a central projection such that $p\circ \mathfrak{J}$ and $(\unit-p)\circ \mathfrak{J} = p_{n+1}\circ \mathfrak{J}$ satisfy the SFP. \Cref{p direct sums of factors with the SFP} assures that $\mathfrak{J}$ satisfies the SFP, which concludes the induction argument.   
\end{proof}

We can now proceed with the description of all associating traces on a JBW$^*$-algebra without commutative summands. 

\begin{theorem}\label{teo bilinear mapping JBW* algebra} Let $\mathfrak{J}$ be a JBW$^*$-algebra with no direct summands of type $I_1$, and let $B: \mathfrak{J} \times \mathfrak{J} \rightarrow \mathfrak{J}$ be a symmetric bilinear map satisfying $[B(x,x),\mathfrak{J},x]=0,$ for all $x \in \mathrm{J}$. Then $B$ admits the following (unique) representation 
    \begin{equation*}
        B(x,x) = \lambda \circ x^2 + \mu(x) \circ x + \nu(x,x)\circ \mathbf{1}
    \end{equation*}
    where $\lambda \in Z (\mathfrak{J})$, $\mu: \mathfrak{J}\rightarrow Z (\mathfrak{J})$ is linear and $\nu : \mathfrak{J}\times \mathfrak{J} \rightarrow Z (\mathfrak{J})$ is bilinear. If $\mathfrak{J}$ is a JBW$^*$-algebra of type $I_2$ the element $\lambda$ is always zero. 
\end{theorem}

\begin{proof} We appeal, once again, to the general structure theory of JBW$^*$-algebras to assure the existence of three (possibly zero) central projections $p_1,p_2$ and $p_3$ such that $p_1 + p_2 + p_3 = \unit$, $p_1 \circ \mathfrak{J}$ is a JW$^*$-algebra with no direct summands of type $I_1$ and $I_2$, $p_2 \circ \mathfrak{J}$ is a JW$^*$-algebra of type $I_2$, and $p_3 \circ \mathfrak{J}$ is a purely exceptional JBW$^*$-algebra (cf. \cite[Theorems 5.1.5, 5.3.5 and 7.2.7]{HOS}). The JBW$^*$-algebras $p_1 \circ \mathfrak{J}$, $p_2 \circ \mathfrak{J}$, and $p_3 \circ \mathfrak{J}$ satisfy the SFP by \Cref{prop:JBWspecial}, \Cref{t bilinear associating trace on spin}, and \Cref{cor excep bil}, respectively. \Cref{p direct sums of n factors with the SFP} now asserts that $\mathfrak{J}$ also satisfies the SFP, which is equivalent to the desired statement. The final claim is a consequence of \Cref{t bilinear associating trace on spin}.
\end{proof}

\section{Applications to preservers of operator commutativity}\label{sec: applications preservers of operator commutativity}

Our previous studies on bilinear maps with associating trace will be now applied to study linear bijections between JBW$^*$-algebras preserving operator commutativity among elements. \smallskip

The result will be obtained after a series of technical lemmata. Our first stop will be a generalization of \cite[Lemma 5]{BreMei1993}, whose proof is new even in the associative setting.

\begin{lemma}\label{l Jordan version of Brsar Meiers central} Let $\mathfrak{J}$ be a JB$^*$-algebra with no  one-dimensional Jordan  representations. If $c\in Z(\mathfrak{J})$ satisfies $c\circ \mathfrak{J}\subseteq Z(\mathfrak{J})$, then $c = 0$. The conclusion also holds if $\mathfrak{J}$ is a JBW$^*$-algebra with no type $I_1$ direct summand. 
\end{lemma}

\begin{proof} Let $\pi : \mathfrak{J} \rightarrow \mathfrak{J}_{\pi}$ be a Jordan factor representation. Since $[\pi (c), \pi(\mathfrak{J}), \pi (a)] = \pi [c, \mathfrak{J}, a] =0,$ for all $a\in \mathfrak{J}$, it follows from the weak$^*$-density of $\pi (\mathfrak{J})$ in the factor JBW$^*$-algebra $\mathfrak{J}_{\pi}$ and the separate weak$^*$-continuity of the Jordan product of the latter JBW$^*$-algebra that $\pi (c)\in Z(\mathfrak{J}_{\pi}) = \mathbb{C}\unit$.   Therefore there exists a complex $\lambda_{\pi}$ such that $\pi (c) = \lambda_{\pi} \unit$, and thus, by hypothesis, $\lambda_{\pi} \pi (\mathfrak{J})= \pi (c) \circ \pi (\mathfrak{J}) \subseteq \pi (Z(\mathfrak{J})) \subseteq Z(\mathfrak{J}_{\pi}) = \mathbb{C}\unit,$ which combined with the weak$^*$-density of $\pi (\mathfrak{J})$ shows that $\lambda_{\pi}  \mathfrak{J}_{\pi} \subseteq \mathbb{C}\unit$. By observing that $\mathfrak{J}$ admits no  one-dimensional Jordan  representations, we conclude that $\lambda_{\pi}=0$. We have therefore shown that $\pi (c) =0$ for every Jordan factor representation $\pi : \mathfrak{J} \rightarrow \mathfrak{J}_{\pi}$. The existence a faithful family of Jordan factor representations for $\mathfrak{J}$ (cf. \cite[Corollary 5.7]{AlfShulStor79GelfandNeumark}) assures that $c =0$, as desired. \smallskip

For the final statement, note that, by \Cref{cor:3_lib including spin part} and \Cref{lemma:elem_homo}, a JBW$^*$-algebra with no type $I_1$ direct summand  admits no  one-dimensional Jordan representations.  
\end{proof}

It is known that for each JB$^*$-algebra $\mathfrak{J}$, the centre of the JB-algebra  $\mathfrak{J}_{sa}$ coincides with the self-adjoint part of the centre of $\mathfrak{J}$. Therefore the next result is a direct consequence of the previous \Cref{l Jordan version of Brsar Meiers central}.

\begin{lemma}\label{l Jordan version of Brsar Meiers central for JB-algebras} Let $\mathfrak{J}$ be a JB$^*$-algebra with no  one-dimensional Jordan representations. If $c\in Z(\mathfrak{J}_{sa})$ satisfies $c\circ \mathfrak{J}_{sa}\subseteq Z(\mathfrak{J}_{sa})$, then $c = 0$. The conclusion also holds if $\mathfrak{J}$ is a JBW$^*$-algebra with no type $I_1$ direct summands. 
\end{lemma}

The next result is inspired by \cite[Proof of Theorem 3]{BreMei1993}. The proof, which is included here for completeness reasons, relies on \Cref{l Jordan version of Brsar Meiers central} and \Cref{t linear associating maps on JBW*-algebras without type I1}.

\begin{lemma}\label{l second technical lemma in last section} Let $\mathfrak{M}$ and $\mathfrak{J}$ be JBW$^*$-algebras with no central summands of type $I_1$. Suppose $\Phi: \mathfrak{M} \rightarrow \mathfrak{J}$ is a linear bijection preserving operator commutativity in both directions, that is, $[x,\mathfrak{M},y] = 0$  if, and only if, $[\Phi(x),\mathfrak{J},\Phi(y)] = 0,$ for all $x,y\in \mathfrak{M}$. Then $\Phi$ maps $Z(\mathfrak{M})$ onto $Z(\mathfrak{J}),$ and there exists an associative isomorphism $\alpha: Z(\mathfrak{M}) \to Z(\mathfrak{J})$ satisfying $\Phi (z \circ a) - \alpha (z) \circ \Phi (a)\in Z(\mathfrak{J}),$ for all $z\in Z(\mathfrak{M})$, $a\in \mathfrak{M},$ \ and $\Phi^{-1} (\tilde{z} \circ b) - \alpha^{-1} (\tilde{z}) \circ \Phi^{-1} (b)\in Z(\mathfrak{M}),$ for all $\tilde{z}\in Z(\mathfrak{J})$, $b\in \mathfrak{J}$.   
\end{lemma}

\begin{proof} Let us fix an arbitrary $z\in Z(\mathfrak{M})$ and consider the mapping $\Phi_{z}: \mathfrak{J}\to \mathfrak{J}$ given by $\Phi_{z} (b):= \Phi (z\circ \Phi^{-1} (b))$. Clearly $z\circ \Phi^{-1} (b)$ operator commutes with $\Phi^{-1}(b)$, and thus, $\Phi (z\circ \Phi^{-1} (b))$ operator commutes with $b$, that is, $[\Phi_{z} (b),\mathfrak{J},b]$ $=$ $[\Phi (z\circ \Phi^{-1} (b)),\mathfrak{J},b] =0 $, for each $b$, which shows that $\Phi_{z}$ is associating. \Cref{t linear associating maps on JBW*-algebras without type I1} implies the existence of a unique element $\alpha(z)\in Z(\mathfrak{\mathfrak{J}})$ and a unique linear mapping $\mu_z : \mathfrak{J}\to Z(\mathfrak{J})$ (both depending on $z$) satisfying 
\begin{equation}\label{eq definition of alpha 0608} 
	\Phi (z\circ \Phi^{-1} (b))=\Phi_{z}(b) = \alpha(z) \circ b + \mu_{z}(b), \hbox{ for all } b\in \mathfrak{J}.
\end{equation} By defining $\alpha: Z(\mathfrak{M})\to Z(\mathfrak{J})$, $z\mapsto \alpha (z)$, we deduce from the bijectivity of $\Phi$ that  
	$$\Phi (z\circ a) = \alpha(z) \circ \Phi(a) + \mu_{z}\Phi (a), \hbox{ and } \ \Phi (z\circ a) - \alpha(z) \circ \Phi(a)\in Z(\mathfrak{J})$$ for all $a\in \mathfrak{M}.$\smallskip
	
By replacing $\Phi$ with $\Phi^{-1}$ in the above arguments, we deduce the existence of a mapping $\widetilde{\alpha}: Z(\mathfrak{J})\to Z(\mathfrak{M})$ with the property that for each $\widetilde{z}\in Z(\mathfrak{J})$, its image,  $\widetilde{\alpha} (\widetilde{z}),$ is the unique  element satisfying 
\begin{equation}\label{eq delf of alpha tilde} \Phi^{-1} (\widetilde{z}\circ \Phi(a))= \widetilde{\alpha} (\widetilde{z}) \circ a + {\mu}_{\widetilde{z}}(a), \hbox{ for all } a\in \mathfrak{M},
\end{equation}
 where ${\mu}_{\widetilde{z}}$ is a linear mapping from $\mathfrak{M}$ to $Z(\mathfrak{M})$. \smallskip
 
By combining \eqref{eq definition of alpha 0608} and \eqref{eq delf of alpha tilde} we arrive to $$ \widetilde{z}\circ \Phi(a) = \Phi \Phi^{-1} (\widetilde{z}\circ \Phi(a)) =  \alpha(\widetilde{\alpha} (\widetilde{z}))\circ \Phi(a) + \Phi {\mu}_{\widetilde{z}}(a), $$ for all $a\in \mathfrak{M}$, which assures that $$\left( \widetilde{z} - \alpha(\widetilde{\alpha} (\widetilde{z}))  \right) \circ \Phi(a) = \Phi {\mu}_{\widetilde{z}}(a) \in Z(\mathfrak{J}), \hbox{ for each } a\in \mathfrak{M}.$$ The surjectivity of $\Phi$ and \Cref{l Jordan version of Brsar Meiers central} lead to conclusion that $  \widetilde{z} = \alpha(\widetilde{\alpha} (\widetilde{z})),$ for all $\widetilde{z}\in Z(\mathfrak{J})$. We similarly obtain $ {z} = \widetilde{\alpha}({\alpha} ({z})),$ for each ${z}\in Z(\mathfrak{M}).$ This shows that $\alpha$ is a bijection and $\widetilde{\alpha} = \alpha^{-1}$. \smallskip

In order to prove that $\alpha$ is an associative homomorphism, we apply the defining properties of $\alpha$ to get 
$$\begin{aligned}
\left(\alpha(t z_1 + s z_2)- t\alpha(z_1) - s \alpha( z_2) \right)\circ \Phi (a) & = \Phi ((t z_1 + s z_2)\circ a) -\mu_{t z_1 + s z_2}\Phi (a)  \\
&- t \Phi (z_1 \circ a) - t \mu_{z_1}\Phi (a) \\
&- s \Phi (z_2\circ a) -s \mu_{ z_2} \Phi (a),
\end{aligned}$$ for all $a\in \mathfrak{M}$ $z_1,z_2\in Z(\mathfrak{M})$, $s,t\in \mathbb{C}$. This implies that $$\begin{aligned}
\Big( \alpha(t z_1 + s z_2)- t\alpha(z_1) - & s \alpha( z_2) \Big)\circ \mathfrak{J} 
\subseteq Z(\mathfrak{J}),
\end{aligned}$$ and \Cref{l Jordan version of Brsar Meiers central} gives $\alpha(t z_1 + s z_2)=  t\alpha(z_1) + s \alpha( z_2)$.\smallskip

Finally, by the previous properties, 
$$\begin{aligned}
\Big( \alpha(z_1) \circ \alpha(z_2) - \alpha (z_1\circ z_2) )\Big) \circ \Phi (a) &=  \alpha(z_1) \circ \Big(\Phi (z_2\circ a )- \mu_{z_2} \Phi(a)  \Big) \\
&- \Phi ((z_1\circ z_2) \circ a )- \mu_{z_1\circ z_2} \Phi(a) \\
&= \Phi (z_1\circ(z_2\circ a) ) -\mu_{z_1} (\Phi (z_2\circ a )) \\
&- \Phi ((z_1\circ z_2) \circ a )- \mu_{z_1\circ z_2} \Phi(a),
\end{aligned}$$  for all $a\in \mathfrak{M}$ $z_1,z_2\in Z(\mathfrak{M}),$ proving that $\Big( \alpha(z_1) \circ \alpha(z_2) - \alpha (z_1\circ z_2) )\Big) \circ \mathfrak{J} \subseteq Z(\mathfrak{J}),$ and \Cref{l Jordan version of Brsar Meiers central} concludes that $ \alpha (z_1\circ z_2) = \alpha(z_1) \circ \alpha(z_2).$
\end{proof}

When in the proof of \Cref{l second technical lemma in last section}, \Cref{t linear associating maps on JBW*-algebras without type I1} and \Cref{l Jordan version of Brsar Meiers central} are replaced with \Cref{t linear associating maps on JBW-algebras without type I1} and \Cref{l Jordan version of Brsar Meiers central for JB-algebras}, respectively, the arguments remain valid to get next result. 

\begin{lemma}\label{l second technical lemma in last section JB-algebras} Let $\mathfrak{M}$ and $\mathfrak{J}$ be JBW$^*$-algebras with no central summands of type $I_1$. Suppose $\Phi: \mathfrak{M}_{sa} \rightarrow \mathfrak{J}_{sa}$ is a linear bijection satisfying $[x,\mathfrak{M}_{sa},y] = 0$  if, and only if, $[\Phi(x),\mathfrak{J}_{sa},\Phi(y)] = 0,$ for all $x,y\in \mathfrak{M}_{sa}$. Then $\Phi$ maps $Z(\mathfrak{M}_{sa})$ onto $Z(\mathfrak{J}_{sa}),$ and there exists an associative isomorphism $\alpha: Z(\mathfrak{M}_{sa}) \to Z(\mathfrak{J}_{sa})$ satisfying $$\Phi (z \circ a) - \alpha (z) \circ \Phi (a)\in Z(\mathfrak{J}_{sa}), \hbox{ and } \Phi^{-1} (\tilde{z} \circ b) - \alpha^{-1} (\tilde{z}) \circ \Phi^{-1} (b)\in Z(\mathfrak{M}_{sa}), $$ for all $z\in Z(\mathfrak{M}_{sa})$, $a\in \mathfrak{M}_{sa},$ $\tilde{z}\in Z(\mathfrak{J}_{sa})$, and $b\in \mathfrak{J}_{sa}$.
\end{lemma}

The next lemma is built upon arguments and ideas from \cite{BreEreVil}.

\begin{lemma}\label{l based on Bresar Eremita Villena} Let $\mathfrak{M}$ be a JB$^*$-algebra, and let $\mathfrak{J}$ be a JBW$^*$-algebra with no central summands of type $I_1$ or $I_2$. Suppose $T: \mathfrak{M}\to \mathfrak{J}$ is a surjective linear mapping, and $\varepsilon : \mathfrak{M}\times \mathfrak{M}\to \mathfrak{J}$ is a bilinear mapping satisfying \begin{equation}\label{eq lemma BEV 1504}
 \varepsilon (a,b)\circ \Big[[T(a)^2,c,d],e,[T(a),c,d]\Big] =0,		
\end{equation} for all $a,b\in \mathfrak{M}$, $c,d,e\in \mathfrak{J}$. Then $\varepsilon \equiv 0$. 
\end{lemma}

\begin{proof} Let $\pi: \mathfrak{J}\to \mathfrak{J}_{\pi}$ be a Jordan factor representation.  Since $\mathfrak{J}$ does not contain summands of type $I_1$ or $I_2$, \Cref{cor:3_lib} assures that $\mathfrak{J}_{\pi}$ is not a spin factor (cf. \Cref{r quadratic}). The mappings $\overline{\varepsilon} = \pi \varepsilon$ and $\overline{T} = \pi T$ satisfy the same identity in \eqref{eq lemma BEV 1504} for all $a,b\in \mathfrak{M}$, $e,c,d\in \mathfrak{J}$. Since $Z(\mathfrak{J}_{\pi}) = \mathbb{C} \unit$, it follows that for all $a,b\in \mathfrak{M}$, $c,d,e\in \mathfrak{J}$, we have \begin{equation}\label{eq ideas in the first paragraph} \overline{\varepsilon} (a,b)= 0, \hbox{ or } \Big[[\overline{T}(a)^2,\pi(c),\pi(d)],\pi(e),[\overline{T}(a),\pi(c),\pi(d)]\Big] =0.
	\end{equation}

If $\Big[[\overline{T}(a)^2,\pi(c),\pi(d)],\pi(e),[\overline{T}(a),\pi(c),\pi(d)]\Big] =0,$ for all $a\in \mathfrak{M}$, $e,$ $c,$ $d\in\mathfrak{J}$, we deduce from the surjectivity of $T$, the strong$^*$-density of the closed unit ball of $\pi(\mathfrak{J})$ in the closed unit ball of $\mathfrak{J}_{\pi}$, and the joint strong$^*$-continuity of the Jordan product of $\mathfrak{J}_{\pi}$ on bounded sets (cf \cite[Theorem]{RodPa91} or \cite[Theorem 9]{PeRo2001}), that $$\Big[[\overline{a}^2,\overline{c},\overline{d}],\overline{e},[\overline{a},\overline{c},\overline{d}]\Big] =0, \hbox{ for every } \overline{a},\overline{c},\overline{d},\overline{e}\in \mathfrak{J}_{\pi}.$$ Since $\mathfrak{J}_{\pi}$ is a prime non-degenerate Jordan algebra, \cite[Lemma 5.2]{BreEreVil} implies that $\mathfrak{J}_{\pi}$ is a spin factor, which is impossible. Therefore, there exist $a_1$ in $\mathfrak{M}$, $c_1,$ $d_1,$ and $e_1$ in $\mathfrak{J}$ such that $\Big[[\overline{T}(a_1)^2,\pi(c_1),\pi(d_1)],\pi(e_1),[\overline{T}(a_1),\pi(c_1),\pi(d_1)]\Big] \neq 0.$ \smallskip

Suppose that $\overline{\varepsilon}\neq 0$. Then there exist $a_2,b_2\in\mathfrak{J}$ with $\overline{\varepsilon} (a_2,b_2) \neq 0$. Consider the mappings $S: \mathfrak{M}^3\to  \mathfrak{J}_{\pi}$,  $F:  \mathfrak{M}\to  \mathfrak{J}_{\pi} $ given by 
$$S(x_1,x_2,x_3) :=  	\Big[[\overline{T}(x_1)\circ \overline{T}(x_2) ,\pi(c_1),\pi(d_1)],\pi(e_1),[\overline{T}(x_3),\pi(c_1),\pi(d_1)]\Big],$$ $$F(x_1) := \overline{\varepsilon} (x_1,b_2).$$ It follows from the hypotheses (see \eqref{eq ideas in the first paragraph}) that, for each $a\in \mathfrak{M}$, $F(a)=0$ or $S(a,a,a) =0$.  Lemma 3.7 in \cite{BreEreVil} assures that $F(a)=0$ for all $a\in \mathfrak{M},$ or $S(a,a,a) =0$ for all $a\in \mathfrak{M}$, but both conclusions are impossible by what we have just seen ($S(a_1,a_1,a_1), F(a_2)\neq 0$). \smallskip

We have therefore shown that $\overline{\varepsilon} = \pi {\varepsilon}= 0$ for every Jordan factor representation $\pi : \mathfrak{J}\to \mathfrak{J}_{\pi}$. Having in mind that $\mathfrak{J}$ admits a faithful family of Jordan factor representations (cf. \cite[Corollary 5.7]{AlfShulStor79GelfandNeumark}), we deduce that $\varepsilon =0$ as desired. 
\end{proof}

We are now in a position to establish our main conclusion on linear bijections preserving operator commutativity. 

\begin{theorem}\label{t linear bijections preserving oper commut} Let $\mathfrak{M}$ and $\mathfrak{J}$ be JBW$^*$-algebras with no central summands of type $I_1$ and $I_2$. Suppose that $\Phi: \mathfrak{M} \rightarrow \mathfrak{J}$ is a linear  bijection preserving operator commutativity in both directions, that is, $$[x,\mathfrak{M},y] = 0  \Leftrightarrow [\Phi(x),\mathfrak{J},\Phi(y)] = 0,$$ for all $x,y\in \mathfrak{M}$. Then there exist an invertible element $z_0$ in $Z(\mathfrak{J})$, a Jordan isomorphism  $J: \mathfrak{M} \rightarrow \mathfrak{J}$, and a linear mapping  $\beta: \mathfrak{M}\to Z(\mathfrak{J})$ satisfying  
    $$ \Phi(x) = z_0 \circ J(x) + \beta(x), $$ for all $x\in \mathfrak{M}.$ Furthermore, this decomposition of $\Phi$ is unique.
\end{theorem}

\begin{proof} The proof is divided into five main steps.\smallskip

\noindent $(1)$ \textit{First step:} By hypotheses $\Phi$ and $\Phi^{-1}$ preserve operator commutativity in both directions, in particular,  $\Phi^{-1}(y)^2$ and $\Phi^{-1}(y)$ operator commute for every $y \in \mathfrak{J}.$  We define a bilinear mapping $B : \mathfrak{J}\times \mathfrak{J} \rightarrow \mathfrak{J}$ given by \begin{equation}\label{eq def of B(y,y)}
	B(x,y) = \Phi \Big( \Phi^{-1}(x)\circ\Phi^{-1}(y)\Big),
\end{equation}
for every $x, y \in \mathfrak{J}$. The mapping $B$ is well defined and symmetric by definition. Furthermore, by the assumptions on $\Phi$ we have  
        $$0= [\Phi^{-1}(y), \mathfrak{M}, \Phi^{-1}(y)^2]=0 \Leftrightarrow 0= [y, \mathfrak{J}, \Phi(\Phi^{-1}(y)^2)] = [y, \mathfrak{J}, B(y,y)],$$ for every $y \in \mathfrak{J},$ that is, $B$ is a symmetric bilinear mapping whose trace is associating. By applying \Cref{teo bilinear mapping JBW* algebra} we derive the existence of an element $\lambda $ in $Z(\mathfrak{J})$, a linear mapping $\mu_1: \mathfrak{J} \rightarrow Z(\mathfrak{J})$, and a symmetric bilinear mapping $\nu_1: \mathfrak{J}\times \mathfrak{J} \rightarrow Z(\mathfrak{J})$ satisfying 
        \begin{equation}\label{eq first application of the quadratic associating map 0908}
            \Phi(\Phi^{-1}(y)^2) = B(y,y) = \lambda \circ y^2 + \mu_1(y)\circ y + \nu_1(y,y),
        \end{equation} for every $y \in \mathfrak{J}$. Let $\alpha: Z(\mathfrak{M})\to Z(\mathfrak{J})$ denote the associative homomorphism given by \Cref{l second technical lemma in last section}. Having in mind the properties of this mapping $\alpha$ and the identity \eqref{eq first application of the quadratic associating map 0908} we arrive to  $$\left\{\begin{aligned}
        \Phi^{-1}(y)^2 &= \alpha^{-1}(\lambda) \circ \Phi^{-1} (y^2) + \left( \Phi^{-1} (\lambda y^2) - \alpha^{-1}(\lambda) \circ \Phi^{-1} (y^2) \right) \\    
        + \alpha^{-1} &(\mu_1(y)) \circ \Phi^{-1}(y) + \left(\Phi^{-1} (\mu_1(y)\circ y)- \alpha^{-1} (\mu_1(y)) \circ \Phi^{-1}(y) \right) \\
        + \Phi^{-1} & \nu_1 (y,y)
        \end{aligned}\right.$$ for all $y\in \mathfrak{M}$, and equivalently, 
        $$
        \left\{\begin{aligned}
        	x^2 &= \alpha^{-1}(\lambda) \circ \Phi^{-1} (\Phi(x)^2) + \left( \Phi^{-1} (\lambda \Phi(x)^2) - \alpha^{-1}(\lambda) \circ \Phi^{-1} (\Phi(x)^2) \right) \\    
        	&+ \alpha^{-1} (\mu_1(\Phi(x))) \circ x + \left(\Phi^{-1} (\mu_1(\Phi(x))\circ \Phi(x))- \alpha^{-1} (\mu_1(\Phi(x))) \circ x \right) \\
        	&+ \Phi^{-1}  \nu_1 (\Phi(x),\Phi(x)),
        \end{aligned}\right.
        $$ for every $x \in \mathfrak{M}$. Let us denote by $\gamma$ the mapping from $ \mathfrak{J}$ into $Z(\mathfrak{J})$ given by $\gamma (x):= \Phi^{-1} (\lambda \Phi(x)^2) - \alpha^{-1}(\lambda) \circ \Phi^{-1} (\Phi(x)^2)  + \Phi^{-1} (\mu_1(\Phi(x))\circ \Phi(x))- \alpha^{-1} (\mu_1(\Phi(x))) \circ x  + \Phi^{-1}  \nu_1 (\Phi(x),\Phi(x))$. According to this notation we can write 
\begin{equation}\label{eq 3007 1}
	\begin{aligned}
		x^2 &= \alpha^{-1}(\lambda) \circ \Phi^{-1} (\Phi(x)^2) +  \alpha^{-1} (\mu_1(\Phi(x))) \circ x + \gamma (x),
	\end{aligned}
\end{equation}  for all $x \in \mathfrak{M}$.    \\


\noindent$(2)$ \textit{Second step}: For each Jordan factor representation $\pi: \mathfrak{M}\to \mathfrak{M}_{\pi}$, the element $\pi (\alpha^{-1}(\lambda))$ is invertible in $Z(\mathfrak{M}_{\pi})$. Assume otherwise, that $\pi (\alpha^{-1}(\lambda)) =0$ for certain Jordan factor representation $\pi: \mathfrak{M}\to \mathfrak{M}_{\pi}$.  Having in mind the assumptions on $\mathfrak{M},$ we can find an element $w \in \mathfrak{M}$ and elementary operators $\mathcal{E}_0, \mathcal{E}_1, \mathcal{E}_2$ on $\mathfrak{M},$ and elementary operators $\widehat{\mathcal{E}}_j \in \element(\mathfrak{M}_{\pi})$ satisfying $\mathcal{E}_i(w^j) = \delta_{ij} \unit,$ and $\widehat{\mathcal{E}}_j (\pi (w)^{i}) = \delta_{i,j} \unit,$ for all $ i,j \in \{0,1,2\}$ (cf. \Cref{cor:3_lib} and 
\Cref{lemma:elem_homo}(ii)). If in \eqref{eq 3007 1} we replace $x$ with $w$ and we apply $\pi$ and $\widehat{\mathcal{E}}_2$ consecutively, we arrive to 
$$\begin{aligned}
	\unit = \widehat{\mathcal{E}}_2 \pi (w^2) &= \pi(\alpha^{-1}(\lambda)) \circ \widehat{\mathcal{E}}_2(\pi(\Phi^{-1} (\Phi(w)^2))) +  \pi(\alpha^{-1} (\mu_1(\Phi(x)))) \circ \widehat{\mathcal{E}}_2(\pi(w)) \\
	&+ \pi(\gamma (x))\circ \widehat{\mathcal{E}}_2( \unit) =0,
\end{aligned}$$ which is impossible.\\

\noindent$(3)$ \textit{Third step}: We claim that $\alpha^{-1}(\lambda)$ (equivalently, $\lambda$) is invertible. To prove the statement, let us begin by picking a faithful family of Jordan factor representations $\{\pi_i : \mathfrak{M} \rightarrow \mathfrak{M}_{\pi_i}\}_i$ (cf. \cite[Corollary 5.7]{AlfShulStor79GelfandNeumark}). The mapping $\displaystyle \pi_0: \mathfrak{M} \rightarrow \stackrel{\infty}{\bigoplus}_{i} \mathfrak{M}_{\pi_i},$ $\pi_0(a) = (\pi_i(a))_i$ is an isometric unital Jordan $^*$-monomorphism. By $(2)$, $\pi_i (\alpha^{-1}(\lambda))$ is an invertible central element in $\mathfrak{M}_{\pi_i}$ for every $i$. Thus, it suffices to proves that the set $\{\| \pi_i (\alpha^{-1}(\lambda))\| : i\}$ is bounded (from below). Arguing by contradiction, we assume the existence of a sequence $(\pi_{i_n} (\alpha^{-1}(\lambda)))_{i_n}$ converging to zero in norm. \smallskip

For each natural $n$ we can find elementary operators $\widehat{\mathcal{E}}^n_j \in \element(\mathfrak{M}_{\pi_{i_n}})$ satisfying  $\widehat{\mathcal{E}}^n_j \pi_{i_n} =\pi_{i_n} {\mathcal{E}}_j,$ and $\|\widehat{\mathcal{E}}^n_j\|\leq 10,$ for all $j \in \{0,1,2\},$ $n\in \mathbb{N}$ (see 
\Cref{lemma:elem_homo metric with norm}), where $w\in \mathfrak{M}_{sa}$ and the elementary operators ${\mathcal{E}}_j$ are those considered in the previous step.\smallskip 

Now, we deduce from \eqref{eq 3007 1} that 
\begin{equation}\label{eq 8.5 on 0812}\left\{ \begin{aligned}
		\pi_{i_n} (x^2) &= \pi_{i_n}(\alpha^{-1}(\lambda)) \circ \pi_{i_n}(\Phi^{-1} (\Phi(x)^2)) +  \pi_{i_n}(\alpha^{-1} (\mu_1(\Phi(x)))) \circ \pi_{i_n} (x) \\
		&+ \pi_{i_n}(\gamma (x)),
	\end{aligned}\right.
\end{equation} for all $x\in \mathfrak{M}$, $n\in \mathbb{N}$. \smallskip

Let $\mathcal{U}$ be a free ultrafilter over $\mathbb{N}$, and consider the ultraproduct, $( \mathfrak{M}_{\pi_{i_n}})_{\mathcal{U}},$ of the family $\{\mathfrak{M}_{\pi_{i_n}}\}_{n\in \mathbb{N}}$. Elements in $(\mathfrak{M}_{\pi_n})_{\mathcal{U}}$ are written in the form $[x_n]_{\mathcal{U}},$ and the sequence $(x_n)_n$ is called
a \emph{representing family} or a \emph{representative} of $[x_n]_{\mathcal{U}}$. It is known that  $\|[x_n]_{\mathcal{U}}\|=\lim_{\mathcal{U}}\|x_n\|$ independently of the chosen representative \cite{Hein80}. The Banach space $( \mathfrak{M}_{\pi_{i_n}})_{\mathcal{U}}$ is a JB$^*$-algebra with product and involution defined by 
$$ [x_n]_{\mathcal{U}}\circ [y_n]_{\mathcal{U}} = [x_n\circ y_n]_{\mathcal{U}}, \hbox{ and } [x_n]_{\mathcal{U}}^* = [x_n^*]_{\mathcal{U}},$$ respectively (the argument in \cite[Proposition 3. 1]{Hein80} works here, or see \cite[Corollary 10]{Dineen86} for a more general statement). The mappings $$\pi_{\mathcal{U}}: \mathfrak{M}\to ( \mathfrak{M}_{\pi_{i_n}})_{\mathcal{U}}, \ \pi_{\mathcal{U}} (x) := [\pi_{i_n} (x)]_{\mathcal{U}},$$ and  $$\mathcal{E}_j^{\mathcal{U}}: ( \mathfrak{M}_{\pi_{i_n}})_{\mathcal{U}} \to ( \mathfrak{M}_{\pi_{i_n}})_{\mathcal{U}}, \ \mathcal{E}_j^{\mathcal{U}} ([x_n]_{\mathcal{U}}) := [\widehat{\mathcal{E}}^n_j(x_n)]_{\mathcal{U}},$$ are well-defined bounded linear operators (cf. \cite[Definition 2.2]{Hein80}), and $\pi_{\mathcal{U}}$ is a Jordan $^*$-homomorphism. By \eqref{eq 8.5 on 0812} and the assumption made on the sequence $(\pi_{i_n} (\alpha^{-1}(\lambda)))_n$, we derive that 
$$\begin{aligned}
\pi_{\mathcal{U}} (x)^2 =\pi_{\mathcal{U}} (x^2) &=[\pi_{i_n}(\alpha^{-1}(\lambda))]_{\mathcal{U}} \circ [\pi_{i_n}(\Phi^{-1} (\Phi(x)^2))]_{\mathcal{U}} \\
&+   [\pi_{i_n}(\alpha^{-1} (\mu_1(\Phi(x))))]_{\mathcal{U}} \circ [\pi_{i_n} (x)]_{\mathcal{U}} + [\pi_{i_n}(\gamma (x))]_{\mathcal{U}} \\
&= [\pi_{i_n}(\alpha^{-1} (\mu_1(\Phi(x))))]_{\mathcal{U}} \circ [\pi_{i_n} (x)]_{\mathcal{U}} + [\pi_{i_n}(\gamma (x))]_{\mathcal{U}},
\end{aligned}$$ for all $x\in \mathfrak{M}$, where $[\pi_{i_n}(\alpha^{-1} (\mu_1(\Phi(x))))]_{\mathcal{U}}$ and $[\pi_{i_n}(\gamma (x))]_{\mathcal{U}}$ lie in $Z(( \mathfrak{M}_{\pi_{i_n}})_{\mathcal{U}})$. By replacing $x$ with $w$ and applying $\mathcal{E}_2^{\mathcal{U}}$ we get $$\begin{aligned} \unit_{( \mathfrak{M}_{\pi_{i_n}})_{\mathcal{U}}} &= [\pi_{i_n}(\unit)]_{\mathcal{U}}= [\pi_{i_n} {\mathcal{E}}_2 (w^2)]_{\mathcal{U}} = [\widehat{\mathcal{E}}^n_j(\pi_{i_n} (w^2))]_{\mathcal{U}} = 
\mathcal{E}_2^{\mathcal{U}} (\pi_{\mathcal{U}} (w^2)) \\
&=  [\pi_{i_n}(\alpha^{-1} (\mu_1(\Phi(w))))]_{\mathcal{U}} \circ \mathcal{E}_2^{\mathcal{U}}([\pi_{i_n} (w)]_{\mathcal{U}}) + [\pi_{i_n}(\gamma (w))]_{\mathcal{U}} \circ \mathcal{E}_2^{\mathcal{U}} ([\unit]_{\mathcal{U}}) \\
&=  [\pi_{i_n}(\alpha^{-1} (\mu_1(\Phi(w))))]_{\mathcal{U}} \circ [\widehat{\mathcal{E}}^n_2(\pi_{i_n} (w))]_{\mathcal{U}} + [\pi_{i_n}(\gamma (w))]_{\mathcal{U}} \circ [\widehat{\mathcal{E}}^n_2 \pi_{i_n} (\unit)]_{\mathcal{U}}\\
&=  [\pi_{i_n}(\alpha^{-1} (\mu_1(\Phi(w))))]_{\mathcal{U}} \circ [\pi_{i_n} (\mathcal{E}_2(w))]_{\mathcal{U}} + [\pi_{i_n}(\gamma (w))]_{\mathcal{U}} \circ [\pi_{i_n} (\mathcal{E}_2 (\unit))]_{\mathcal{U}}\\
&= 0, 
\end{aligned} $$ which is impossible.\smallskip

\noindent$(4)$ \textit{Fourth step}: There exist an invertible element $z_0$ in $Z(\mathfrak{J})$, a Jordan isomorphism  $J: \mathfrak{M} \rightarrow \mathfrak{J}$, and a linear mapping  $\beta: \mathfrak{M}\to Z(\mathfrak{J})$ satisfying  
$$ \Phi(x) = z_0 \circ J(x) + \beta(x), $$ for all $x\in \mathfrak{M}.$\smallskip

For this final step we follow the construction in \cite{BreEreVil}. We begin by observing that, from \eqref{eq first application of the quadratic associating map 0908} and the bijectivity of $\Phi$,  we can write $$ \Phi(x^2) = \lambda \circ \Phi(x)^2 + \mu_1(\Phi(x))\circ \Phi(x) + \nu_1(\Phi(x),\Phi(x)),$$ for every $x \in \mathfrak{M}$. Set $J: \mathfrak{M}\to \mathfrak{J}$, $J(x) =  \lambda \circ \Phi(x) + \frac12 \mu_1(\Phi(x))$. It is not hard to see from the formula displayed two lines above, and the properties of central elements, that $$\begin{aligned}
J(x^2)- J(x)^2 &= \lambda \circ \Phi(x^2) + \frac12 \mu_1(\Phi(x^2)) - \left(  \lambda \circ \Phi(x) +\frac12  \mu_1(\Phi(x)) \right)^2 \\
=& \lambda \circ \left( \lambda \circ \Phi(x)^2 + \mu_1(\Phi(x))\circ \Phi(x) + \nu_1(\Phi(x),\Phi(x)) \right)+\frac12 \mu_1(\Phi(x^2))  \\
&-  \lambda^2 \circ \Phi(x)^2 - \frac14 \mu_1(\Phi(x))^2 - \left(\lambda \circ \mu_1 (\Phi(x))\right) \circ \Phi(x) \\
=&    \lambda \circ  \nu_1(\Phi(x),\Phi(x)) + \frac12 \mu_1(\Phi(x^2)) - \frac14 \mu_1(\Phi(x))^2\in Z(\mathfrak{J}), 
\end{aligned}$$ for all $x\in \mathfrak{M}$, that is, $$ J(x)^2- J(x^2) =\varepsilon (x,x)\in Z(\mathfrak{J}),$$ is a symmetric bilinear mapping taking values in the centre of $\mathfrak{J}$.  It follows that $$\begin{aligned}
J(x\circ y )&=\frac12 \Big(J((x+y)^2)- J(x^2)- J(y^2) \Big) \\
&= \frac12 \Big(J(x+y)^2-\varepsilon (x+y,x+y) - J(x)^2- \varepsilon (x,x) - J(y)^2 -\varepsilon (y,y) \Big)\\
&= J(x)\circ J(y) - \varepsilon (x,y), \ \ x,y\in \mathfrak{M}.
\end{aligned}$$ Now, the identities $$\begin{aligned}
J((x^2\circ y)\circ x) &= J(x^2\circ y )\circ J(x) -  \varepsilon(x^2\circ y,x) \\
&= \left(J(x^2)\circ J(y)\right)\circ J(x) - \varepsilon(x^2,y)\circ J(x) -  \varepsilon(x^2\circ y,x) \\
&= \left(J(x)^2\circ J(y)\right)\circ J(x)- \left(\varepsilon(x,x)\circ J(y)\right)\circ J(x) - \varepsilon(x^2,y)\circ J(x) \\
&-  \varepsilon(x^2\circ y,x),
\end{aligned}$$ and 
$$\begin{aligned}
	J((x\circ y)\circ x^2) &= J(x\circ y )\circ J(x^2) -  \varepsilon(x\circ y,x^2) \\
	&= \left(J(x)\circ J(y)\right)\circ J(x^2) - \varepsilon(x,y)\circ J(x^2) -  \varepsilon(x\circ y,x^2) \\
	&= \left(J(x)\circ J(y)\right)\circ J(x)^2 - \left(J(x)\circ J(y)\right)\circ \varepsilon (x,x)  - \varepsilon(x,y)\circ J(x^2) \\
	&-  \varepsilon(x\circ y,x^2),
\end{aligned}$$ combined with the Jordan identity ($(x\circ y)\circ x^2= (x^2\circ y)\circ x$), assure that $$\begin{aligned}  &(\varepsilon(x,y)\circ\lambda) \circ \Big( \lambda \circ \Phi(x)^2 +\mu_1(\Phi(x))\circ \Phi(x)  + \nu_1(\Phi(x),\Phi(x)) \Big)\\ &+\frac12  \varepsilon(x,y)\circ \mu_1(\Phi(x^2)) =\varepsilon(x,y)\circ \left(\lambda \circ \Phi(x^2) + \frac12 \mu_1(\Phi(x^2))\right) \\ & =  \varepsilon(x,y)\circ J(x^2)= \varepsilon(x^2,y)\circ J(x) + \varepsilon(x^2\circ y,x) - \varepsilon(x\circ y,x^2)\\
&= \varepsilon(x^2,y)\circ \left( \lambda \circ \Phi(x) + \frac12 \mu_1(\Phi(x)) \right) + \varepsilon(x^2\circ y,x) - \varepsilon(x\circ y,x^2),
\end{aligned}$$ and thus, the invertibility of $\lambda$ implies that $$ \varepsilon(x,y)\circ \Phi(x)^2 \in Z(\mathfrak{J}) \circ J(x) +  Z(\mathfrak{J}), \hbox{ for all } x,y\in \mathfrak{M}.$$ We consequently have 
$$0= \Big[ [\varepsilon(x,y)\circ \Phi(x)^2, b,c],\mathfrak{J},[\Phi(x), b,c] \Big] =\varepsilon(x,y)\circ  \Big[ [ \Phi(x)^2, b,c],\mathfrak{J},[\Phi(x), b,c] \Big],$$ for all $x,y\in\mathfrak{M}$, $b,c\in \mathfrak{J}$. We are in a position to apply \Cref{l based on Bresar Eremita Villena}, which assures that $\varepsilon(x,y)=0$, for all $x,y\in\mathfrak{M}$, and thus $J$ is a Jordan homomorphism. \smallskip

Observe finally that, by the invertibility of the central element $\lambda$ we can write $$ \Phi(x)= \lambda^{-1} \circ J(x)-\frac12 \lambda^{-1} \circ \mu_1(\Phi(x)),\ \  x\in \mathfrak{M}.$$ Taking $z_0 = \lambda^{-1}$ and $\beta (x) = -\frac12 \lambda^{-1} \circ \mu_1(\Phi(x))$ we get the desired factorization for $\Phi$.\smallskip

It only remains to prove that $J$ is bijective. If $a\in \ker(J)$, we have $\Phi (a) = \beta(a)\in Z(\mathfrak{J})$, and thus $a\in Z(\mathfrak{M})$. We have shown that $\ker(J)\subseteq Z(\mathfrak{M})$. Having in mind that $\ker(J)$ is a Jordan ideal of $\mathfrak{M}$, we have $\ker(J)\circ \mathfrak{M} \subseteq \ker(J) \subseteq Z(\mathfrak{M})$. \Cref{l Jordan version of Brsar Meiers central} proves that $\ker(J) =\{0\}$. \smallskip

Since $J(x) = \lambda \circ \Phi(x) + \frac12 \mu_1(\Phi(x))$, and $\Phi$ maps the centre onto the centre, we deduce that $J(Z(\mathfrak{M}))  \subseteq Z(\mathfrak{J})$. Let $\alpha: Z(\mathfrak{M})\to Z(\mathfrak{J})$ denote the associative homomorphism given by \Cref{l second technical lemma in last section} we employed in the first step. It is not hard to see from the identities linking $\alpha$, $\Phi$ and $J$, in particular $\Phi (z x) = \alpha (x) \circ \Phi (x)$ for all $x\in \mathfrak{M}$, $z\in Z(\mathfrak{M})$,  that 
$$z_0 \circ \Big( (\alpha(z)-J(z))\circ J(x) \Big) = \beta(z\circ x)- \alpha(z) \circ \beta (x),$$ and thus  $$(\alpha(z)-J(z))\circ \Phi(x) = \beta(z\circ x)- \alpha(z) \circ \beta (x) + (\alpha(z)-J(z))\circ \beta(x)\in Z(\mathfrak{J}),$$  for all $x\in \mathfrak{M}$, $z\in Z(\mathfrak{M})$. It follows from the surjectivity of $\Phi$ and \Cref{l Jordan version of Brsar Meiers central} that $\alpha(z)=J(z)$ for all $z\in Z(\mathfrak{M})$. If we write $$\begin{aligned}
\Phi (x) &= z_0 \circ J(x) + \beta (x) = J(\alpha^{-1}(z_0)) \circ J(x) + J(\alpha^{-1}(\beta (x))) \\
&=J(\alpha^{-1}(z_0) \circ x + \alpha^{-1}(\beta (x))), \ \ x\in \mathfrak{M},  
\end{aligned}$$ we see that the surjectivity of $\Phi$ induces the same property on $J$.\smallskip

\noindent$(5)$ \textit{Fifth, and last, step}: The Jordan isomorphism $J$ and the mapping $\beta$ are unique. Suppose we can write $\Phi = z_1 \circ J_1 + \beta_1,$ where $z_1$ is an invertible element in $Z(\mathfrak{J})$, $J: \mathfrak{M}\to \mathfrak{J}$ is a Jordan isomorphism, and $\beta: \mathfrak{M}\to Z(\mathfrak{J})$ is a linear mapping. We shall show that $z_0 =z_1$, $J = J_1,$ and $\beta = \beta_1$.\smallskip

Arguing as in the fourth step we see that $\alpha (x) \circ \Phi (x) = \Phi (z x) $ for all $x\in \mathfrak{M}$, $z\in Z(\mathfrak{M})$,  and thus, by the assumptions on $J_1$ and $z_1,$ we derive that $$z_1 \circ \Big( (\alpha(z)-J_1(z))\circ J_1(x) \Big) = \beta_1(z\circ x)- \alpha(z) \circ \beta_1 (x),$$ and thus  $$(\alpha(z)-J_1(z))\circ \Phi(x) = \beta_1(z\circ x)- \alpha(z) \circ \beta_1 (x) + (\alpha(z)-J_1(z))\circ \beta_1(x)\in Z(\mathfrak{J}),$$ for all $x\in \mathfrak{M}$, $z\in Z(\mathfrak{M})$. It follows from the surjectivity of $\Phi$ and \Cref{l Jordan version of Brsar Meiers central} that $\alpha(z)=J_1(z),$ for all $z\in Z(\mathfrak{M})$. Having this property in mind we get $$\begin{aligned}
\Phi (x) &= z_1\circ J_1 (x) + \beta_1 (x) = J_1 (\alpha^{-1}(z_1))\circ J_1 (x) + J_1(\alpha^{-1}(\beta_1 (x))) \\
&= J_1 \Big( \alpha^{-1}(z_1)\circ x + \alpha^{-1} (\beta_1 (x)) \Big)= J_1 \Big( \alpha^{-1}(z_1)\circ x + J_1^{-1} (\beta_1 (x)) \Big),
\end{aligned}$$   for all $x\in \mathfrak{M}$. It is not hard to check from this that  $$\begin{aligned}
\Phi^{-1} (y) &=\alpha^{-1}(z_1^{-1})\circ \Big( J_1^{-1}(y) - J_1^{-1} (\beta_1 (\Phi^{-1}(y))) \Big)\\
&=J_1^{-1}(z_1^{-1})\circ \Big( J_1^{-1}(y) - J_1^{-1} (\beta_1 (\Phi^{-1}(y))) \Big) \\
&=J_1^{-1}\Big(z_1^{-1}\circ \Big( y - \beta_1 (\Phi^{-1}(y)) \Big) \Big),
\end{aligned}$$ for all $y\in \mathfrak{J}$.\smallskip

We turn now our attention to the bilinear mapping $B : \mathfrak{J}\times \mathfrak{J} \rightarrow \mathfrak{J}$ defined in the first step (see \eqref{eq def of B(y,y)}). Having in mind the expression for the mapping $\Phi^{-1}$ in the previous paragraph we arrive to $$\begin{aligned}
 B(y,y) &= \Phi \Big( \Phi^{-1}(y)^2\Big) = \Phi \Big( J_1^{-1}\Big(z_1^{-2}\circ \Big( y - \beta_1 (\Phi^{-1}(y)) \Big)^2 \Big) \Big)\\
 &= z_1\circ \Big(z_1^{-2}\circ \Big( y - \beta_1 (\Phi^{-1}(y)) \Big)^2 \Big) +\beta_1 \Big( J_1^{-1}\Big(z_1^{-2}\circ \Big( y - \beta_1 (\Phi^{-1}(y)) \Big)^2 \Big) \Big)\\
 &= z_1^{-1} \circ y^2 - 2 \left(z_1^{-1}\circ \beta_1 (\Phi^{-1}(y)) \right) \circ y + z_1^{-1} \circ \beta_1 (\Phi^{-1}(y))^2 \\
 &+\beta_1 \Big( J_1^{-1}\Big(z_1^{-2}\circ \Big( y - \beta_1 (\Phi^{-1}(y)) \Big)^2 \Big) \Big),
\end{aligned}$$ for all $y\in \mathfrak{J}$, where $z_1^{-1} \circ \beta_1 (\Phi^{-1}(y))^2 +\beta_1 \Big( J_1^{-1}\Big(z_1^{-2}\circ \Big( y - \beta_1 (\Phi^{-1}(y)) \Big)^2 \Big) \Big)$ lies in $Z(\mathfrak{J})$. The uniqueness of the decomposition of $B(\cdot,\cdot)$ in \eqref{eq first application of the quadratic associating map 0908} guaranteed by \Cref{teo bilinear mapping JBW* algebra} assures that $z_1^{-1} = \lambda = z_0^{-1}$, and $\mu_1 (y) = - 2 \left(z_1^{-1}\circ \beta_1 (\Phi^{-1}(y)) \right)$, for all $y\in \mathfrak{J}$. It then follows that 
$$\begin{aligned}
J(x) &= \lambda\circ \Phi(x) +\frac12 \mu_1(\Phi(x)) = \lambda\circ \left(z_1\circ J_1 (x) + \beta_1 (x)\right) - z_1^{-1}\circ \beta_1 (x) \\
&= J_1 (x) + z_1^{-1}\beta_1 (x) - z_1^{-1}\circ \beta_1 (x) = J_1 (x),
\end{aligned}$$ for all $x\in \mathfrak{M}$. This shows that $z_1= z_0$ and $J_1 = J$, which gives $\beta = \beta_1$ and finishes the proof.
\end{proof}

\begin{remark}\label{r optimality of the hypotheses in the previous theorem} We analyse now the optimality of the hypotheses in \Cref{t linear bijections preserving oper commut} affirming that the involved JBW$^*$-algebras admits no central summands of type $I_1$ and $I_2$. Trivially every linear mapping on a commutative C$^*$-algebra preserves operator commutativity. Let $V$ be a spin factor. We have seen in \Cref{r associator-commutator of an element in a spin} that elements $a,b$ in $V$ operator commute if and only if $b$ is a linear combination of $a$ and $\unit$. Therefore, any linear mapping on a spin factor preserves operator commutativity.  Let $\Phi : V\to V$ be any linear bijection such that $\Phi (\unit)\notin \mathbb{C} \unit$. We have already justified that $\Phi$ preserves operator commutativity in both directions. 
Recall that $Z(V) = \mathbb{C} \unit$. If the mapping $\Phi$ were written in the form   
	$$ \Phi(x) = z_0 \circ J(x) + \beta(x), $$ for all $x\in V,$ where $z_0$ is an invertible element in $Z(V) = \mathbb{C} \unit$, a Jordan isomorphism  $J: V \rightarrow V$, and a linear mapping  $\beta: V\to Z(V)$, it would follow that $\Phi (\unit) = z_0 \unit + \beta (\unit) \in \mathbb{C}\unit$, which is impossible. 
\end{remark}

Let $\mathfrak{M}$ and $\mathfrak{J}$ be two JB$^*$-algebras.  Following the standard notation, for each mapping $F: \mathfrak{M}\to \mathfrak{J}$, we write $F^{\sharp}:\mathfrak{M}\to \mathfrak{J}$ defined by $F^{\sharp} (x) = F(x^*)^*$ ($x\in \mathfrak{M}$). Clearly, $F$ is linear if, and only if, $F^{\sharp}$ is, and $F^{\sharp \sharp} = F$. The mapping $F$ is called \emph{symmetric} or \emph{$\sharp$-symmetric} if $F^{\sharp} = F$. Observe that symmetric Jordan homomorphisms between JB$^*$-algebras are usually called Jordan $^*$-homomorphisms in the literature. Every real linear mapping between the self-adjoint parts of two JBW$^*$-algebras admits a natural extension to a symmetric complex linear mapping between both JBW$^*$-algebras.\smallskip

Our next result will consist in showing that when the mapping $\Phi$ in \Cref{t linear bijections preserving oper commut} is symmetric (i.e., $\Phi^{\sharp} = \Phi$) we can actually assume that $J$ is a Jordan $^*$-isomorphism. The results seems to be new also in the case of JBW$^*$-algebra factors (cf. \cite[\S 5.1]{BreEreVil}), and the arguments are completely different to those in \cite{BreMei1993}, where the result is proved in the case of von Neumann algebras.

\begin{corollary}\label{c linear bijections preserving oper commut symmetric} Let $\mathfrak{M}$ and $\mathfrak{J}$ be JBW$^*$-algebras with no central summands of type $I_1$ and $I_2$. Suppose that $\Phi: \mathfrak{M} \rightarrow \mathfrak{J}$ is a symmetric linear bijection preserving operator commutativity in both directions, that is, $\Phi^{\sharp}=\Phi$ and $$[x,\mathfrak{M},y] = 0  \Leftrightarrow [\Phi(x),\mathfrak{J},\Phi(y)] = 0,$$ for all $x,y\in \mathfrak{M}$. Then there exist a unique invertible element $z_0$ in $Z(\mathfrak{J})_{sa}$, a unique Jordan $^*$-isomorphism  $J: \mathfrak{M} \rightarrow \mathfrak{J}$, and a unique symmetric linear mapping  $\beta: \mathfrak{M}\to Z(\mathfrak{J})$ satisfying  
	$$ \Phi(x) = z_0 \circ J(x) + \beta(x), $$ for all $x\in \mathfrak{M}.$ 
\end{corollary}
  
\begin{proof} By \Cref{t linear bijections preserving oper commut} there exist a unique invertible element $z_0$ in $Z(\mathfrak{J})$, a unique Jordan isomorphism  $J: \mathfrak{M} \rightarrow \mathfrak{J}$, and a unique linear mapping  $\beta: \mathfrak{M}\to Z(\mathfrak{J})$ satisfying  
	$$ \Phi(x) = z_0 \circ J(x) + \beta(x), $$ for all $x\in \mathfrak{M}.$ Since $\Phi^{\sharp} = \Phi$ we get $$ z_0 \circ J(x) + \beta(x) = \Phi(x)= \Phi^{\sharp} (x) = z_0^* \circ J^{\sharp}(x) + \beta^{\sharp}(x),$$ for all $x\in \mathfrak{M}$. The uniqueness of the decomposition in the just quoted \Cref{t linear bijections preserving oper commut} implies that $z_0 = z_0^*$, $J^{\sharp} = J,$ and $\beta^{\sharp} = \beta$, which finishes the proof. 
\end{proof}  

JBW-algebras constitute a suitable mathematical model in quantum mechanics as its elements describe observable physical quantities. Operator commutativity for elements in a JBW-algebra can be seen as a generalization of commutativity (or compatibility in the
language of quantum mechanics) for elements in $B(H)_{sa}$. Therefore, the study of linear bijections between JBW-algebras preserving operator commutativity in both directions is an interesting question, worth to be consider by itself. Unfortunately, the result is not a mere consequence of what we proved for JBW$^*$-algebras and requires some extra arguments.  

\begin{theorem}\label{t linear bijections preserving oper commut JBW-algebras} Let $\mathfrak{M}$ and $\mathfrak{J}$ be JBW$^*$-algebras with no central summands of type $I_1$ and $I_2$. Suppose that $\Phi: \mathfrak{M}_{sa} \rightarrow \mathfrak{J}_{sa}$ is a linear  bijection preserving operator commutativity in both directions, that is, $$[x,\mathfrak{M}_{sa},y] = 0  \Leftrightarrow [\Phi(x),\mathfrak{J}_{sa},\Phi(y)] = 0,$$ for all $x,y\in \mathfrak{M}_{sa}$. Then there exist an invertible element $z_0$ in $Z(\mathfrak{J}_{sa})$, a Jordan isomorphism  $J: \mathfrak{M}_{sa} \rightarrow \mathfrak{J}_{sa}$, and a linear mapping  $\beta: \mathfrak{M}_{sa}\to Z(\mathfrak{J}_{sa})$ satisfying  
	$$ \Phi(x) = z_0 \circ J(x) + \beta(x), $$ for all $x\in \mathfrak{M}_{sa}.$ Furthermore, this decomposition of $\Phi$ is unique.
\end{theorem}

\begin{proof} Let us begin with some simple observations. The natural complex linear extension $\Phi_{\mathbb{C}}: \mathfrak{M} \rightarrow \mathfrak{J},$ $\Phi_{\mathbb{C}} (a + i b ):= \Phi(a) + i \Phi(b)$, $a+ib\in \mathfrak{M}$, is a symmetric linear bijection. However, it is not clear at this stage that $\Phi_{\mathbb{C}}$ preserves operator commutativity. For this reason we must rebuild and modify part of the proof of \Cref{t linear bijections preserving oper commut}.\smallskip
	
\Cref{l second technical lemma in last section JB-algebras} assures that $\Phi$ maps $Z(\mathfrak{M}_{sa})$ onto $Z(\mathfrak{J}_{sa}),$ and the existence an associative isomorphism $\alpha: Z(\mathfrak{M}_{sa}) \to Z(\mathfrak{J}_{sa})$ satisfying $\Phi (z \circ a) - \alpha (z) \circ \Phi (a)\in Z(\mathfrak{J}_{sa}),$ for all $z\in Z(\mathfrak{M}_{sa})$, $a\in \mathfrak{M}_{sa},$ \ and $\Phi^{-1} (\tilde{z} \circ b) - \alpha^{-1} (\tilde{z}) \circ \Phi^{-1} (b)\in Z(\mathfrak{M}_{sa}),$ for all $\tilde{z}\in Z(\mathfrak{J}_{sa})$, $b\in \mathfrak{J}_{sa}$. Since $Z(\mathfrak{J}_{sa}) = Z(\mathfrak{J})_{sa}$ and $Z(\mathfrak{M}_{sa}) = Z(\mathfrak{M})_{sa}$, the mapping $\Phi_{\mathbb{C}}$ maps $Z(\mathfrak{M})$ onto $Z(\mathfrak{J})$. The mapping $\alpha_{\mathbb{C}} : Z(\mathfrak{M})\to Z(\mathfrak{J})$, $\alpha_{\mathbb{C}} (a+ i b) = \alpha (a) + i \alpha (b)$, is an associative isomorphism satisfying 
\begin{equation}\label{eq properties of the mapping alpha complex}\left\{ \begin{aligned}
		\Phi_{\mathbb{C}} (z \circ x) - \alpha_{\mathbb{C}} (z) \circ \Phi_{\mathbb{C}} (x)\in & Z(\mathfrak{J}),\hbox{ and } \\ \Phi_{\mathbb{C}}^{-1} (\tilde{z} \circ y) - \alpha_{\mathbb{C}}^{-1} (\tilde{z}) \circ \Phi_{\mathbb{C}}^{-1} (y)&\in Z(\mathfrak{M}),
	\end{aligned}\right.
\end{equation} for all $z\in Z(\mathfrak{M})$, $x\in \mathfrak{M},$ $\tilde{z}\in Z(\mathfrak{J})$, $y\in \mathfrak{J}$.\smallskip

We claim that 
\begin{equation}\label{eq complex linear extension defines an associating bilinear}
[\Phi_{\mathbb{C}} (x), \mathfrak{J}, \Phi_{\mathbb{C}} (x^2)] =0 = [\Phi_{\mathbb{C}}^{-1} (y), \mathfrak{M}, \Phi_{\mathbb{C}}^{-1} (y^2)],
\end{equation} for all $x\in\mathfrak{M}$, $y\in \mathfrak{J}.$ Namely, given $x = a+ i b\in \mathfrak{M}$ with $a,b\in \mathfrak{M}_{sa}$, the hypothesis on $\Phi$ implies that $$\begin{aligned} 0 =  [\Phi (a\pm b) , \mathfrak{J}_{sa}, \Phi((a\pm b)^2) ] &= [\Phi (a) , \mathfrak{J}_{sa}, \Phi(a^2) ] \pm [\Phi ( b) , \mathfrak{J}_{sa}, \Phi(b^2) ]\\
&\pm 2 [\Phi (a) , \mathfrak{J}_{sa}, \Phi( a\circ  b) ] + 2 [\Phi (b) , \mathfrak{J}_{sa}, \Phi( a\circ  b) ] \\
	&+ [\Phi (a) , \mathfrak{J}_{sa}, \Phi(b^2) ] \pm [\Phi (b) , \mathfrak{J}_{sa}, \Phi(a^2) ] \\
	&=  \pm 2 [\Phi (a) , \mathfrak{J}_{sa}, \Phi( a\circ  b) ] + 2 [\Phi (b) , \mathfrak{J}_{sa}, \Phi( a\circ  b) ] \\
	&+ [\Phi (a) , \mathfrak{J}_{sa}, \Phi(b^2) ] \pm [\Phi (b) , \mathfrak{J}_{sa}, \Phi(a^2)],
\end{aligned} $$ which assures that 
$$ 2 [\Phi (a) , \mathfrak{J}_{sa}, \Phi( a\circ  b) ] + [\Phi (b) , \mathfrak{J}_{sa}, \Phi(a^2)]  = 0,$$ and  $$2 [\Phi (b) , \mathfrak{J}_{sa}, \Phi( a\circ  b) ] + [\Phi (a) , \mathfrak{J}_{sa}, \Phi(b^2) ] = 0.$$ We therefore deduce that 
$$\begin{aligned}
	[\Phi_{\mathbb{C}} (a+ i b) , \mathfrak{J}, \Phi_{\mathbb{C}}((a+i b)^2) ] &= [\Phi (a) , \mathfrak{J}, \Phi(a^2) ] - i [\Phi (b) , \mathfrak{J}, \Phi (b^2) ] \\
	&+ 2 i  [\Phi(a) , \mathfrak{J}, \Phi(a \circ  b) ] - 2  [\Phi ( b) , \mathfrak{J}, \Phi(a\circ  b) ] \\
	&- [ \Phi(a) , \mathfrak{J}, \Phi(b^2) ] - i [\Phi (b) , \mathfrak{J}, \Phi(a^2) ] = 0,
\end{aligned} $$ which proves the first identity in \eqref{eq complex linear extension defines an associating bilinear}. The second identity in the claim follows via similar arguments.\smallskip

It is not hard to check, thanks to \eqref{eq complex linear extension defines an associating bilinear}, that the trace of the symmetric bilinear mapping \begin{equation}\label{eq BC is associating} B: \mathfrak{J}\times\mathfrak{J}\to \mathfrak{J}, \ B(y,b)= \Phi_{\mathbb{C}}\left( \Phi_{\mathbb{C}}^{-1} (y)\circ \Phi_{\mathbb{C}}^{-1} (b)\right) \hbox{ is associating},
\end{equation} that is, $$\begin{aligned}
[B(y,y),\mathfrak{J},y] &= [\Phi_{\mathbb{C}}\left( \Phi_{\mathbb{C}}^{-1} (y)\circ \Phi_{\mathbb{C}}^{-1} (y)\right), \mathfrak{J}, y]\\ 
&= [\Phi_{\mathbb{C}}\left( \Phi_{\mathbb{C}}^{-1} (y)\circ \Phi_{\mathbb{C}}^{-1} (y)\right), \mathfrak{J}, \Phi_{\mathbb{C}}( \Phi_{\mathbb{C}}^{-1}( y))]=0.
\end{aligned}$$ 

Let us find $w\in \mathfrak{J}_{sa}$ and elementary operators $\mathcal{E}_{j}\in \element_{\mathfrak{J}_{sa}} (\mathfrak{J})$ satisfying  $\mathcal{E}_{j} (w^i) = \delta_{ij} \unit,$ for all $i,j\in \{0,1,2\}$ (cf. \Cref{cor:3_lib}). We have gathered all the necessary ingredients in \eqref{eq properties of the mapping alpha complex} and \eqref{eq BC is associating} in order to repeat, line by line, the proof of \Cref{t linear bijections preserving oper commut} but replacing $\Phi$ and $\alpha$ there with  $\Phi_{\mathbb{C}}$ and $\alpha_{\mathbb{C}},$ respectively, to deduce the existence of a unique invertible element $z_0$ in $Z(\mathfrak{J})$, a unique linear mapping $\mu_1: \mathfrak{J} \rightarrow Z(\mathfrak{J})$, and a unique symmetric bilinear mapping $\nu_1: \mathfrak{J}\times \mathfrak{J} \rightarrow Z(\mathfrak{J})$ satisfying \begin{equation}\label{eq B final 0609} \Phi_{\mathbb{C}}(\Phi_{\mathbb{C}}^{-1}(y)^2) = B(y,y) = z_0^{-1} \circ y^2 + \mu_1 (y)\circ y + \nu_1(y,y),
\end{equation}
 for all $y\in \mathfrak{J},$ (cf. \Cref{teo bilinear mapping JBW* algebra}),  $z_0^{-1} = {\mathcal{E}_2}(B(w,w))$ (cf. \Cref{prop: uniq_sol}), and taking $\tilde{J}: \mathfrak{M}\to \mathfrak{J}$, $\tilde{J}(x) =  z_0^{-1} \circ \Phi_{\mathbb{C}}(x) + \frac12 \mu_1(\Phi_{\mathbb{C}}(x))$ it follows that $\tilde{J}$ is a Jordan isomorphism and  
$$ \Phi_{\mathbb{C}}(x) = z_0 \circ \tilde{J}(x) + \tilde{\beta}(x), \hbox{ for all } x\in \mathfrak{M},$$ for a certain linear mapping  $\tilde{\beta}: \mathfrak{M}\to Z(\mathfrak{J}).$ \smallskip

By recalling that all elementary operators  $\mathcal{E}_{j}\in \element (\mathfrak{J})$ can be assumed to be symmetric linear mappings (cf. \Cref{cor:3_lib}) we derive that $$\begin{aligned}
(z_0^{-1})^* = {\mathcal{E}_2}(B(w,w))^* &= {\mathcal{E}_2}(B(w,w)^*)= {\mathcal{E}_2}\left( ( \Phi_{\mathbb{C}}(\Phi_{\mathbb{C}}^{-1}(w)^2))^*\right)\\
& = {\mathcal{E}_2}\left(  \Phi_{\mathbb{C}}(\Phi_{\mathbb{C}}^{-1}(w)^2)\right) = z_0^{-1},
\end{aligned}$$ which proves that $z_0 = z_0^*\in Z(\mathfrak{J}_{sa})$. \smallskip

On the other hand, $\Phi_{\mathbb{C}}$ being a symmetric linear mapping implies that $B(y,y)=B(y^*,y^*)^* $, and thus by \eqref{eq B final 0609}, $$\begin{aligned}
 &z_0^{-1} \circ y^2 + \mu_1 (y)\circ y + \nu_1(y,y) = B(y,y) \\ &= B(y^*,y^*)^* = z_0^{-1} \circ y^2 + \mu_1 (y^*)^*\circ y + \nu_1(y^*,y^*)^*,
\end{aligned}$$ for all $y\in \mathfrak{J}$. The uniqueness of the decomposition in \eqref{eq B final 0609} asserts that $ \mu_1 (y^*)^* =  \mu_1 (y)$ and $\nu_1(y,y) = \nu_1(y^*,y^*)^*$ for all $y\in \mathfrak{J}$. We therefore have $$\begin{aligned}
\tilde{J}(x^*) &=  z_0^{-1} \circ \Phi_{\mathbb{C}}(x^*) + \frac12 \mu_1(\Phi_{\mathbb{C}}(x^*))=  z_0^{-1} \circ \Phi_{\mathbb{C}}(x)^* + \frac12 \mu_1(\Phi_{\mathbb{C}}(x))^*   \\ &= (z_0^{-1} \circ \Phi_{\mathbb{C}}(x) + \frac12 \mu_1(\Phi_{\mathbb{C}}(x)))^* = \tilde{J}(x)^*, \ x\in \mathfrak{M}, 
\end{aligned}$$ and hence $\tilde{J}$ is a Jordan $^*$-isomorphism. It trivially holds that $\tilde{\beta}$ is a symmetric linear mapping. By defining $J = \tilde{J}|_{\mathfrak{J}_{sa}}$ and $\beta = \tilde{\beta}|_{\mathfrak{J}_{sa}}$ we obtain the desired decomposition for $\Phi$. \smallskip

Finally, the uniqueness of the decomposition of $\Phi$ essentially follows from the uniqueness of the decomposition in  \eqref{eq B final 0609} as explicitly shown in the final part of the proof of \Cref{t linear bijections preserving oper commut}.
\end{proof}

We conclude this note with a consequence of our previous \Cref{t linear associating maps on JBW-algebras without type I1} and \Cref{P Topping for non-self adjoint} which is a novelty in the setting of von Neumann algebras. Let us first recall a precedent related to our objective. Let $A$ and $B$ be von Neumann algebras, and suppose additionally that $A$ is a factor not of type $I_1$ or $I_2$.  L. Molnár established in \cite[Theorem 4]{Mol2010} that for every linear bijection $\Phi : A_{sa}\to B_{sa}$ preserving commutativity in both directions, 
there exist a non-zero real number $\lambda$, a Jordan $^*$-isomorphism $J : A \to B$ and a real-linear functional $\beta : A_{sa} \to \mathbb{R}$ such that $$ \Phi (a) = c J(a) + \beta (a) \unit, \hbox{ for all } a \in A_{sa}.$$ We can now relax the hypothesis concerning the von Neumann algebra $A$ and optimize the conclusion in the just quoted result by Molnár.

\begin{corollary}\label{c self-adjoint parts of von Neumann algebras} Let $A$ and $B$ be von Neumann algebras with no central summands of type $I_1$ and $I_2$. Suppose that $\Phi: {A}_{sa} \rightarrow {B}_{sa}$ is a linear  bijection preserving commutativity in both directions. Then there exist an invertible element $z_0$ in $Z({B}_{sa})$, a Jordan isomorphism  $J: {A}_{sa} \rightarrow {B}_{sa}$, and a linear mapping  $\beta: {A}_{sa}\to Z({B}_{sa})$ satisfying  
	$$ \Phi(x) = z_0\ J(x) + \beta(x), \hbox{ for all } x\in {A}_{sa}.$$ Furthermore, this decomposition of $\Phi$ is unique.
\end{corollary} 

\medskip\medskip

\textbf{Acknowledgements}\medskip

Authors supported by MICIU/AEI/10.13039/501100011033 and ERDF/EU grant PID2021-122126NB-C31, by ``Maria de Maeztu'' Excellence Unit IMAG, reference CEX2020-001105-M funded by MICIU/AEI/10.13039/501100011033, and by Junta de Andaluc{\'i}a grants FQM185 and FQM375. First author supported by grant FPU21/00617 at University of Granada founded by Ministerio de Universidades (Spain). Second author partially supported by MOST China project number G2023125007L.

\smallskip\smallskip

\noindent\textbf{Data Availability} Statement Data sharing is not applicable to this article as no datasets were generated or analysed during the preparation of the paper.\smallskip\smallskip

\noindent\textbf{Declarations} 
\smallskip\smallskip

\noindent\textbf{Conflict of interest} The authors declare that they have no conflict of interest.

\end{document}